\documentclass[11pt]{article}

\newcommand{\tn}{\textnormal}
 \newcommand{\lp}{\left(}
\newcommand{\rp}{\right)}

\newcommand{\lb}{\left\lbrace}
\newcommand{\rb}{\right\rbrace} 
	
\newcommand{\mc}{\mathcal}

\usepackage {epsfig,amsmath,euscript}

\usepackage[latin1]{inputenc}

\usepackage[english]{babel}
\usepackage{color}
\usepackage{subcaption}
\usepackage{amsmath}
\usepackage[normalem]{ulem}
\usepackage{amsthm}
\usepackage{amssymb}

\usepackage{ae}
\usepackage{float}
\usepackage[T1]{fontenc}
\usepackage{graphicx}
\usepackage{epstopdf}
\usepackage{longtable}
\usepackage{color}
\definecolor{rred}{rgb}{0.7,0.0,0.2}
\definecolor{bblue}{rgb}{0.2,0.0,0.7}

\newcommand{\secref}[1]{Section \ref{sec:#1}}

\newtheorem{definition}{Definition}
\newtheorem{proposition}{Proposition}
\newtheorem{remark}{Remark}
\newtheorem{theorem}{Theorem}

\newcommand{\seclab}[1]{\label{sec:#1}}

\newcommand{\eqlab}[1]{\label{eq:#1}}
\renewcommand{\eqref}[1]{(\ref{eq:#1})}

\newcommand{\figref}[1]{Fig.~\ref{fig:#1}}
\newcommand{\figlab}[1]{\label{fig:#1}}
\newcommand{\propref}[1]{Proposition~\ref{proposition:#1}}
\newcommand{\proplab}[1]{\label{proposition:#1}}

\newcommand{\defnref}[1]{Definition~\ref{definition:#1}}
\newcommand{\defnlab}[1]{\label{definition:#1}}

\newcommand{\remref}[1]{Remark~\ref{remark:#1}}
\newcommand{\remlab}[1]{\label{remark:#1}}
\newcommand{\thmref}[1]{Theorem~\ref{theorem:#1}}
\newcommand{\thmlab}[1]{\label{theorem:#1}}

\newcommand{\asuref}[1]{Assumption~\ref{assumption:#1}}
\newcommand{\asulab}[1]{\label{assumption:#1}}

\newtheorem{asu}{Assumption}

\makeatletter
\newcommand{\tpitchfork}{%
	\vbox{
		\baselineskip\z@skip
		\lineskip-.52ex
		\lineskiplimit\maxdimen
		\m@th
		\ialign{##\crcr\hidewidth\smash{$-$}\hidewidth\crcr$\pitchfork$\crcr}
	}%
}
\makeatother

\title{Geometric singular perturbation analysis of the multiple-timescale Hodgkin-Huxley equations}
\author{Panagiotis Kaklamanos, Nikola Popovi\'c, and Kristian Uldall Kristiansen\footnote{P. Kaklamanos and N. Popovi\'c: Maxwell Institute for Mathematical Sciences and School of Mathematics, University of Edinburgh; K. U. Kristiansen: Department of Applied Mathematics and Computer Science, Technical University of Denmark}}
\date{}

\begin{document}
	\maketitle
	
	\pagestyle{myheadings}
	\thispagestyle{plain}
	
	\begin{abstract}
		{We present a novel and global three-dimensional reduction of a non-dimensionalised version of the four-dimensional Hodgkin-Huxley equations} [J. Rubin  and  M.  Wechselberger, Giant  squid--hidden  canard: the 3D geometry  of the Hodgkin-Huxley model, Biological Cybernetics, 97 (2007), pp. 5--32] that is based on geometric singular perturbation theory (GSPT). We investigate the dynamics of the resulting {three-dimensional} system in two parameter regimes in which the flow evolves on three distinct timescales. Specifically, we demonstrate that the system exhibits bifurcations of oscillatory dynamics and complex mixed-mode oscillations (MMOs), in accordance with the geometric mechanisms introduced in [P. Kaklamanos, N. Popovi\'c, and K. U. Kristiansen, Bifurcations of mixed--mode oscillations in three--timescale systems: An extended prototypical example, Chaos: An Interdisciplinary Journal of Nonlinear Science, 32 (2022), p. 013108], {and we classify the various firing patterns {in terms of} the external applied current}. {While such patterns have been documented in [S. Doi, S. Nabetani, and S. Kumagai, Complex nonlinear dynamics of the Hodgkin-Huxley equations induced by time scale changes, Biological Cybernetics, 85 (2001), pp. 51--64] for the multiple-timescale Hodgkin-Huxley equations, we elucidate the geometry that underlies the transitions between them, which had not been previously emphasised.}
	\end{abstract}
	
	\maketitle
	
	\section{Introduction} 
	The Hodgkin-Huxley (HH) equations comprise a model that describes the generation of action potentials in the squid giant axon \cite{HHmain}. The particular importance of these equations is due to the fact that they constitute one of the most successful mathematical models for the quantitative description of biological phenomena, as the underlying formalism is directly applicable to many types of neurons and other cells \cite{doi2001complex, izhikevich2007dynamical}.
	
	The three currents that the squid axon carries are the voltage-gated persistent potassium ($ K^+ $) 
	current, the voltage-gated transient sodium ($ Na^+ $) current, and the Ohmic leak current, which are given by the expressions
	\begin{align}
	I_{Na}(V,m,h)=g_{Na}\lp V-{E}_{Na}\rp m^3h, \quad I_{K}(V,n) = {g}_K\lp V-{E}_K\rp n^4,\quad\text{and}\quad I_L(V)={g}_L\lp V-{E}_L\rp,\eqlab{currents}
	\end{align} 
	respectively; here, $ V $ denotes the membrane potential, in units of mV, and $ m $ and $ h $ are the activation and
	inactivation variables of the $ Na^+ $ ion channel, respectively, while $ n $ is the activation variable of the $ K^+ $ channel. The conductance density $g_x$ and the Nernst potentials $E_x$ ($x=Na,K,L$) are given in units of mS/cm$^2$ and mV, respectively. 
	With these definitions, the original HH equations \cite{HHmain, doi2001complex,izhikevich2007dynamical} read 
	\begin{align}
	    \begin{aligned}
	    C\dot{V} &= {I}-I_{Na}(V,m,h)-I_{K}(V,n)-I_L(V), \\ 
		\dot{m}&= \frac{1}{\tau_m\hat{t}_m\lp V\rp} \lp m_\infty \lp V\rp - m\rp, \\
		\dot{h}&= \frac{1}{\tau_h \hat{t}_h\lp V\rp} \lp h_\infty \lp V\rp - h\rp, \\
		\dot{n}&= \frac{1}{\tau_n\hat{t}_n\lp V\rp} \lp n_\infty \lp V\rp - n\rp,
	    \end{aligned}
	    \eqlab{original}
	\end{align}
	where {$C$ denotes the capacitance density}, in units of $\mu$F/cm$^2$, and $ I $ is the applied current, in units of $ \mu \tn{A}/\tn{cm}^2 $, {while the constants $\tau_x$ and the voltage-dependent parameters $\hat{t}_x(V)$ ($x=m,h,n$) are associated with the characteristic timescales of the corresponding channel gates.} Moreover,
		\begin{gather}
		\hat{t}_x\lp V\rp = \frac{1}{\alpha_x\lp V\rp+\beta_x\lp V\rp}\quad\text{and}\quad x_\infty \lp V\rp= \frac{\alpha_x\lp V\rp}{\alpha_x\lp V\rp+\beta_x\lp V\rp}, \quad\text{with }x = m,h,n, \eqlab{taus}
		\end{gather}
	where the functions 
	$ \alpha_x\lp V\rp $ and $ \beta_x\lp V\rp $ are defined as
	\begin{gather*}
	\alpha_m\lp V\rp = \frac{\lp 25-V\rp/10}{{\rm e}^{\lp 25-V\rp/10}-1}, \quad \alpha_h\lp v\rp=\frac{7}{100}{\rm e}^{-V/20}, \quad \alpha_n\lp v\rp = \frac{\lp 10-V\rp/100}{{\rm e}^{\lp 10-V\rp/10}-1}\\
	\beta_m\lp V\rp = 4{\rm e}^{-V/18}, \quad \beta_h\lp v\rp = \frac{1}{1+{\rm e}^{\lp 30-V\rp/10}}, 
	\quad\text{and}\quad \beta_n\lp v\rp=\frac{1}{4}{\rm e}^{-V/80} 
	\end{gather*}
	in the original form of \eqref{original}.
	In \cite{doi2001complex}, the parameters $ I $, $ \tau_m $, $ \tau_h $, and $ \tau_n $ are considered bifurcation parameters, while the values of the remaining parameters are fixed as
	\begin{gather*}
	C=1\ \mu\tn{F/cm}^2, \quad {g}_{Na}=120\ \tn{mS/cm}, \quad {g}_K=36\  \tn{mS/cm}, \quad {g}_L = 0.3\ \tn{mS/cm},\\
	{E}_{Na} = 115\ \tn{mV}, \quad {E}_K=-12\ \tn{mV}, \quad\text{and}\quad {E}_L = 10.599\ \tn{mV}.
	\end{gather*}
	
	When either $ \tau_h $ or $ \tau_n $ is large, it is demonstrated in \cite{doi2001complex} that the equations in \eqref{original} evolve on multiple timescales, which gives rise to mixed-mode dynamics. {Mixed-mode oscillations (MMOs) are trajectories that consist of small-amplitude oscillations (SAOs) alternating with large-amplitude excursions (LAOs); {see \figref{HStimeseries} for an illustration of the corresponding time series of $v$}. {In the two-timescale context, mixed-mode dynamics has been described via the so-called ``generalised canard mechanism" which is based on a combination of Fenichel's geometric singular perturbation theory (GSPT) \cite{fenichel1979geometric} and the desingularisation technique known as ``blow-up" \cite{krupa2001extending}. In that context, SAOs typically occur in the neighbourhood of folded singularities that are of folded node type \cite{wechselberger2005existence,brons2006mixed,desroches2012mixed}; trajectories that escape the vicinities of these folded singularities after undergoing SAOs are then reinjected by a suitable global return mechanism, giving rise to MMO trajectories. In the context of multiple-scale systems with more than two scales, mixed-mode dynamics frequently arises due to the presence of folded singularities of folded saddle-node type \cite{krupa2008mixed,krupa2010local,letson2017analysis} which give rise to complex and irregular oscillatory behaviour, as is also the case here. Thus, for instance,} the oscillatory trajectories illustrated in panels (a) and (b) of \figref{HStimeseries} exhibit \textit{double epochs} of {{slow dynamics}}. {Here, we refer to time series with double epochs as featuring slow segments or SAOs which are centred around two distinct values of $v$ and which we hence denote by slow dynamics ``above'' and ``below'', respectively. By contrast, in panels (d) and (e), trajectories only exhibit epochs of slow dynamics ``below" which are separated by large excursions (spikes); finally, the trajectory in panel (f) exhibits no epochs of slow dynamics. {As will become apparent from our analysis, in that case only dynamics on the fastest and an additional, ``intermediate" timescale occur. (Here, we note that the scaling of the $t$-axes in \figref{HStimeseries} differs by panel; correspondingly, the ``below'' segments in panel (f) are in fact steeper than the ``above'' segments in panels (a), (b), and (c).)}. 

	\begin{figure}[p!ht!]
		\centering
		\begin{subfigure}[b]{0.45\textwidth}
			\centering
			\includegraphics[scale = 0.45]{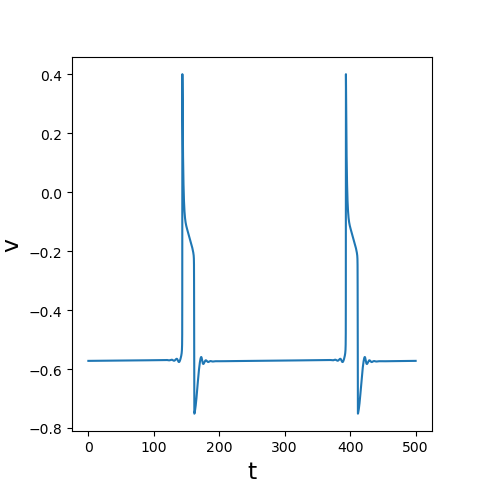}
			\caption{$ {I} = 20.051 $}
		\end{subfigure}
		~
		\begin{subfigure}[b]{0.45\textwidth}
			\centering
			\includegraphics[scale = 0.45]{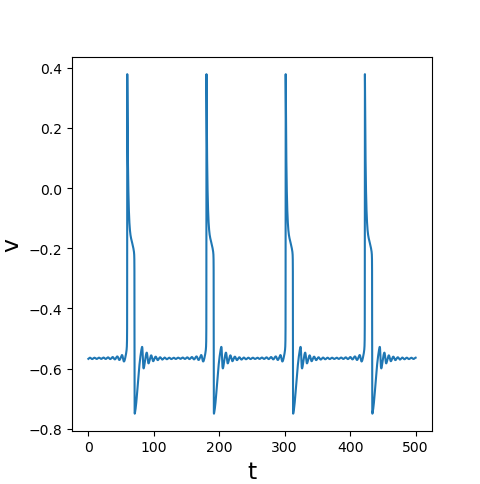}
			\caption{$ {I} = 23.051 $}
		\end{subfigure}
		\\
		\centering
		\begin{subfigure}[b]{0.45\textwidth}
			\centering
			\includegraphics[scale = 0.45]{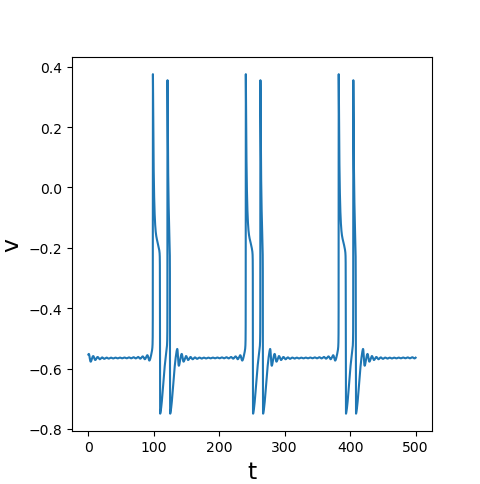}
			\caption{$ {I} = 23.5 $}
		\end{subfigure}
		~
		\begin{subfigure}[b]{0.45\textwidth}
			\centering
			\includegraphics[scale = 0.45]{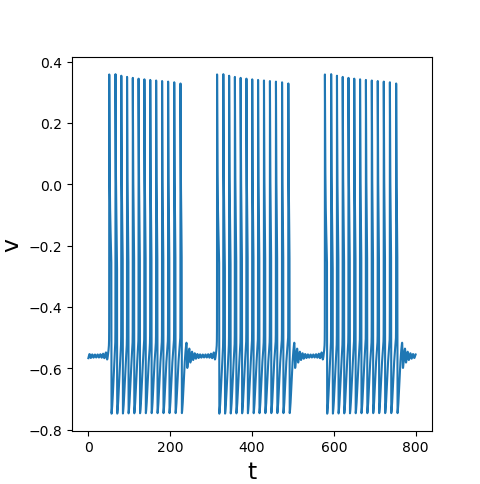}
			\caption{$ {I} = 26.03452346 $}
		\end{subfigure}
		\\
		\centering
		\begin{subfigure}[b]{0.45\textwidth}
			\centering
			\includegraphics[scale = 0.45]{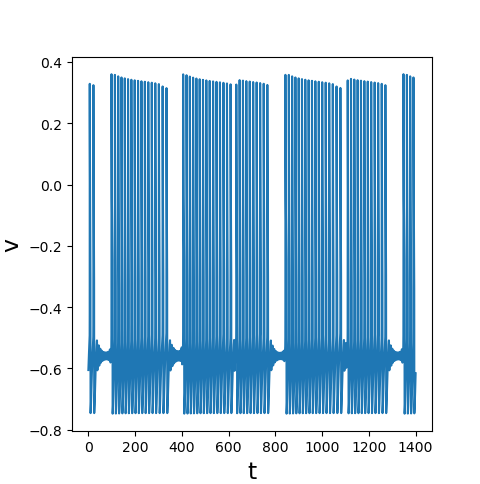}
			\caption{$ {I} = 26.1209956 $}
		\end{subfigure}
		~
		\begin{subfigure}[b]{0.45\textwidth}
			\centering
			\includegraphics[scale = 0.45]{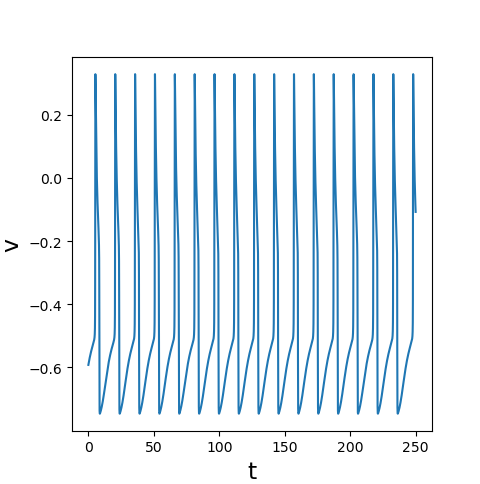}
			\caption{$ {I} = 26.2 $}
		\end{subfigure}
		\caption{Time series of the $ v $-variable in the dimensionless HH model, Equation~\eqref{hh}, for $ \varepsilon =0.0083 $, $\tau_h=40$,  and varying values of $ {I} =k_vg_{Na}\bar{I}$: panels (a) and (b) illustrate oscillatory trajectories with double {epochs of slow dynamics}; panel (c) gives an example of a {``transitional''} MMO with double {epochs of slow dynamics} separated by LAOs; panels (d) and (e) show MMO trajectories with single SAO epochs. Finally, panel (f) displays relaxation oscillation. Qualitatively similar time series have been documented in \cite[Figure 2]{doi2001complex}; the underlying geometric mechanisms in phase space that are responsible for the qualitative properties of the associated MMOs are described in \secref{hslow}.}
		\figlab{HStimeseries}
	\end{figure}
	
	In \cite{kaklamanos2022bifurcations}, geometric mechanisms were described that distinguish oscillatory trajectories with double epochs of {{slow dynamics} -- as shown in  panels (a) and (b) -- from those with single epochs thereof}, as seen in panels (d) and (e), as well as from relaxation oscillation in a general family of three-timescale systems; {see panel (f)}. A further distinction can be made between the trajectories with double epochs in panels (a) and (b) and the {``transitional''} mixed-mode dynamics observed in panel (c) in which {epochs of slow dynamics} above and below are separated by LAOs; see \secref{hslow} for details.} While the above qualitative distinction was documented in \cite{doi2001complex}, the mechanisms that are responsible for transitions between the various firing patterns were not emphasised.

	To our knowledge, the first attempt at studying the multiple-timescale HH model, Equation~\eqref{original}, on the basis of GSPT was made by Rubin and Wechselberger \cite{rubin2007giant,rubin2008selection}. They applied the scaling
	\begin{gather}
	\begin{gathered}
	v = \frac{V}{k_v}, \quad \tau = \frac{t}{k_t}, \quad \bar{I} =\frac{I}{k_v g_{Na}}, \\
	\bar{E}_x = \frac{E_x}{k_v}, \quad\text{and}\quad \bar{g}_x = \frac{g_x}{g_{Na}}\quad \tn{ for }~x=Na,K,L,
	\end{gathered}\eqlab{scaling}
	\end{gather}
	with $k_v = 100$ \tn{mV} and $k_t = 1$ ms, and considered the following dimensionless version of Equation~\eqref{original},
		\begin{align}
		\begin{aligned}
		\epsilon\dot{v} &= \bar{I}-\lp v-\bar{E}_{Na}\rp m^3h-\bar{g}_K\lp v-\bar{E}_K	\rp n^4-\bar{g}_L\lp v-\bar{E}_L\rp, \\ 
		\dot{m}&= \frac{1}{\tau_m\hat{t}_m\lp v\rp} \lp m_\infty \lp v\rp - m\rp, \\
		\dot{h}&= \frac{1}{ \tau_h\hat{t}_h\lp v\rp} \lp h_\infty \lp v\rp - h\rp, \\
		\dot{n}&= \frac{1}{\tau_n\hat{t}_n\lp v\rp} \lp n_\infty \lp v\rp - n\rp, 
		\end{aligned}
		\eqlab{hh}
		\end{align}
	where 
	\begin{align}
	\epsilon = \frac{C}{k_t\cdot g_{Na}} \simeq 0.0083. \eqlab{eps}
	\end{align}
	Here, the functions $ \hat{t}_x\lp v \rp $ and $ x_\infty\lp v \rp $ are
	defined as in \eqref{taus}, with
	\begin{gather*}
	\alpha_m\lp v\rp = \frac{\lp v+40\rp/10}{1-{\rm e}^{-\lp v+40\rp/10}}, \quad \alpha_h\lp v\rp=\frac{7}{100}{\rm e}^{-\lp v+65\rp/20}, \quad \alpha_n\lp v\rp = \frac{\lp v+55\rp/100}{1-{\rm e}^{-\lp v+55\rp/10}},\\
	\beta_m\lp v\rp = 4{\rm e}^{-\lp v+65\rp/18}, \quad \beta_h\lp v\rp = \frac{1}{1+{\rm e}^{-\lp v+35\rp/10}}, 
	\quad\text{and}\quad\beta_n\lp v\rp=\frac{1}{4}{\rm e}^{-\lp v+65\rp/80};
	\end{gather*}
	in addition, the parameters in \eqref{hh} now read
	\begin{gather}
	\begin{gathered}
	\bar{I} = \frac{I}{k_ug_{Na}}, \quad \bar{g}_K=0.3, \quad \bar{g}_L = 0.0025,\\
	\bar{E}_{Na} = 0.5, \quad \bar{E}_K=-0.77, \quad\text{and}\quad \bar{E}_L = -0.544.
	\end{gathered}\eqlab{parbars}
	\end{gather}
	The dimensionless Equation~\eqref{hh} captures the different firing patterns that emerge in
	Equation~\eqref{original} as $I$ is varied, whereas the remaining parameters are fixed, as documented 
	in \cite{doi2001complex} and illustrated in \figref{HStimeseries}.
	
	We note that, due to the non-dimensionalisation in \eqref{hh} by Rubin and Wechselberger \cite{rubin2007giant}, the corresponding $ I $-values differ from those in \cite{doi2001complex}; moreover, we remark that for ease of comparison and numerical convenience, \figref{HStimeseries} refers to $I$ rather than to its rescaled counterpart $\bar I$ in Equation~\eqref{hh}.
	
	Under the above assumptions, the parameter $\epsilon$ defined in \eqref{eps} can be considered a small perturbation parameter. Moreover, the order of magnitude of the characteristic timescale $[\hat{t}_m(v)]^{-1}$ is larger than those of $[\hat{t}_h(v)]^{-1}$ and $[\hat{t}_n(v)]^{-1}$, as can be seen in \figref{txs}.
	
	\begin{figure}[ht!]
		\centering
		\begin{subfigure}[b]{0.49\textwidth}
			\centering
			\includegraphics[scale = 0.22]{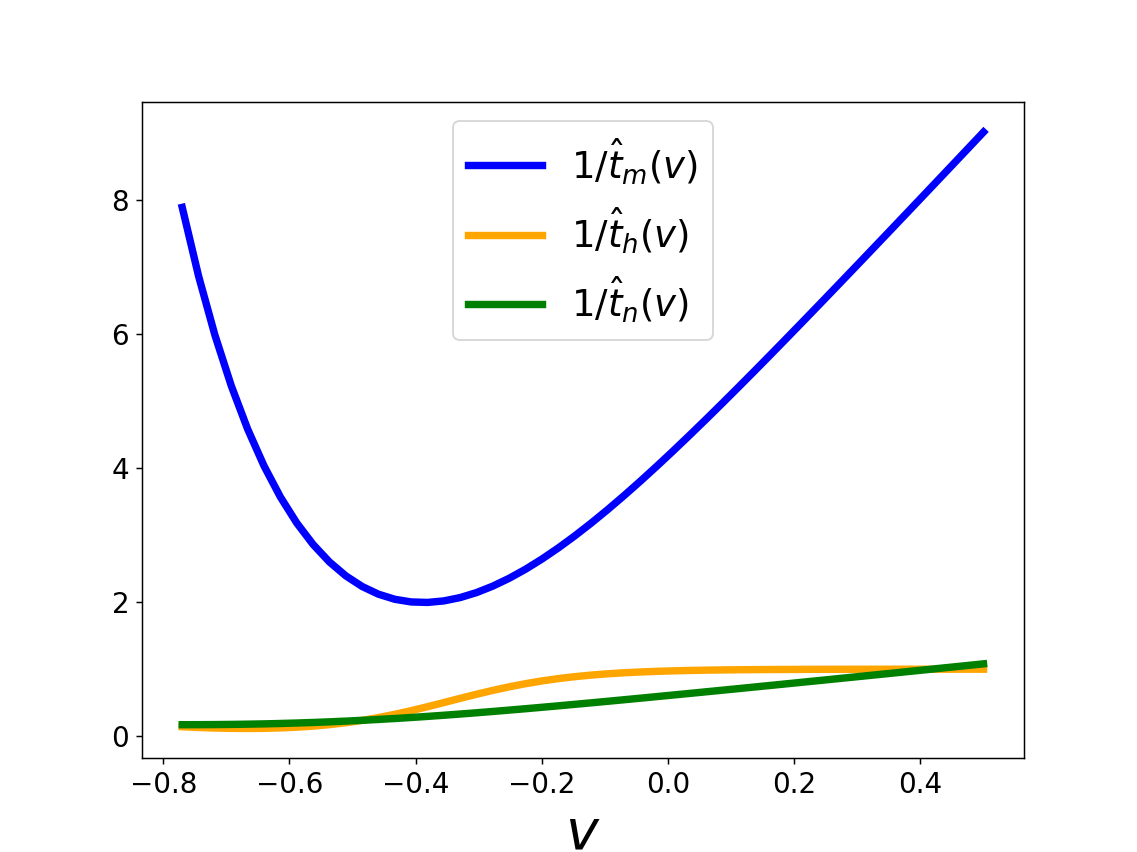}
			\caption{}
		\end{subfigure}
		\begin{subfigure}[b]{0.49\textwidth}
			\centering
			\includegraphics[scale = 0.22]{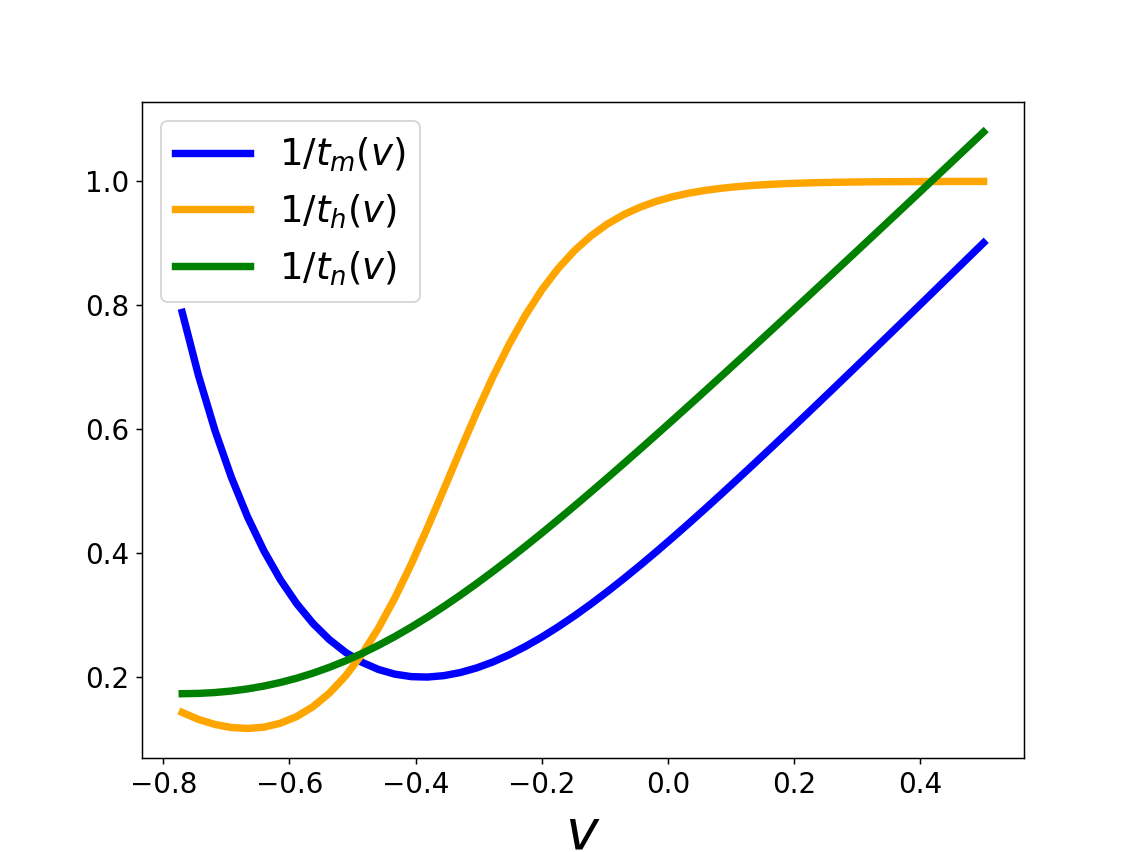}
			\caption{}
		\end{subfigure}
		\caption{Characteristic timescales (a) $[\hat{t}_x(v)]^{-1}$ and (b) $[{t}_x(v)]^{-1}$ for $x = m,h,n$.}
		\figlab{txs}
	\end{figure}
	
	{Following \cite{rubin2007giant},} we consider the rescaling
	\begin{align}
	T_x = \max_{v\in(\bar{E}_K, \bar{E}_{Na})} \frac{1}{\hat{t}_x(v)}\quad\text{for }x=m,h,n\eqlab{Tx}
	\end{align}
	and we define
	\begin{align}
	{{t}_x(v)} := {T_x\hat{t}_x(v)}. \eqlab{tx}
	\end{align}
	Denoting further
	\begin{align}
	V\lp v,m,h,n\rp&:= \bar{I}-\lp v-\bar{E}_{Na}\rp m^3h-\bar{g}_K\lp v-\bar{E}_K\rp n^4-\bar{g}_L\lp v-\bar{E}_L\rp, \eqlab{V4}\\
	M(v,m) &:=\frac{m_\infty \lp v\rp - m}{t_m\lp v\rp}, \eqlab{M4} \\
	H(v,h) &:=\frac{h_\infty \lp v\rp - h}{t_h\lp v\rp}, \eqlab{H4} \\
	N(v,n) &:=\frac{n_\infty \lp v\rp - n}{t_n\lp v\rp},\quad\tn{and} \eqlab{N4}\\
	\delta_x &:= \frac{T_x}{\tau_x},\quad\text{with }x = m,h,n,
	\end{align}
	one can write
		\begin{align}
		\begin{aligned}
		\epsilon\dot{v} &= V\lp v,m,h,n\rp \\ 
		\dot{m}&= \delta_mM(v,m), \\
		\dot{h}&= \delta_hH(v,h), \\
		\dot{n}&= \delta_nN(v,n). 
		\end{aligned}
		\eqlab{hh1f}
		\end{align}
	In the following, we will suppress the dependence of the system in \eqref{hh1f} on the parameter $\bar{I}$, and will only remark on it as required.
	
	From \figref{txs}, it is apparent that $T_m \simeq 10$, while $T_h\simeq T_n \simeq 1$. Based on this observation, in \cite{rubin2007giant} it was assumed that $\delta_m^{-1}=\mc{O}(\epsilon)$; thus, \eqref{hh1f} was written as
		\begin{align}
		\begin{aligned}
		\epsilon\dot{v} &= V\lp v,m,h,n\rp,\\ 
		\epsilon\dot{m}&= M(v,m),\\
		\dot{h}&= \delta_hH(v,h),\\
		\dot{n}&= \delta_nN(v,n), 
		\end{aligned}
		\eqlab{hh2f}
		\end{align}
	where $(v,m)$ are the fast variables and $(h,n)$ are the slow ones. On the basis of centre manifold theory, Rubin and Wechselberger \cite{rubin2007giant} then derived the following three-dimensional reduction of \eqref{hh2f},
		\begin{align}
		\begin{aligned}
		\epsilon\dot{v} &= \bar{I}-\lp v-\bar{E}_{Na}\rp m^3_{\infty}(v)h-\bar{g}_K\lp v-\bar{E}_K\rp n^4-\bar{g}_L\lp v-\bar{E}_L\rp = V\lp v, m_\infty,h,n\rp,  \\ 
		\dot{h}&= \delta_hH(v,h),\\
		\dot{n}&= \delta_nN(v,n) ,
		\end{aligned}
		\eqlab{RW}
		\end{align}
	where the variable $m$ has been eliminated. However, in their analysis, they considered $ 0<\epsilon\ll 1 $ and $ \delta_x = \mc{O}(1) $ {($ x = h,n $)}.  As a consequence, only MMO trajectories with single epochs of SAOs were documented in previous works \cite{rubin2007giant,rubin2008selection}, in contrast to \cite[Figure~2]{doi2001complex}; cf.~panels (d), (e), and (f) of  \figref{HStimeseries}. As argued in \cite{desroches2012mixed}, in the two-timescale case -- where SAOs arise due to the presence of a folded node singularity \cite{wechselberger2005existence} -- the parameter regimes in which MMOs with double SAO epochs are apparent seem very narrow. By contrast, such MMOs become more prominent in a three-timescale context, which further attests to our claim that the physiologically relevant mixed-mode dynamics of the HH model, Equation~\eqref{hh}, cannot be fully understood via a standard two-timescale analysis.
	
	\begin{remark}
		The timescale separation at all levels can be controlled by merely varying the parameters $\tau_x$ ($x = m,h,n$) in \eqref{hh} or, equivalently, by variation of $\delta_x$ in \eqref{hh1f}. 
	\end{remark}
	
	\bigskip
	
	\subsection*{{Main results}}
	{In this work, we incorporate the assumption that the variable $ m $ is slower than $ v $, but faster than $ h $ and $ n $, in accordance with \cite[Remark 1]{rubin2007giant}, which we formulate as follows:
		\begin{asu}
			{Let $\varepsilon :=\delta_m^{-1}$, and define $\gamma>0$ such that
			\begin{align*}
			\epsilon = \gamma\varepsilon.
			\end{align*}
			Then, we assume that 
			\begin{align*}
			0<\epsilon \ll \gamma,\varepsilon\ll 1
			\end{align*}
			in \eqref{hh1f}, with $\epsilon, \gamma$, and $\varepsilon$ positive and sufficiently small.}
			\asulab{gamma}
		\end{asu}
	{Our first main result is the following:} for $\epsilon \simeq 0.0083$ and $\delta_m\simeq10$, we have $\varepsilon \simeq 0.1$ and $\gamma = 0.083$}. Then, {Assumption 1 allows} us to apply GSPT to derive a global three-dimensional reduction of Equation~\eqref{hh} of the form
	\begin{align}\eqlab{red-gen}
	\varepsilon \dot{v} = U(v,h,n;\bar{I}, \gamma,\varepsilon,\delta_h, \delta_n), \quad
	\dot{h}= \delta_hH(v,h), \quad
	\dot{n}= \delta_nN(v,n);
	\end{align}
    see \thmref{3reduction} in \secref{reduction} below. 
 	Motivated by \cite{doi2001complex}, we then consider the two regimes where either $ h $ or $ n $ in Equation~\eqref{red-gen} is assumed to be evolving on the slowest timescale. In other words, we first take $ \delta_h >0 $ to be small, with $ \delta_n = \mc{O}(1) $; then, we consider the regime where $ \delta_h = \mc{O}(1) $, with $ \delta_n >0 $ small, and we show that these two regimes are not fundamentally different in terms of their singular geometry and of the resulting mixed-mode dynamics. 
 	
 	{Under these additional assumptions, our second main result is that transitions from mixed-mode dynamics with double epochs of slow dynamics via single-epoch MMOs to relaxation oscillation in the three-timescale system, Equation~\eqref{red-gen}, can be explained by the geometric mechanisms introduced in \cite{kaklamanos2022bifurcations}. {Specifically, we identify invariant manifolds in \eqref{red-gen} and decompose the dynamics into segments evolving on fast, intermediate, and slow timescales. Then, qualitative properties of oscillatory trajectories of \eqref{red-gen} can be deduced based on the geometric configuration of these manifolds and the reduced flows thereon; in particular, we show that these properties are determined by the relative positions of so-called folded singularities, which can be either ``connected", ``aligned", or ``remote", in accordance with the classification in \cite{kaklamanos2022bifurcations}. The rescaled applied current $\bar I$ will emerge as the relevant bifurcation parameter}. We illustrate the underlying geometry schematically in \figref{fig11}; see the corresponding caption for details. In sum, we thus obtain a full description of the dynamics of \eqref{red-gen} under the assumption of timescale separation between the three variables therein, which is motivated by the physiological properties of the particular system under consideration.}
 	{We emphasise that, although the original HH equations in \eqref{original} have been studied extensively and are by now fairly well understood, the methodology presented here extends to a wide class of multiple-scale systems that can be expressed in an HH-type formalism, such as models for calcium-induced calcium release, cardiac dynamics, circadian and homeostatic processes, and other physiological phenomena \cite{cloete2020dual,diekman2021circadian,forger2017biological, iosub2015calcium, kugler2018early, pages2019model, takano2021highly}, to mention but a few. }

	\begin{figure}[ht!]
			\centering
	    \begin{subfigure}[b]{0.3\textwidth}
		\centering
		\includegraphics[scale = 0.37]{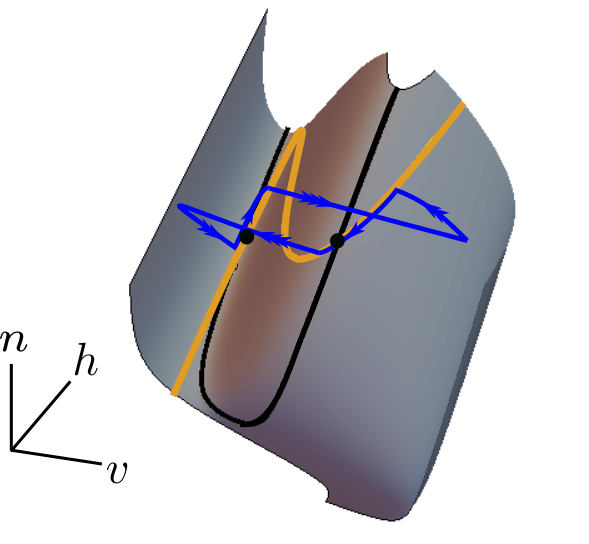}
		\caption{}
	\end{subfigure}
	\quad
	\begin{subfigure}[b]{0.3\textwidth}
				\centering
				\includegraphics[scale = 0.37]{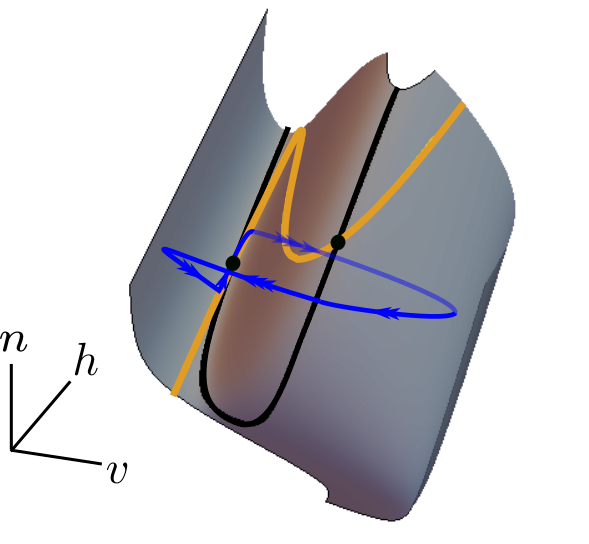}
				\caption{}
			\end{subfigure}
		\quad
		\begin{subfigure}[b]{0.3\textwidth}
			\centering
			\includegraphics[scale = 0.37]{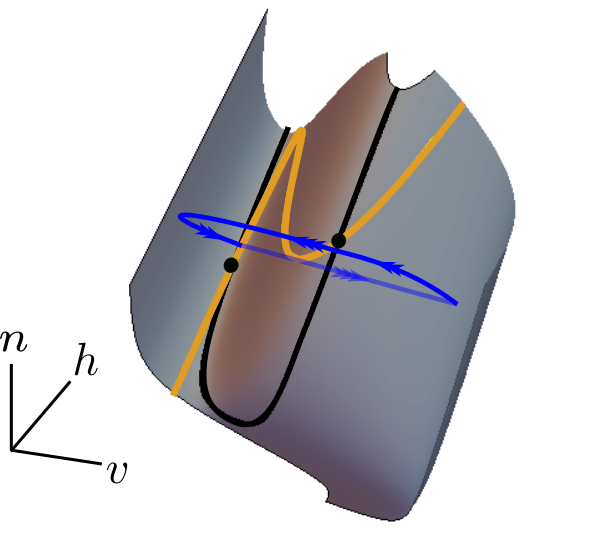}
			\caption{}
		\end{subfigure}
		\caption{{Schematic illustration of oscillatory trajectories in the $(v,h,n)$-phase space of \eqref{red-gen}: (a) MMO with double epochs of slow dynamics; (b) MMO with a single epoch of {slow dynamics}; (c) relaxation oscillation. Trajectories are attracted towards a two-dimensional invariant manifold on a fast timescale; they then evolve close to that manifold on an intermediate timescale. {The locations of folded singularities \cite{kaklamanos2022bifurcations,letson2017analysis,szmolyan2001canards} (black dots) relative to one another and the reduced flow on the two-dimensional manifold determine whether trajectories are attracted to one-dimensional invariant manifolds (orange curves), whereon they evolve on a slow timescale, or whether they reach the vicinity of a non-hyperbolic set (black curves), where they undergo fast jumps. As the parameter $\bar{{I}}$ is varied, the configuration of these geometric objects and the flows thereon changes, resulting in the different qualitative behaviours outlined above; cf. also \figref{ibif} below. The main distinction is due to the folded singularities being either ``connected'', ``aligned'', or ``remote'' \cite{kaklamanos2022bifurcations}, as described in detail in \secref{hslow} and \secref{nslow}.}}}
		\figlab{fig11}
		\end{figure}
	
	\bigskip
	
	The article is organised as follows. In \secref{reduction}, we present a novel, global geometric reduction of Equation~\eqref{hh} after elimination of the variable $ m $. We relate the dynamics of the resulting, three-dimensional slow-fast system, Equation~\eqref{red-gen}, to that of Equation~\eqref{hh1f}, as well as to the local centre manifold reduction proposed by Rubin and Wechselberger \cite{rubin2007giant}, Equation~\eqref{RW}. We then consider the scenario where Equation~\eqref{red-gen} exhibits dynamics on three timescales, in which case the full system in \eqref{hh} -- or, equivalently, in \eqref{hh1f} -- exhibits dynamics on four timescales. In \secref{hslow} and \secref{nslow}, the variables $ h $ and $ n $ are taken to be slow, respectively; we demonstrate that the geometric mechanisms proposed in \cite{kaklamanos2022bifurcations} can explain bifurcations of MMOs that have been previously documented, but not emphasised in the context of GSPT, in the literature \cite{doi2001complex}. We conclude the article in \secref{HHconclusion} with a summary, and an outlook to future work. Finally, in Appendix~\ref{partpert}, we give a brief overview of SAO-generating mechanisms in a three-timescale context.
	
	\section{Global multiple-timescale reduction}
	\seclab{reduction}
	Rescaling time in \eqref{hh1f} as $ \sigma =  t/\epsilon$, we obtain the ``fast formulation''
		\begin{align}
		\begin{aligned}
		{v}' &= V\lp v,m,h,n\rp, \\ 
		{m}' &= \epsilon \delta_m M(v,m),\\
		{h}' &= \epsilon\delta_h H(v,h) ,\\
		{n}' &= \epsilon\delta_n N(v,n),
		\end{aligned}
		\eqlab{hhmod-fast0}
		\end{align}
	where the prime denotes differentiation with respect to the new time $\sigma$. {With \asuref{gamma},} the system in \eqref{hhmod-fast0} is then written as
		\begin{align}
		\begin{aligned}
		{v}' &= V\lp v,m,h,n\rp, \\ 
		{m}' &= \gamma  M(v,m),\\
		{h}' &= \gamma \varepsilon\delta_h H(v,h) ,\\
		{n}' &= \gamma \varepsilon\delta_n N(v,n). 
		\end{aligned}
		\eqlab{hhmod-fast}
		\end{align}
	We will restrict our analysis to a compact set $ \mc{D}\subset(\bar{E}_{K}, \bar{E}_{Na})\times(0,1)^3$,
	as is also done in \cite{rubin2007giant}; recall \eqref{Tx}. {The constraint that $v\in(\bar{E}_{K}, \bar{E}_{Na})$ is a physiological one, set by the underlying Nernst potentials: when $v=\bar{E}_K$, respectively when $v=\bar{E}_{Na}$, the corresponding current, denoted by $\bar{I}_K$, respectively by $\bar{I}_{Na}$, is zero; see \cite{izhikevich2007dynamical} for details. A mathematical interpretation of this restriction is given in \remref{nernst} below.}
	
	Equation~\eqref{hhmod-fast} is singularly perturbed with respect to the small parameter $\gamma$ and features three distinct timescales if $ \delta_h, \delta_n = \mc{O}(1) $, with $v$ the fast variable, $m$ the intermediate one, and $h$ and $n$ the slow ones. The singular limit of $ \gamma=0$ gives the one-dimensional ``layer problem'' with respect to the fastest timescale,
		\begin{align}
		\begin{aligned}
		v' &= V\lp v,m,h,n\rp, \\ 
		m' &= 0,\\
		h' &= 0,\\
		n' &= 0;
		\end{aligned}
		\eqlab{hhlay}
		\end{align}
	{the $v$-nullcline of \eqref{hhlay}} defines the three-dimensional critical manifold ${\mc{M}_1}$ as
	\begin{align}
	V\lp v,m,h,n\rp =0. \eqlab{M0}
	\end{align}
	Since
	\begin{align}
	\partial_vV\lp v,m,h,n\rp &= -m^3h-\bar{g}_Kn^4 - \bar{g}_L<0 \quad \tn{for} \quad (v,m,h,n)\in \mc{D}, \eqlab{dV}
	\end{align}
	the manifold ${\mc{M}_1}$ is normally hyperbolic and attracting everywhere in $ \mc{D} $. Let
	\begin{align}
	\mu(v,h,n) := \left[\frac{\bar{I}-\bar{g}_K\lp v-\bar{E}_K\rp n^4-\bar{g}_L\lp v-\bar{E}_L\rp}{\lp v-\bar{E}_{Na}\rp h}\right]^{\frac{1}{3}}; \eqlab{mu}
	\end{align}
	then, by \eqref{M0}, ${\mc{M}_1} $ can be written as a graph of $ m $ over $(v,h,n)$, with
	\begin{align*}
	m = \mu(v,h,n).
	\end{align*}
	Therefore, for $\gamma>0$ sufficiently small, the dynamics of Equation~\eqref{hhmod-fast} can be effectively reduced to a three-dimensional multiple-timescale system, which constitutes our first main result in this work. 

	\begin{theorem}	\thmlab{3reduction}
		There exists $\gamma_0$ sufficiently small such that, for every $\gamma\in(0,\gamma_0)$, the HH equations in \eqref{hhmod-fast} admit a three-dimensional attracting {slow} manifold {$\mathcal{M}_{1\gamma} $} that is diffeomorphic, and {$ \mc{O}\lp \gamma\rp$-close to}, {the compact manifold $\mathcal{M}_1$} in the Hausdorff distance. The manifold $ \mathcal{M}_{1\gamma} $ can be written as a graph $m = \mu(v,h,n)+\mc{O}(\gamma)$ and is locally invariant under the flow of
			\begin{align}
			\begin{aligned}
			\begin{split}
			{v}' &= {\frac{m_\infty(v)-\mu(v,h,n)}{t_m(v)\partial_v\mu(v,h,n)} -\varepsilon\delta_hH(v,h)\frac{\partial_h\mu(v,h,n)}{\partial_v\mu(v,h,n)}-\varepsilon\delta_nN(v,n)\frac{\partial_n\mu(v,h,n)}{\partial_v\mu(v,h,n)}+ \mathcal{O}\lp \gamma\rp} \\ &=:U(v,h,n;\gamma,\varepsilon,\delta_h,\delta_n),
			\end{split} \\ 
			{h}' &= \varepsilon\delta_hH(v,h),  \\
			{n}' &= \varepsilon\delta_nN(v,n). 
			\end{aligned}
			\eqlab{3red}
			\end{align}
	\end{theorem}
	\begin{proof}
		The existence of $ {\mathcal{M}_{1\gamma}} $ and its diffeomorphic relation to, and distance from, $ \mathcal{M}_{1} $ follow from Fenichel's First Theorem \cite{fenichel1979geometric}, since the critical manifold $ \mathcal{M}_{1} $ is normally hyperbolic everywhere in $\mc{D}$, by \eqref{M0} and \eqref{dV}. In particular, $ {\mathcal{M}_{1\gamma}} $ is given as a  graph of $ m $ over $ (v,h,n) $, with
		\begin{align}
		m = \mu	(v,h,n)+\gamma\mu_{1}(v,h,n;\gamma,\varepsilon),
		\eqlab{mu2}
		\end{align} 
		{which is smooth jointly in $ v,h,n, \gamma $, and $ \varepsilon $.} The attractivity of $ {\mathcal{M}_{1\gamma}} $ follows from Fenichel's Second Theorem \cite{fenichel1979geometric}, since, by \eqref{dV}, the manifold $ \mathcal{M}_{1} $ is attracting everywhere. (See also \cite{hek2010geometric} for a concise discussion of Fenichel's First and Second Theorems.) 
		
		Finally, Equation~\eqref{3red} for the flow on the slow manifold $ {\mathcal{M}_{1\gamma}} $ is obtained as follows. {The equations for $h$ and $n$ in \eqref{3red}} follow immediately from \eqref{hhmod-fast}, as the dynamics of $h$ and $n$ is independent of $m$. To derive {the $v$-equation in \eqref{3red}}, we first differentiate \eqref{mu2} with respect to $\sigma$. {(Here, we reiterate that $\sigma=t/\epsilon$ is the independent variable in \eqref{hhmod-fast0} and \eqref{hhmod-fast}.)} Next, we substitute {the right-hand sides of the equations for $m$, $h$, and $n$ in \eqref{hhmod-fast}} into the resulting expression and then solve for ${v}'$, collecting terms. 
	\end{proof}
	We remark that, for $\gamma = 0.083$ and $\bar{g}_L=0.0025$, the latter may be considered negligible with regard to the existence and hyperbolicity of $\mc{M}_{\gamma}$; for instance, the inequality in \eqref{dV} would also hold for $\bar{g}_L=0$ therein, by the definition of $\mathcal{D}$. 
	\thmref{3reduction} implies that, for $\gamma$ positive and sufficiently small, the dynamics of the full, four-dimensional HH model, Equation~\eqref{hhmod-fast}, is captured by the three-dimensional reduction in \eqref{3red}, where we have eliminated the variable $ m $ via the graph representation in \eqref{mu2}. {We remark that \eqref{3red} is a regular perturbation problem in $ \gamma $ in the slow formulation of GSPT; as the dynamics of  \eqref{3red} with $\gamma>0$ sufficiently small is therefore qualitatively similar, and $\mc{O}(\gamma)$-close, to that of \eqref{3red}$_{\gamma=0}$, we will restrict to considering the limit of $ \gamma=0 $ for simplicity.} {In the language of GSPT, we will hence approximate the slow flow on ${\mathcal{M}_{1\gamma}}$ by the reduced flow of {\eqref{hhmod-fast}} on $\mc{M}_1$, where the former is a regular perturbation of the latter \cite{fenichel1979geometric}:}
		\begin{align}
		\begin{aligned}
		{v}' &= \frac{m_\infty(v)-\mu(v,h,n)}{t_m(v)\partial_v\mu(v,h,n)} -\varepsilon\delta_hH(v,h)\frac{\partial_h\mu(v,h,n)}{\partial_v\mu(v,h,n)}-\varepsilon\delta_nN(v,n)\frac{\partial_n\mu(v,h,n)}{\partial_v\mu(v,h,n)}, \\ 
		{h}' &= \varepsilon \delta_hH(v,h),  \\
		{n}' &= \varepsilon \delta_nN(v,n). 
		\end{aligned}
		\eqlab{3redfast}
		\end{align}
	
	 {On the other hand, the systems in \eqref{3red} and \eqref{3redfast} are singular perturbation problems in standard form with respect to the parameter $\varepsilon$; hence, their dynamics for $ \varepsilon>0 $ sufficiently small can be studied via GSPT. Moreover, the $\mc{O}(\varepsilon)$-terms {on the right-hand side of \eqref{3red}} contribute to local qualitative phenomena, which in turn affect the resulting global mixed-mode dynamics. These terms can therefore not be neglected in our analysis; see Appendix \ref{partpert} and \cite{kaklamanos2022bifurcations} for details.}
	
	{{In the following, we will refer to the system in \eqref{3redfast} as \textit{our three-dimensional reduction of the HH equations.}} That system is written in the fast formulation of GSPT; the corresponding layer problem is obtained by setting $ \varepsilon=0 $ therein:}
	\begin{align}
	\begin{aligned}
	{v}' &= \frac{m_\infty(v)-\mu(v,h,n)}{t_m(v)\partial_v\mu(v,h,n)}, \\ 
	{h}' &= 0, \\
	{n}' &= 0. 
	\end{aligned}
	\eqlab{3redlay}
	\end{align}
	The set of equilibria of \eqref{3redlay} defines the two-dimensional critical manifold $\mc{M}_2$ as
	\begin{align}
	m_\infty(v)-\mu(v,h,n)= 0,
	\eqlab{em1}
	\end{align}
	which is equivalent to requiring that
	\begin{align}
	V(v,m_\infty(v), h,n) =0; \eqlab{em1V}
	\end{align}
	recall \eqref{V4} for the definition of $ V(v,m, h,n) $. {We emphasise that Equation~\eqref{3redfast} and the centre manifold reduction in \eqref{RW} that was introduced in \cite{rubin2007giant} admit the same critical manifold $ \mathcal{M}_2 $ as is defined by \eqref{em1} -- or, equivalently, by \eqref{em1V}. In other words, the algebraic constraint in \eqref{em1} corresponds to the $ v $-nullcline in \eqref{RW}, recall \eqref{V4}; see \cite{rubin2007giant} for details. However, although the factor of $ [\partial_v\mu(v,h,n)]^{-1} $ in {the $v$-equation of \eqref{3redfast}} is positive, there are qualitative differences in the dynamics of the two systems due to the $ \mathcal{O}(\varepsilon) $-terms that occur in {\eqref{3redfast}}, as discussed in Appendix~\ref{partpert}.}
	
	The manifold $ {\mc{M}_2} $ is normally hyperbolic everywhere except on {$ \mathcal{F}_{\mathcal{M}_2}$}, where
	\begin{align}
	\partial_v  \left[V(v,m_\infty(v), h,n) \right]= 0.
	\eqlab{Fm1V}
	\end{align}
	For future reference, we note that
    \begin{align}
    \begin{aligned}
    \partial_h \left[V(v,m_\infty(v),h,n)\right]&= -\lp v-\bar{E}_{Na}\rp m_\infty^3(v)>0\quad\text{and} \\
    \partial_n \left[V(v,m_\infty(v),h,n)\right]&= -4\bar{g}_K\lp v-\bar{E}_{k}\rp n^3<0 
    \end{aligned}
    \eqlab{parders}
    \end{align}
	for $ v \in (\bar{E}_K, \bar{E}_{Na})$, and we define
	\begin{align}
	\nu(v,h)&:= \bigg[\frac{\bar{I}-\lp v-\bar{E}_{Na}\rp m^3_{\infty}(v)h-\bar{g}_L\lp v-\bar{E}_L\rp}{\bar{g}_K\lp v-\bar{E}_K	\rp }\bigg]^\frac{1}{4}\quad\text{and} \eqlab{Fn}\\
	\eta(v,n)&:= \frac{\bar{I}-\bar{g}_K\lp v-\bar{E}_K	\rp n^4-\bar{g}_L\lp v-\bar{E}_L\rp}{\lp v-\bar{E}_{Na}\rp m_\infty^3(v)} \eqlab{Gh}
	\end{align}
	by solving \eqref{em1V} for $ n $ and $h$, respectively. {(We note that \eqref{parders} is consistent with the $Na^+$ current being depolarising in the HH formalism, while the $K^+$ current is hyperpolarising \cite{izhikevich2007dynamical}.)} The manifold $ {\mc{M}_2} $ can therefore be written globally either as a graph over $(v,h)$, with 
	\begin{align*}
	n= \nu(v,h),
	\end{align*} 
	or as a graph over $(v,n)$, with
	\begin{align*}
	h= \eta(v,n).
	\end{align*} 
	Finally, a parametric expression for {$ \mathcal{F}_{\mathcal{M}_2}$} can be obtained from \eqref{Fn} {or \eqref{Gh}}, by solving the system \{\eqref{em1V},\eqref{Fm1V}\}. Elementary calculation shows that $ \mc{F}_{{\mc{M}_2}} $ is a $ U $-shaped curve, which can be viewed as the tangential connection of two fold curves $ \mc{L}^\mp $; see \figref{crimaHH}. {This tangential connection occurs at a point on $\mc{F}_{{\mc{M}_2}}$ where $\partial_v^2[V(v,m_\infty(v), h,n)]= 0$.}
	
	It follows that the manifold $ {\mc{M}_2} $ is separated by {$ \mathcal{F}_{\mathcal{M}_2}$} into a normally attracting portion $ \mc{S}^a $, where $ \partial_v  \left[V(v,m_\infty(v), h,n) \right]<0 $, and a normally repelling portion $ \mc{S}^r $, where $ \partial_v  \left[V(v,m_\infty(v), h,n) \right]>0 $; cf.~again \figref{crimaHH}. We denote the normally hyperbolic submanifold of {$\mathcal{M}_2$} by $ \mc{S} =\mc{S}^a\cup\mc{S}^r $.
	
	\begin{remark}
		It is now also apparent that the restriction to $ v \in\lp \bar{E}_K, \bar{E}_{Na}\rp$ mathematically stems from Equations~\eqref{Fn} and \eqref{Gh}, since the denominators in those expressions tend to zero when $ v\to \bar{E}_{Na}$ or $ v\to \bar{E}_K$, respectively. {It therefore follows that the planes $\lb v = \bar{E}_K \rb$ and $\lb v = \bar{E}_{Na} \rb$  are vertical asymptotes of $\mc{M}_2$. Strictly speaking, in order to be able to apply GSPT, we will in fact restrict our analysis to sufficiently large compact subsets of $\mc{M}_2$.} \remlab{nernst}
	\end{remark}
	
	\begin{figure}[ht!]
		\centering
		\begin{subfigure}[b]{0.45\textwidth}
			\centering
			\includegraphics[scale = 0.25]{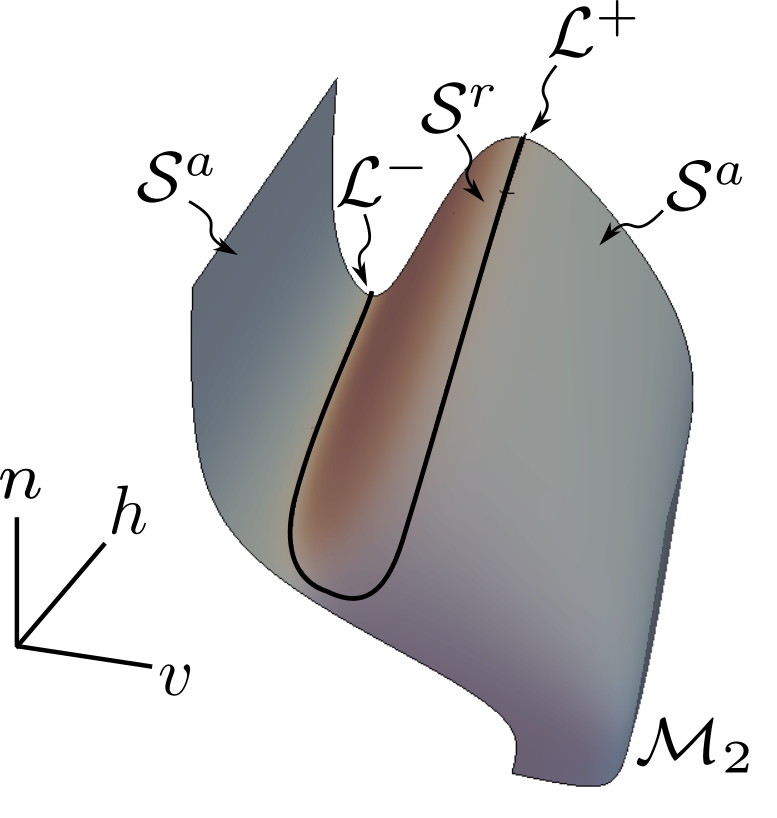}
		\end{subfigure}
		\caption{The critical manifold {$\mathcal{M}_2$} of the three-dimensional reduction, Equation~\eqref{3red}.}
		\figlab{crimaHH}
	\end{figure}
	
	Rescaling time in {the three-dimensional reduction} in \eqref{3redfast} as $t = \varepsilon\sigma$ gives the slow formulation
	\begin{align}
    \begin{aligned}
	\varepsilon \dot{v} &= \frac{m_\infty(v)-\mu(v,h,n)}{t_m(v)\partial_v\mu(v,h,n)} -\varepsilon\delta_hH(v,h)\frac{\partial_h\mu(v,h,n)}{\partial_v\mu(v,h,n)}-\varepsilon\delta_nN(v,n)\frac{\partial_n\mu(v,h,n)}{\partial_v\mu(v,h,n)},  \\ 
	\dot{h} &= \delta_hH(v,h),  \\
	\dot{n} &= \delta_nN(v,n). 
	\end{aligned}
	\eqlab{3redslow}
	\end{align}
	The {reduced flow of \eqref{3redslow}} on $ {\mc{M}_2} $ is obtained by setting $ \varepsilon=0 $ in \eqref{3redslow}:
	\begin{align}
	\begin{aligned}
	0 &= \frac{m_\infty(v)-\mu(v,h,n)}{t_m(v)\partial_v\mu(v,h,n)}, \\ 
	\dot{h} &= \delta_hH(v,h),  \\
	\dot{n} &= \delta_nN(v,n); 
	\end{aligned}
	\eqlab{3redred}
	\end{align}
	then, with regard to the reduced flow on $ \mc{S}^a$, we have the following result:
	\begin{proposition}
		The {reduced flow of \eqref{3redslow}} on $ \mc{S}^a $ is smoothly topologically equivalent to the flow of
		\begin{align}
		\begin{aligned}
		\dot{v} &= \delta_h\partial_h  \left[V(v,m_\infty(v), h,n) \right]\lvert_{n = \nu(v,h)}H(v,h)+\delta_n\partial_n  \left[V(v,m_\infty(v), h,n) \right] \lvert_{n = \nu(v,h)}N(v,\nu(v,h)),   \\
		\dot{h} &=-\delta_h\partial_v  \left[V(v,m_\infty(v), h,n) \right]\lvert_{n = \nu(v,h)} H(v,h) 
		\end{aligned}
		\eqlab{em1red}
		\end{align}
		and to that of
		\begin{align}
		\begin{aligned}
		\dot{v} &= \delta_h\partial_h  \left[V(v,m_\infty(v), h,n) \right]\lvert_{h = \eta(v,n)}H(v,\eta(v,n))+\delta_n\partial_n  \left[V(v,m_\infty(v), h,n) \right]\lvert_{h = \eta(v,n)} N(v,n),  \\
		\dot{n} &=-\delta_n\partial_v  \left[V(v,m_\infty(v), h,n) \right]\lvert_{h = \eta(v,n)} N(v,n).
		\end{aligned}
		\eqlab{em1redn}
		\end{align}
		\proplab{em1red}
	\end{proposition}
	\begin{proof}
		Differentiating implicitly the algebraic constraint in \eqref{em1V} which defines {$\mathcal{M}_2$}, we 
		obtain
		\begin{align*}
		-\partial_v  \left[V(v,m_\infty(v), h,n) \right]\dot{v} &= \partial_h  \left[V(v,m_\infty(v), h,n) \right]\dot{h} +\partial_n  \left[V(v,m_\infty(v), h,n) \right] \dot{n}. 
		\end{align*}
		Making use of {the $h$- and $n$-equations of \eqref{3red}}, together with \eqref{Fn}, on $ {\mc{M}_2} $ we obtain
		\begin{align*}
		\begin{split}
		-\partial_v  \left[V(v,m_\infty(v), h,n) \right]\dot{v} &= \partial_h  \left[V(v,m_\infty(v), h,n) \right]\lvert_{n = \nu(v,h)}H(v,h)\\
		&\qquad \qquad +\delta_n\partial_n  \left[V(v,m_\infty(v), h,n) \right]\lvert_{n = \nu(v,h)} N(v,\nu(v,h)), 
		\end{split} \\
		\dot{h} &=\delta_h H(v,h).
		\end{align*}
		Rescaling time in the above by a factor of $-\lp \partial_v  \left[V(v,m_\infty(v), h,n) \right]\rp^{-1}$, which preserves the direction of time on $ \mc{S}^a $, we find \eqref{em1red}, as claimed. The equivalence to \eqref{em1redn} can be shown in a similar fashion.
	\end{proof}
	
	We remark that, although \eqref{em1red} and \eqref{em1redn} are equivalent on $ \mc{S}^a $, the regime where $ h $ is taken to be the slowest variable in \eqref{3redslow} is naturally studied in the context of \eqref{em1red}, which becomes a two-dimensional slow-fast system in the standard form of GSPT in that regime, as will become apparent in \secref{hslow}. Similarly, the $n$-slow regime is conveniently studied in \eqref{em1redn}; cf. \secref{nslow}. (By contrast, if $ n $ is taken to be the slow variable in \eqref{em1red}, respectively if $ h $ is taken to be the slow variable in \eqref{em1redn}, then the corresponding two-dimensional system is a slow-fast system in the non-standard form of GSPT \cite{wechselberger2020geometric}.)
	
	Naturally, and as is the case for \eqref{RW} also, the slow-fast formulation in Equation~\eqref{3redslow}  exhibits dynamics on two timescales when only $\varepsilon$ is assumed to be small and $ \delta_h, \delta_n = \mc{O}(1) $; correspondingly, there is then no separation of scales in the desingularised {reduced flow on $ \mc{S}^a $} in \eqref{em1red}. In the following, we consider the two regimes where either $ \tau_n= \mathcal{O}(1) $ and $ \tau_h $ is large in \eqref{hh} and, hence, $ \delta_h $ is small in \eqref{3redslow}; or $ \tau_h= \mathcal{O}(1) $ and $ \tau_n$ is large and, hence, $ \delta_n $ is small. The system in \eqref{3redslow} then exhibits dynamics on three distinct timescales, with \eqref{em1red} corresponding to the formulation on the ``intermediate'' timescale. In the following, we will in turn investigate the resulting mixed-mode dynamics in those two regimes. 
	
	\section{The $h$-slow regime}
	\seclab{hslow}
	
	In this section, we consider the regime where the variable $h$ is the slowest variable in  \eqref{hhmod-fast}, which is
	realised for $\delta_h>0$ sufficiently small and $\delta_n=1$, that is, when $ \tau_h $ is large and $ \tau_n =\mc{O}(1)$ in the original HH model, Equation~\eqref{original}. In that regime, {the three-dimensional reduction} in \eqref{3redfast} on {$ {\mathcal{M}_{1\gamma}} $} is a slow-fast system in the standard form of GSPT which features similar geometric and dynamical properties to the extended prototypical example introduced in \cite{kaklamanos2022bifurcations}.
	
	In particular, we will classify the mixed-mode dynamics of Equation~\eqref{hhmod-fast} with $ \gamma,\varepsilon $, and $ \delta_h>0 $ small {in terms of} the (rescaled) applied current $\bar I$, by applying the analysis outlined in \cite{kaklamanos2022bifurcations} to the three-dimensional reduction in \eqref{3redfast}. We will show that, in the parameter regime defined in \eqref{parbars}, there exist values $ 0<\bar{I}_{h}^-<\bar{I}_h^a<\bar{I}_h^r<\bar{I}_{h}^+ $ of $\bar I$ that distinguish between the various types of oscillatory dynamics in \eqref{hhmod-fast} for $ \gamma $, $\varepsilon$, and $\delta_h$ positive and sufficiently small and $ \delta_n = \mc{O}(1) $. The resulting classification is illustrated in \figref{ibif}, with the corresponding mixed-mode oscillatory dynamics as shown in \figref{HStimeseries}. While such dynamics was previously documented in \cite{doi2001complex}, the underlying geometric mechanisms that generate these various firing patterns, and the transitions between them, were not emphasised in the context of GSPT. Correspondingly, the $ \bar{I} $-values for which these transitions occur were not identified.
	
	\begin{figure}[ht!]
		\begin{subfigure}[b]{1\textwidth}
			\centering
			\includegraphics[scale = 0.15]{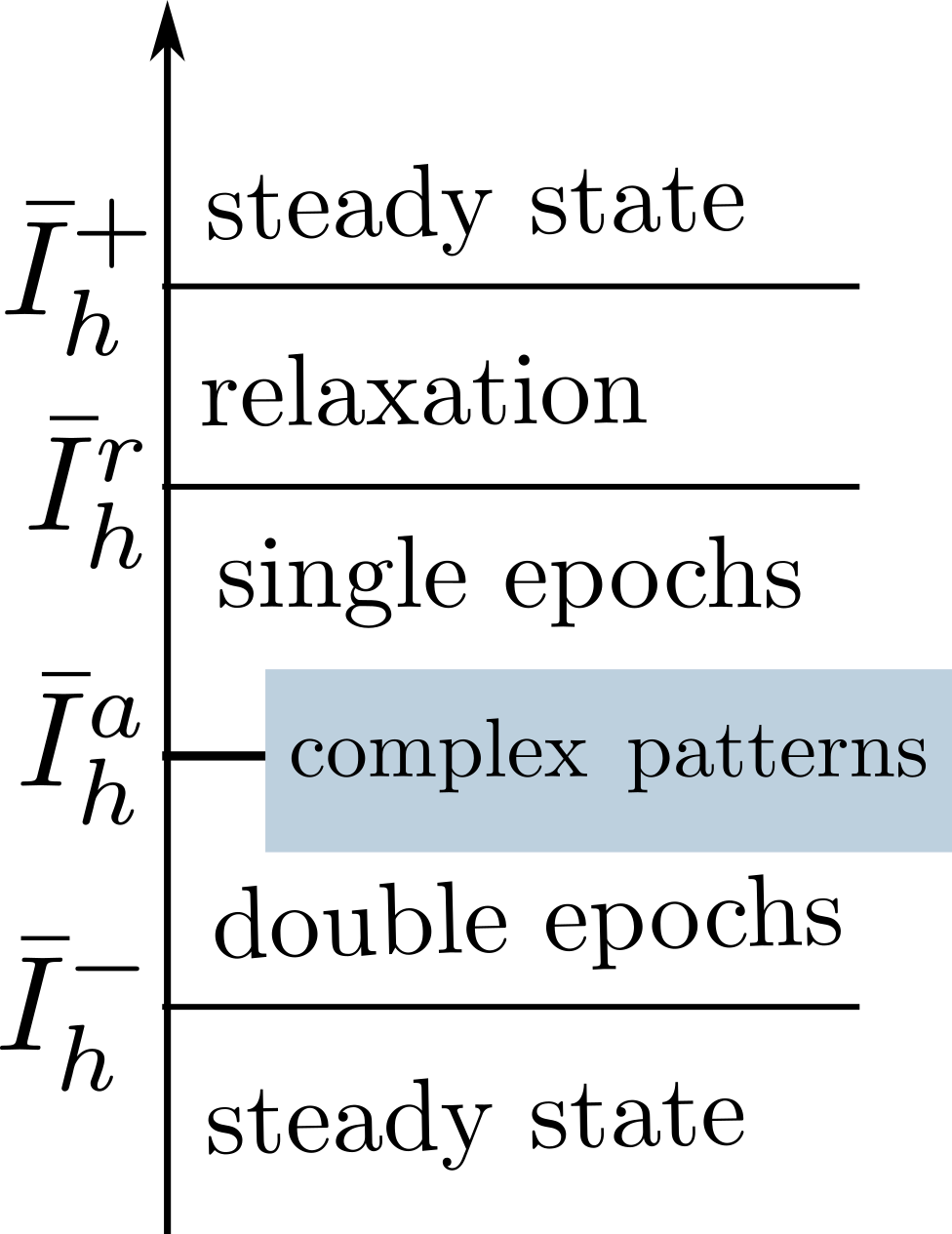}
		\end{subfigure}
		\caption{Bifurcations of mixed-mode oscillatory dynamics {in terms of} the parameter $ \bar{I} $ in the $ h $-slow regime, for $ \varepsilon$ and $\delta_h$ positive and sufficiently small in Equation~\eqref{3red}. {Complex oscillatory patterns are observed in a subset of $(\bar{I}_h^-,\bar{I}_h^r)$ around $\bar{I}_h^a$.}}
		\figlab{ibif}
	\end{figure}
	
	Specifically, the values $ \bar{I}_{h}^-$ and $ \bar{I}_{h}^+ $ {distinguish between oscillatory dynamics and steady-state behaviour}. The value $ \bar I_h^a $ separates oscillatory trajectories with double epochs of slow dynamics from those with single epochs under the perturbed flow of \eqref{3redslow}, with $\varepsilon$ and $\delta_h$ sufficiently small; see panels (a) and (b), as well as (d) and (e), of \figref{HStimeseries}, respectively. We remark that complex mixed-mode patterns featuring double epochs of slow dynamics separated by LAOs are observed during the transition from double-epoch MMOs to those with single epochs, i.e., for $ \bar I $-values close to $ \bar I_h^a $; see \figref{HStimeseries}(c). Moreover, there exists an $\bar I$-value $\bar I_h^r$ which separates MMOs with single epochs of SAOs from relaxation oscillation; see panels (e) and (f) of \figref{HStimeseries}, respectively. 
	
	\subsection{Singular geometry}
	
	In the singular limit of $ \varepsilon =0 $, with $ {\delta_h}>0 $ sufficiently small and $ \delta_n=1 $, the flow in Equation~\eqref{em1red} reads
		\begin{align}
		\begin{aligned}
		\begin{split}
		\dot{v} &= \partial_n\left[V(v,m_\infty(v), h,n) \right]\big\lvert_{n = \nu(v,h)} N(v,\nu(v,h))+\mc{O}\lp \delta_h\rp, 
		\end{split}  \\
		\dot{h} &=-\delta_h\partial_v  \left[V(v,m_\infty(v), h,n) \right]\big\lvert_{n = \nu(v,h)} H(v,h) 
		\end{aligned}
		\eqlab{em1hinter}
		\end{align}
	in the intermediate formulation, whereas on the 
	slow timescale {$ s_h = \delta_h t$}, we can write
	
		\begin{align}
		\begin{aligned}
		\begin{split}
		\delta_h\dot{v} &= \partial_n\left[V(v,m_\infty(v), h,n) \right]\big\lvert_{n = \nu(v,h)} N(v,\nu(v,h))+\mc{O}\lp \delta_h\rp, 
		\end{split}  \\
		\dot{h} &=-\partial_v  \left[V(v,m_\infty(v), h,n) \right]\big\lvert_{n = \nu(v,h)} H(v,h); 
		\end{aligned}
		\eqlab{em1hslow}
		\end{align}
	hence, we obtain a slow-fast system in the standard form of GSPT \cite{fenichel1979geometric}. 
	
	Setting $ \delta_h=0$ in Equation~\eqref{em1hinter} gives the one-dimensional layer problem
	\begin{align}
	\begin{aligned}
	\dot{v} &= \partial_n  \left[V(v,m_\infty(v), h,n) \right]\big\lvert_{n = \nu(v,h)} N(v,\nu(v,h)),  \\
	\dot{h} &=0. 
	\end{aligned}
	\eqlab{em1hlay}
	\end{align}
	Solutions of \eqref{em1hlay} are called the \textit{intermediate fibres} of {the three-dimensional reduction} in \eqref{3redfast} with $ h $ constant; {see \cite{kaklamanos2022bifurcations} for details}. {Equilibria of \eqref{em1hlay}, which are given by
		\begin{align*}
		N(v,\nu(v,h))=0,
		\end{align*}
		define the {critical} manifold $ \mathcal{M}_{h} $, which is the subset of ${\mc{M}_2}$ that can be expressed as
		\begin{align}
		\begin{aligned}
		n_\infty(v)-\nu(v,h)=0;
		\end{aligned} \eqlab{algconMh}
		\end{align}
		recall \eqref{N4} and \eqref{Fn}. We reiterate that $\partial_n  \left[V(v,m_\infty(v), h,n) \right]\neq 0$ by \eqref{parders}, and we remark that the algebraic constraint in \eqref{algconMh} is equivalent to
		\begin{align}
		V(v,m_\infty(v),h,n_\infty(v)) =0, \eqlab{acHV}
		\end{align}
		where $ V(v,m, h,n) $ is as defined in \eqref{V4}.}
	
	
	The manifold $ \mathcal{M}_h $  is normally hyperbolic on the set $ \mathcal{H}$ where 
	\begin{align*}
	\partial_v  \left[V(v,m_\infty(v), h,n_\infty(v)) \right] \neq0,
	\end{align*}
	losing normal hyperbolicity on the set $ \mc{F}_{\mc{M}_h} =\mc{M}_h\backslash\mc{H} $ where \begin{align}
	\partial_v  \left[V(v,m_\infty(v), h,n_\infty(v)) \right] =0.
	\end{align}
	{Computations show that
		\begin{align*}
		\mc{F}_{\mc{M}_h} = \begin{cases}
		\lb p_h^-, p_h^+\rb & \quad\tn{for } \bar{I} < \bar{I}_h^p, \\
		\emptyset & \quad\tn{for }\bar{I} > \bar{I}_h^p,
		\end{cases} 
		\end{align*}
		where 
		\begin{align} 
		\bar{I}_h^p \simeq \frac{120}{k_v g_{Na}}; \eqlab{Ip}
		\end{align}
		cf. \figref{hs-sing} and \figref{hs280}. When $\mathcal{M}_h$ admits two fold points $p_h^\mp$, these separate the normally hyperbolic portion $ \mc{H} $ of $\mc{M}_h$ into two outer  branches $ \mc{H}^{\mp}$, whereon $\partial_v  \left[V(v,m_\infty(v), h,n_\infty(v)) \right] <0$, and a middle branch $ \mc{H}^r$, whereon $ \partial_v  \left[V(v,m_\infty(v), h,n_\infty(v)) \right] >0$; see \figref{hs-rel}. }
	
	{Away from the tangential intersection of $\mc{L}^\mp$, $\mc{S}$ can be decomposed into two attracting sheets $\mc{S}^{a^\mp}$ that are separated by a repelling sheet $\mc{S}^{r}$, as shown in \figref{hs-sing}. Computations show that the points $ p_h^\mp $, when they exist, lie on $ \mathcal{S}^r $. In terms of the desingularised flow in \eqref{em1hinter}, the subsets $\mc{H}^\mp\cap\mc{S}^{a^\mp}$ are characterised as stable, while the branches $\mc{H}^\mp\cap\mc{S}^{r}$ are characterised as unstable; see \cite{kaklamanos2022bifurcations} for a detailed discussion in the context of a prototypical three-timescale system with similar geometric properties.}
	
	\begin{figure}[ht!]
		\centering
		\begin{subfigure}[b]{0.45\textwidth}
			\centering
			\includegraphics[scale = 0.35]{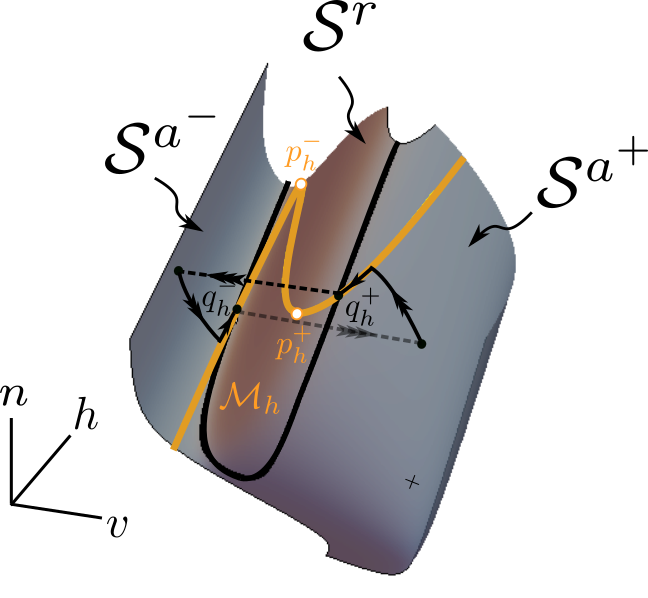}
			\caption{}
		\end{subfigure}
		~
		\begin{subfigure}[b]{0.45\textwidth}
			\centering
			\includegraphics[scale = 0.15]{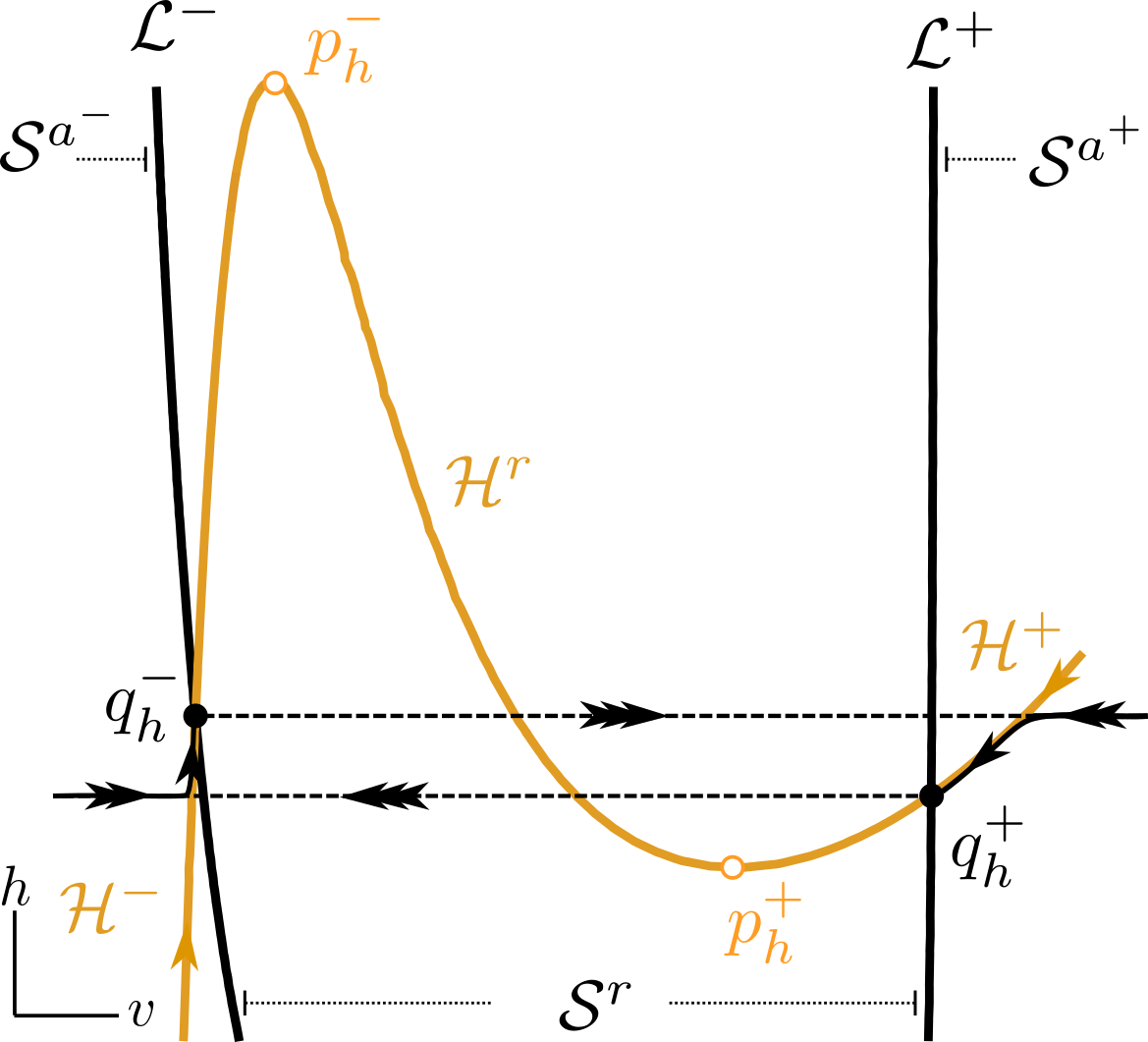}
			\caption{}
		\end{subfigure}
		\caption{Singular geometry of {the three-dimensional reduction} in \eqref{3redfast} for $\bar{I} = 20/(k_vg_{Na})$. In the double singular limit of $ \varepsilon=0=\delta_h $ in \eqref{3redfast}, the manifold $ \mathcal{H} $ consists of two outer branches $ \mathcal{H}^{\mp} $ which are attracting in $\mc{S}^{a^\mp}$, and separated by a repelling branch $ \mathcal{H}^r $. The folded singularities $ p_h^\mp $ of \eqref{3redfast} lie on $ \mc{S}^r $; the flow on $ {\mathcal{H}^{\mp} \cap \mc{S}^a}$ is directed towards $ q^\mp_h $. Singular cycles are obtained by concatenating fast, intermediate, and slow segments of {\eqref{3redlay}}, {the $v$-equation of \eqref{em1hlay}}, and \eqref{redMh}, respectively {-- note that intermediate segments (solid line; double arrowhead) overlie parts of fast segments (dashed line; triple arrowhead) in the projection in panel (b)}.}
		\figlab{hs-sing}
	\end{figure}
	
	In the double singular limit of $ \varepsilon=0=\delta_h $, the \textit{folded singularities} $ Q_h=\lb q_h^-,q_h^+\rb $ \cite{szmolyan2001canards}  of \eqref{3red} are given by 
	\begin{align*}
	q^\mp_h = \mathcal{M}_h\cap\mathcal{L}^\mp;
	\end{align*} 
	see \figref{hs-sing}. In the following, we will denote these singularities in coordinate form as
	\begin{align*}
	q^\mp_h = \lp v_{q_h}^{^\mp}, h_{q_h}^{^\mp}, n_{q_h}^{^\mp}\rp; 
	\end{align*}
	{the above coordinates are obtained by solving \eqref{em1}, \eqref{em1V}, and \eqref{acHV}.}
	\begin{remark}
		We emphasise that the properties of all geometric objects defined above -- i.e., of $ \mathcal{M}_h $, $ \mathcal{F}_{\mathcal{M}_h} $, and $ Q_h $ -- are dependent on the parameter $ I $ or, rather, on its rescaled counterpart $\bar I$.
		\remlab{remdep}
	\end{remark}
	
	\begin{remark}
		We remark that consideration of the reduced flow {of \eqref{3redslow}} on $\mc{S}^a$ in \eqref{em1redn}, instead of \eqref{em1red} with $\delta_h>0$ small, gives a slow-fast system in the non-standard form of GSPT; cf.~\figref{hs-ns} and see \cite{jelbart2020two,wechselberger2020geometric} for details.
		
		\begin{figure}[ht!]
			\centering
			\includegraphics[scale = 0.17]{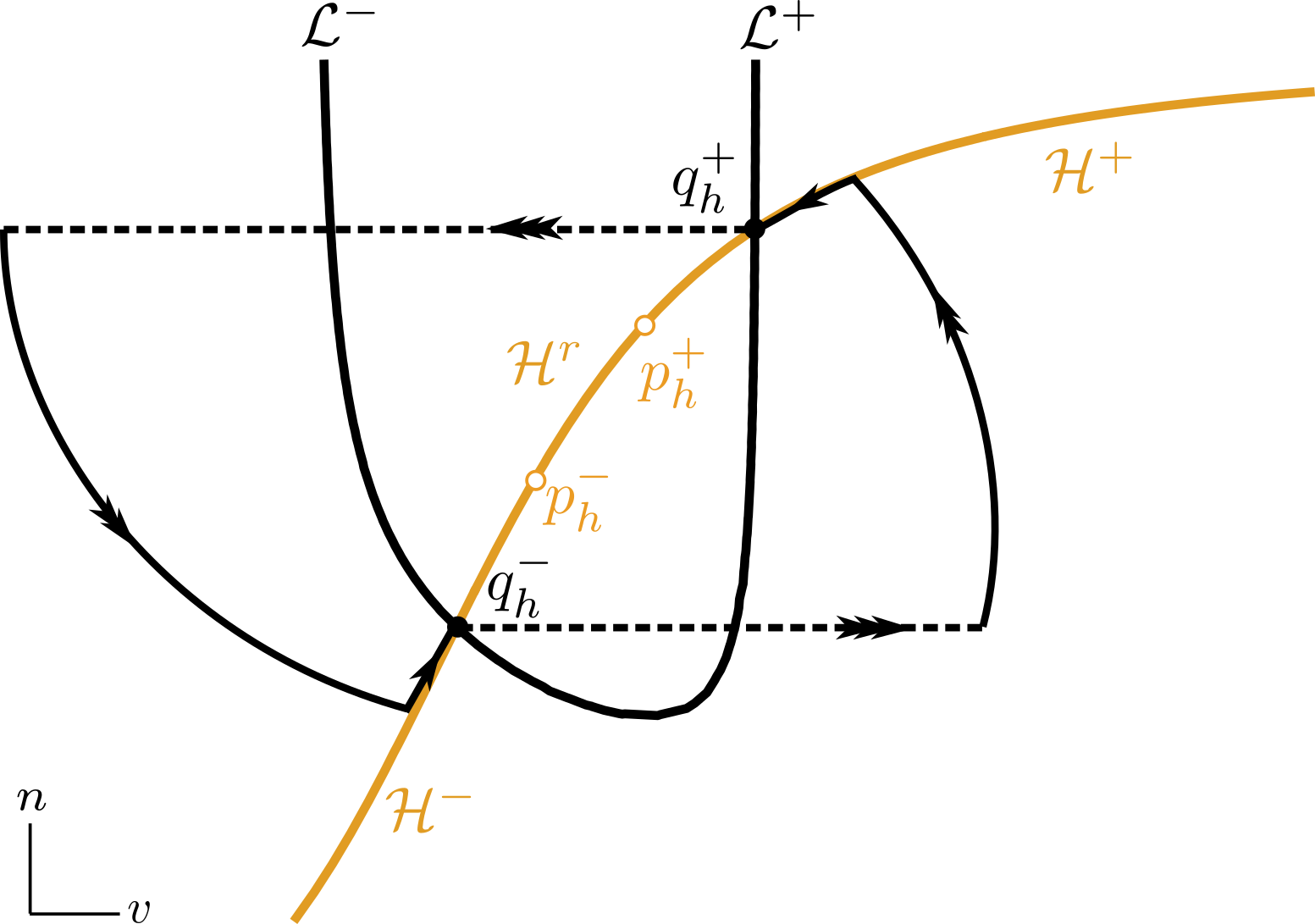}
			\caption{Critical manifold $\mc{M}_h$ in \eqref{em1redn} with $\delta_h>0$ small, which is a slow-fast system in the non-standard form of GSPT; {a singular cycle consisting of fast, intermediate, and slow segments is indicated for illustration of the resulting singular dynamics}.}
			\figlab{hs-ns}
		\end{figure}
		\remlab{hs-ns}
	\end{remark}
	
	Next, setting $ \delta_h=0$ in the slow formulation, Equation~\eqref{em1hslow}, we obtain the one-dimensional flow on $ \mathcal{M}_{h} $:
	
		\begin{align}
		\begin{aligned}
		0 &=  \partial_n  \left[V(v,m_\infty(v), h,n) \right]\lvert_{n = \nu(v,h)} N(v,\nu(v,h)),  \\
		\dot{h} &=-\partial_v  \left[V(v,m_\infty(v), h,n) \right]\lvert_{n = \nu(v,h)} H(v,h). 
		\end{aligned}
		\eqlab{em1hred}
		\end{align}

	Differentiating the algebraic constraint in \eqref{acHV} implicitly, via
	\begin{align*}
	\partial_v  \left[V(v,m_\infty(v), h,n_\infty(v) \right]\dot{v} &= {-}\partial_h   \left[V(v,m_\infty(v), h,n_\infty(v) \right]\dot{h},
	\end{align*}
	and using {the $h$-equation of \eqref{em1hred}} and \eqref{Gh}, we find the following expression for the flow on $ \mc{H}^{\mp} $:
	\begin{align}
	\dot{v} &=\frac{\partial_v  \left[V(v,m_\infty(v), h,n) \right]\partial_h  \left[V(v,m_\infty(v), h,n_\infty(v)) \right]}{\partial_v  \left[V(v,m_\infty(v), h,n_\infty(v) \right]}\bigg\lvert_{\{h = \eta(v,n{_\infty(v)}),n=n_\infty(v)\}} H(v,\eta(v,n{_\infty(v)})). \eqlab{redMh}
	\end{align}

	{When a true equilibrium exists in $\mc{S}^r$}, the {one-dimensional}  flow on $ \mc{H}^{\mp} $ is directed towards $ q_h^\mp $, respectively; cf. 	\figref{hs-sing} and \figref{hs280}. The resulting singular geometry is summarised in \figref{hs-sing}.  {Combining orbit segments from the layer flow of  Equation~\eqref{3redlay}$ _{\varepsilon=0} $, the intermediate fibres of \eqref{em1hlay} on $ {\mc{M}_2} $, and the {one-dimensional} dynamics of \eqref{redMh} on $ \mc{M}_h $, one can construct \textit{singular} or \textit{candidate trajectories}; closed singular trajectories are called \textit{singular cycles}.}
	
	We emphasise that, according to the above, {the three-dimensional reduction} in \eqref{3redfast}, with $ h $ the slow variable, has the {following properties, in analogy to those of the (three-timescale) extended prototypical example introduced in \cite{kaklamanos2022bifurcations} when no stable equilibrium of \eqref{3redred} exists on $ \mc{S}^a $}, namely:
	\begin{itemize}
		\item[P1.] The critical manifold $ {\mc{M}_2} $ is $ S $-shaped away from the tangential connection of $ \mc{L}^\mp $, with two attracting sheets $ \mc{S}^{a^\mp} $ separated by a repelling sheet $ \mc{S}^r $.
		\item[P2.] The manifold $ \mc{M}_h $ is $ S $-shaped, with two attracting branches $ \mc{H}^{\mp} $ separated by a repelling branch $ \mc{H}^r $; moreover, the fold points of $ \mc{M}_h $ lie on the repelling sheet $ \mc{S}^r $.
		\item[P3.] The {one-dimensional} flow on \eqref{redMh} $ \mc{H}^{\mp} $ is directed towards $ q^\mp_h $, respectively. 
	\end{itemize}
	
	{Correspondingly, we will conclude that the global mixed-mode dynamics of the three-timescale HH model, Equation~\eqref{hhmod-fast}, is fundamentally not very different to that of the three-timescale prototypical system studied in \cite{kaklamanos2022bifurcations}; {in the following, we will show how the various, qualitatively different firing patterns illustrated in \figref{HStimeseries} can be explained by the mechanisms described therein, resulting in the $\bar{I}$-dependent classification of \figref{ibif}.}}
	
	\subsection{Perturbed dynamics and MMOs}
	By standard GSPT \cite{fenichel1979geometric} and according to \cite{kaklamanos2022bifurcations}, there exist invariant manifolds $ \mathcal{S}^{a}_{\varepsilon\delta_h} $ and $\mathcal{S}^{r}_{\varepsilon\delta_h}$ that are diffeomorphic to their unperturbed, normally hyperbolic counterparts $ \mathcal{S}^{a}$ and $\mathcal{S}^{r}$, respectively, and that lie $ \mathcal{O}\lp\varepsilon,\delta_h\rp $-close to them in the Hausdorff distance, for $ \varepsilon$ and $\delta_h$ positive and sufficiently small. The perturbed manifolds $ \mathcal{S}^{a}_{\varepsilon\delta_h} $ and $\mathcal{S}^{r}_{\varepsilon\delta_h}$ are locally invariant under the flow of \eqref{3redfast}. Moreover, for $\varepsilon $ and $\delta_h$ sufficiently small, there exist invariant manifolds $ \mathcal{H}^{\mp}_{\varepsilon\delta_h} $ and $\mathcal{H}^{r}_{\varepsilon\delta_h}$ that are diffeomorphic, and $ \mathcal{O}\lp\delta_h\rp $-close in the Hausdorff distance, to their unperturbed, normally hyperbolic counterparts $ \mathcal{H}^{\mp}$ and $\mathcal{H}^{r}$, respectively. The perturbed branches $ \mathcal{H}^{\mp}_{\varepsilon\delta_h} $ and $\mathcal{H}^{r}_{\varepsilon\delta_h} $ are again locally invariant under the flow of \eqref{3redfast}. 
	
	{Away from the tangential connection between $\mc{L}^\mp$,} a trajectory initiated in a point $ (v,h,n)\in \mathcal{S}^{a}_{\varepsilon\delta_h} $ will follow the intermediate flow thereon until it is either attracted to $\mathcal{H}^{\mp}_{\varepsilon\delta_h} $ or until it reaches the vicinity of the fold curve $ \mc{L}^{\mp} $. If trajectories reach the vicinity of $ \mc{L}^{\mp} $ away from the folded singularities $ q^\mp_h $, they ``jump'' to the opposite attracting sheet $\mathcal{S}^{a^\pm}_{\varepsilon\delta_h}$ following the fast flow of \eqref{3redfast}; see \cite{wechselberger2005existence,szmolyan2001canards}. On the other hand, if trajectories are attracted to appropriate subregions of $ \mathcal{H}^{\mp}_{\varepsilon\delta_h} $ or to the vicinity of $ q^-_h $, they undergo {slow dynamics}; details can be found in \cite{kaklamanos2022bifurcations}. We include a brief discussion of SAO-generating mechanisms, as well as a comparison between \eqref{RW} and \eqref{3red} in that regard, in Appendix~\ref{partpert}. 
	
	Therefore, orbits for the perturbed flow of Equation~\eqref{3red}, with $ \varepsilon $ and $\delta_h$ positive and  sufficiently small, can be constructed by combining perturbations of fast, intermediate, and slow segments of singular trajectories. In \cite{kaklamanos2022bifurcations}, we showed that the qualitative characteristics of the resulting MMO trajectories will depend on the properties of the flow on $\mathcal{S}^{a}_{\varepsilon\delta_h}$ and $\mathcal{H}^{\mp}_{\varepsilon\delta_h} $, as well as on the location of the folded singularities $ q^\mp_h $ relative to one another.
	
	{{In the following, we demonstrate how the various $ \bar{I} $-values that distinguish between qualitatively different firing patterns in \eqref{3redfast} with $ \varepsilon,\delta_h >0$ small, as shown in \figref{ibif}, are obtained.}}

	\subsection{Onset and cessation of oscillatory dynamics}
	
	In the singular limit of $\varepsilon=0=\delta_h$ in \eqref{3redfast}, the flow on $\mc{H}$, as defined in Equation~\eqref{redMh}, has a stable equilibrium point on $\mc{S}^{r}$ for $\bar{I}\in (\bar{I}_{h}^-, \bar{I}_{h}^+)$, where we numerically obtain
	\begin{align} 
	\bar{I}_{h}^- \simeq \frac{4.8}{k_v g_{Na}}\quad\text{and}\quad 
	\bar{I}_{h}^+ \simeq \frac{280}{k_v g_{Na}}; \eqlab{Imp}
	\end{align} 
	cf. \figref{hs280}. Moreover, we observe that $ \bar{I}_{h}^-<\bar{I}_h^p<\bar{I}_{h}^+$; recall \eqref{Ip}. Therefore, the {one-dimensional} flow in \eqref{redMh} has a stable equilibrium point on $\mc{H}^{-}$ for $\bar{I}<\bar{I}_{h}^-$, and on the unique normally hyperbolic branch $\mc{H}$ for $\bar{I}> \bar{I}_{h}^+$, see again \figref{hs280}; correspondingly, that equilibrium crosses $ \mc{L}^- $ for both $ \bar{I} = \bar{I}_{h}^- $ and $ \bar{I} = \bar{I}_{h}^+ $.
	
	\begin{figure}[ht!]
		\centering
		\begin{subfigure}[b]{0.3\textwidth}
			\centering
			\includegraphics[scale = 0.22]{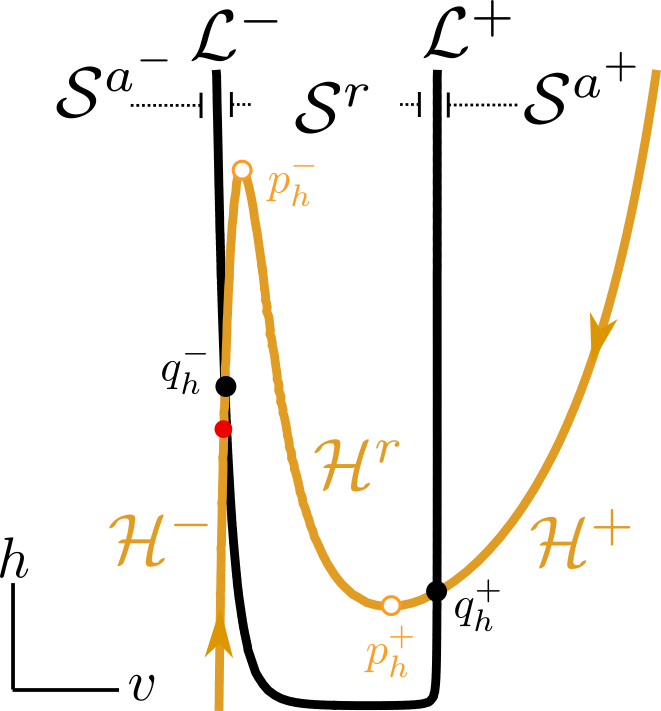}
			\caption{$\bar{I}<\bar{I}_{h}^-$}
		\end{subfigure}
		~
		\begin{subfigure}[b]{0.3\textwidth}
			\centering
			\includegraphics[scale = 0.22]{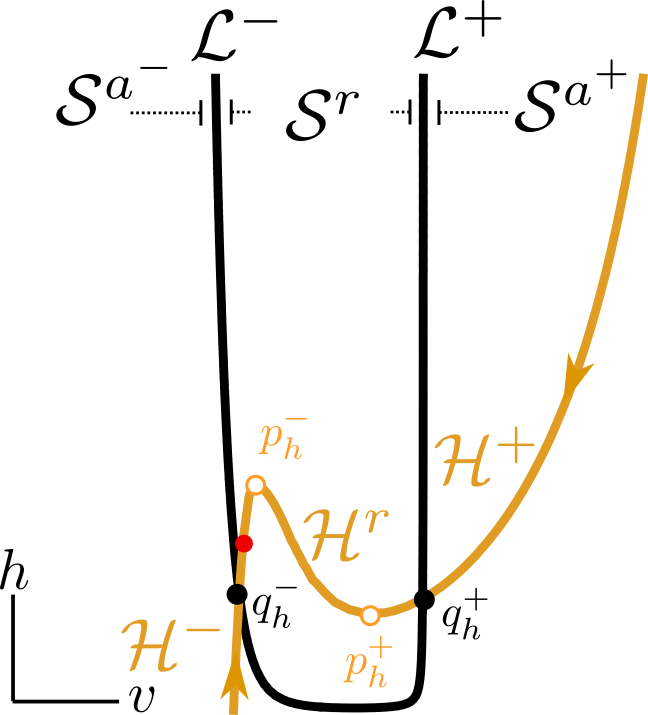}
			\caption{$\bar{I}\in (\bar{I}_{h}^-,\bar{I}_{h}^+)$}
		\end{subfigure}
		~
		\begin{subfigure}[b]{0.3\textwidth}
			\centering
			\includegraphics[scale = 0.22]{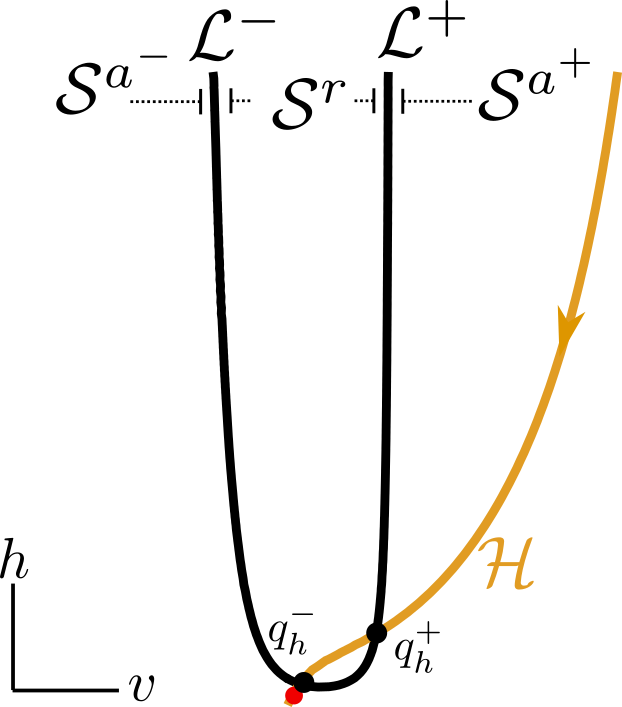}
			\caption{$\bar{I}>\bar{I}_{h}^+$}
		\end{subfigure}
		\caption{A stable equilibrium of the {one-dimensional} flow in \eqref{redMh} (red dot) lies in the attracting portion $\mc{S}^a$ of ${\mc{M}_2}$ for $\bar{I}<\bar{I}_{h}^-$ and $\bar{I}>\bar{I}_{h}^+$, as seen in panels (a) and (c), respectively, and in the repelling portion $\mc{S}^r$ for $\bar{I}\in (\bar{I}_{h}^-,\bar{I}_{h}^+)$; cf.~panel (b). As $ \bar{I} $ increases, the distance between $ \mc{L}^- $ and $ \mc{L}^+ $ decreases, with the leftmost branch of $ \mc{H} $ ``moving downward''. For $\bar{I} >\bar{I}_h^p$, the manifold $\mc{M}_h$ has no fold points; it then consists of a unique branch $\mc{H}$ which is attracting on $\mc{S}^a$ and repelling on $\mc{S}^r$.}
		\figlab{hs280}
	\end{figure}

	{For $\varepsilon,\delta_h>0$ sufficiently small, {the three-dimensional reduction} in \eqref{3redfast} undergoes Hopf bifurcations for $\bar{I}$-values that are $\mc{O}(\varepsilon,\delta_h)$-close to $\bar{I}_h^\mp$, which implies that the ``full" four-dimensional system, Equation~\eqref{hhmod-fast}, undergoes Hopf bifurcations for $\bar{I}$-values  that are $\mc{O}(\gamma,\varepsilon,\delta_h)$-close to $\bar{I}_h^\mp$. The Hopf bifurcation of Equation~\eqref{hhmod-fast} near $\bar{I}_h^-$ is subcritical, while the one close to $\bar{I}_h^+$ is supercritical \cite{doi2001complex,rubin2007giant}.
		It then follows that, for fixed $\bar{I}\in (\bar{I}_h^-, \bar{I}_h^+)$, there exist $\gamma,\varepsilon,\delta_h>0$ sufficiently small such that
		Equation~\eqref{3redfast} -- and, by extension, Equation~\eqref{hhmod-fast} -- features global mixed-mode dynamics.} By numerical sweeping, we find that for $\gamma=0.083$, $\varepsilon = 0.1$, $\delta_h = 0.025$, and $\delta_n=1$, with the values of the remaining parameters as in \eqref{parbars}, the onset of oscillatory dynamics in \eqref{hhmod-fast} occurs at approximately $\bar{I} \simeq 8.1/(k_v g_{Na})$, while the cessation of oscillatory dynamics is observed at approximately $\bar{I} \simeq 272/(k_v g_{Na})$; these values are close to $\bar{I}_{h}^-$ and $\bar{I}_{h}^+$, respectively, as expected. 
	
	{We emphasise that the global mixed-mode dynamics of \eqref{3redfast} is determined mainly by the position of $q^\mp_h$ relative to one another, {as will become apparent in the following,} and that the qualitative properties of that dynamics are not affected by the fold points $p_h^\mp$ being destroyed as $\bar{I}$ increases.}

	\subsection{Double epochs {of slow dynamics}}
	{At the beginning of this section, we introduced the notion of singular cycles, and we showed a particular example in \figref{hs-sing}. However, the configuration of the singular cycle illustrated therein is not the only possible one. Namely, depending on the location of the folded singularities $ q_h^\mp $ relative to one another, singular cycles can pass through only one or through both of them, as shown in \figref{hs-rel}{;  by \remref{remdep}, the corresponding geometry depends on the parameter $\bar{I}$.} We characterise the pair of folded singularities $ q_h^\mp $ based on whether they are ``connected'' by singular cycles or not, by slightly rephrasing the corresponding definition in \cite{kaklamanos2022bifurcations}, as follows. 
		
		\begin{definition}[\cite{kaklamanos2022bifurcations}]
			\textcolor{white}{white text}
			\begin{enumerate}
				\item The folded singularities $ q_h^\mp $ are \textit{``orbitally connected''} if they are contained in the same singular cycle, with the latter containing two slow segments, one on $ \mathcal{H}^{-} $ and one on $ \mathcal{H}^{+} $, i.e., if $ h_{q_h}^->h_{q_h}^+ $; cf.~\figref{hs-rel}(a).
				\item The folded singularities $ q_h^\mp $ are \textit{``orbitally remote''} if a singular cycle that passes through $ q^-_h $ does not pass through $ q^+_h $, i.e., if $ h_{q_h}^-<h_{q_h}^+ $; cf.~\figref{hs-rel}(c).  		
				\item The folded singularities $ q_h^\mp $ are \textit{``orbitally aligned''} if they are neither orbitally remote nor orbitally connected, i.e., if $ h_{q_h}^-=h_{q_h}^+ $; cf.~\figref{hs-rel}(b). 
			\end{enumerate}
			\defnlab{orbconHH}
		\end{definition} 
		
		\begin{figure}[ht!]
			\centering
			\begin{subfigure}[b]{0.3\textwidth}
				\centering
				\includegraphics[scale = 0.15]{pics/hh/hslow265.png}
				\caption{ $ {I} < \bar{I}_h^a $}
			\end{subfigure}
			\begin{subfigure}[b]{0.3\textwidth}
				\centering
				\includegraphics[scale = 0.15]{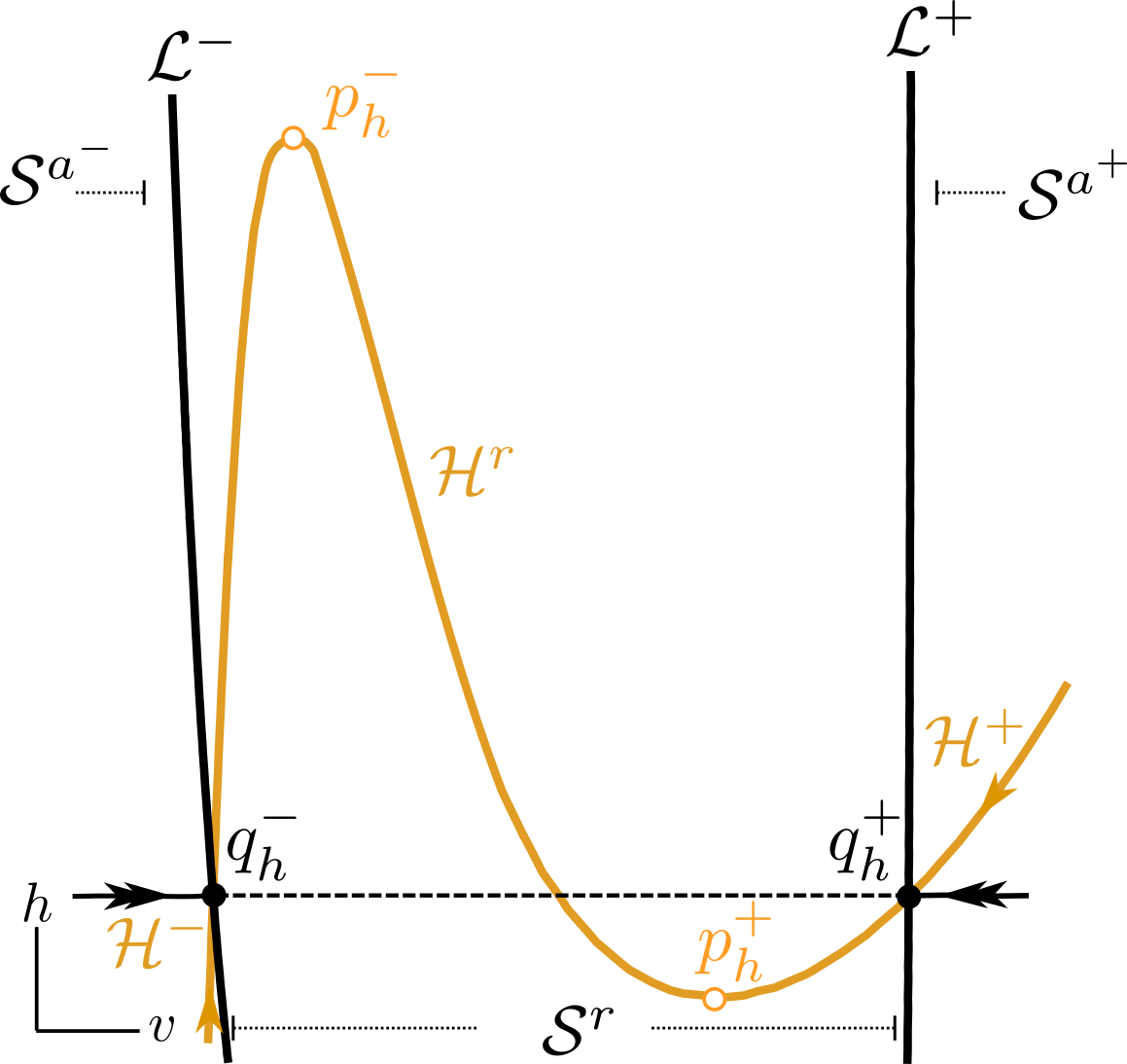}
				\caption{ $ \bar{I} = \bar{I}_h^a  $}
			\end{subfigure}
			\begin{subfigure}[b]{0.3\textwidth}
				\centering
				\includegraphics[scale = 0.15]{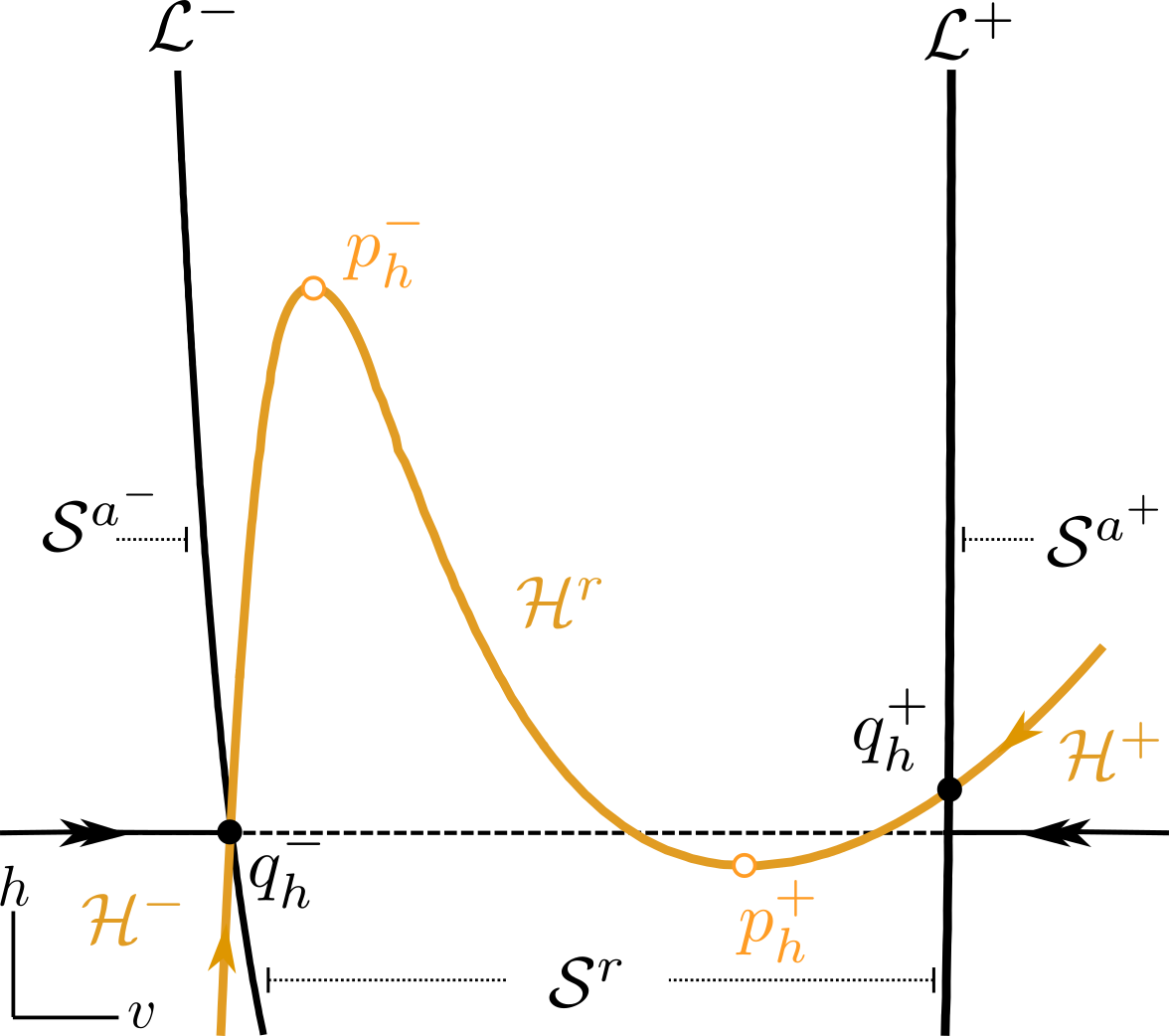}
				\caption{$ \bar{I} > \bar{I}_h^a $}
			\end{subfigure}
			\caption{Orbital connection, or lack thereof, in the limit of $ \varepsilon=0=\delta_h $, in accordance with \defnref{orbconHH}. (a) For $ \bar{I} < \bar{I}_h^a  $, the folded singularities $ q^\mp_h  $ are orbitally connected, in the sense that there exists a singular cycle that passes through both of them and that also contains slow segments on both $ \mc{H}^{\mp} $. (b) For $ \bar{I} = \bar{I}_h^a $, $ q^\mp_h $ are orbitally aligned, in that there exists a singular cycle that passes through both of them but that contains no slow segments on either $ \mc{H}^{\mp}$ for $ \varepsilon=0=\delta_h $. (c) For $ \bar{I} > \bar{I}_h^a  $, $ q^\mp_h $ are orbitally remote, in that a singular cycle that passes through one of the folded singularities does not pass through the other. } 
			\figlab{hs-rel}
		\end{figure}
		
		
		Due to the $\bar I$-dependence of the location of $q_h^\mp$, the classification in \defnref{orbconHH} hence encodes the position of the folded singularities $q_h^\mp$ relative to one another {in terms of} the rescaled applied current $ \bar{I} $. It follows that, in the parameter regime given by \eqref{parbars}, there exists a unique value $\bar{I}_h^a$  such that the singularities $q_h^-$ and $q_h^+$ are orbitally aligned. Numerically, that value is approximately found to equal
		\begin{align*}
		\bar{I}_h^a \simeq \frac{26.49}{{k_v g_{Na}}}; 
		\end{align*} 
		recall \eqref{parbars}. Correspondingly, the folded singularities $ q^\mp_h $ are orbitally connected for $\bar I<\bar I^a_h$ and orbitally remote for $\bar I>\bar I^a_h$. {We have the following result:}

		{\begin{proposition}
				Fix $\bar{I}\in(\bar{I}_h^-,\bar{I}_h^a)$ and let $\delta_n = 1$. Then, there exist $\gamma_0, \varepsilon_0$, and $\delta_{0}$ positive and sufficiently small such that the HH equations in \eqref{hhmod-fast} feature MMOs with double epochs of {slow dynamics} for all $(\gamma, \varepsilon, \delta_h)\in(0,\gamma_0)\times(0, \varepsilon_0)\times(0, \delta_{0})$.
				\proplab{hsdouble}
		\end{proposition}}
		
		\begin{proof}
			By \thmref{3reduction}, there exists $\gamma_0$ sufficiently small such that, for $\gamma\in(0,\gamma_0)$, the dynamics of \eqref{hhmod-fast} near its attracting slow manifold ${\mathcal{M}_{1\gamma}}$ is captured by the dynamics of {the three-dimensional reduction} in \eqref{3redfast}. Given Assumptions~P1 through P3, the proof that, for every fixed $\bar{I}\in (\bar{I}_h^-, \bar{I}_h^a)$ -- i.e., when the folded singularities $q^\mp_h$ are orbitally connected -- there exist $ \varepsilon_0,\delta_{0}>0 $ sufficiently small such that the three-dimensional system in \eqref{3redfast} exhibits MMOs with double epochs of {slow dynamics} for all $(\varepsilon, \delta_h)\in(0, \varepsilon_0)\times(0, \delta_{0})$, is similar to the proof of \cite[Theorem 1]{kaklamanos2022bifurcations}. {(Qualitatively, the statement follows from the fact that trajectories of \eqref{3redfast} are attracted to the vicinity of both $\mc{H}^\mp$ for $\varepsilon,\delta_h>0$ sufficiently small, and that they therefore feature double epochs of {slow dynamics}; cf. \figref{hs-sing} and panels (a) and (b) of \figref{HStimeseries}.)} 
		\end{proof}
		
		\propref{hsdouble} states that for every fixed $\bar{I}\in(\bar{I}_h^-,\bar{I}_h^a)$, there exist $\gamma, \varepsilon, \delta_h>0$ sufficiently small such that {the three-dimensional reduction} in \eqref{3redfast} -- and, by extension, Equation~\eqref{hhmod-fast} -- exhibits MMOs with double epochs of {slow dynamics}. Conversely, for fixed $ \varepsilon,\delta_h>0 $ sufficiently small, \eqref{3redfast} admits MMOs with double epochs of {slow dynamics} for $\bar{I}\in (\bar{I}_h^-+\mc{O}(\varepsilon, \delta_h), \bar{I}_h^a+\mc{O}(\varepsilon, \delta_h))$. By extension, for $ \gamma,\varepsilon, $ and $ \delta_{h}>0 $ sufficiently small, the mixed-mode dynamics of the full four-dimensional system in \eqref{hhmod-fast} is characterised by double epochs of {slow dynamics} for $\bar{I}\in (\bar{I}_q^-+\mc{O}(\gamma,\varepsilon, \delta_h), \bar{I}_h^a+\mc{O}(\gamma,\varepsilon, \delta_h))$; cf. panels (a) and (b) of \figref{HStimeseries} and \figref{HSphase}.
		
		By numerical sweeping, we find that for $\gamma=0.083$, $\varepsilon = 0.1$, and $\delta_h = 0.025$, with the values of the remaining parameters as in \eqref{parbars}, \eqref{hhmod-fast} features MMOs with double epochs of {slow dynamics} without LAOs between them, as illustrated in panels (a) and (b) of \figref{HStimeseries}, until approximately $\bar{I} \simeq 23.09/(k_v g_{Na})$, which is close to $\bar{I}_h^a$. 
		
		\begin{figure}[ht!]
			\centering
			\begin{subfigure}[b]{0.45\textwidth}
				\centering
				\includegraphics[scale = 0.45]{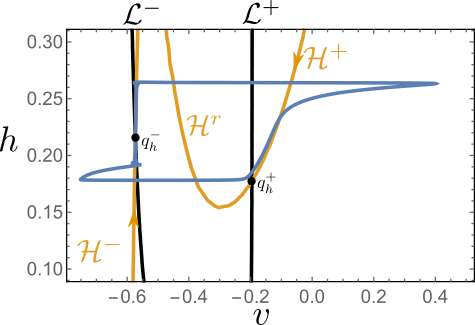}
				\caption{$ {I} = 20.051 $}
			\end{subfigure}
			~
			\begin{subfigure}[b]{0.45\textwidth}
				\centering
				\includegraphics[scale = 0.45]{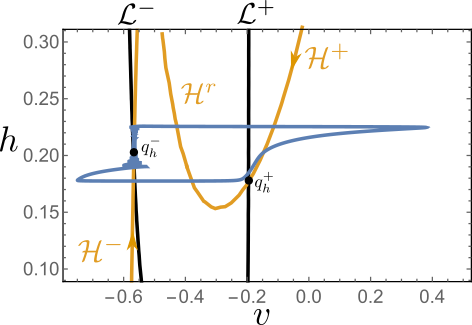}
				\caption{$ {I} = 23.051 $}
			\end{subfigure}
			\\
			\centering
			\begin{subfigure}[b]{0.45\textwidth}
				\centering
				\includegraphics[scale = 0.45]{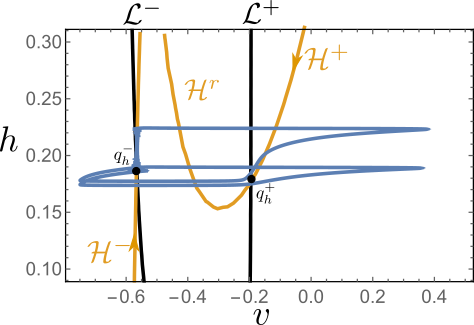}
				\caption{$ {I} = 23.5 $}
			\end{subfigure}
			~
			\begin{subfigure}[b]{0.45\textwidth}
				\centering
				\includegraphics[scale = 0.45]{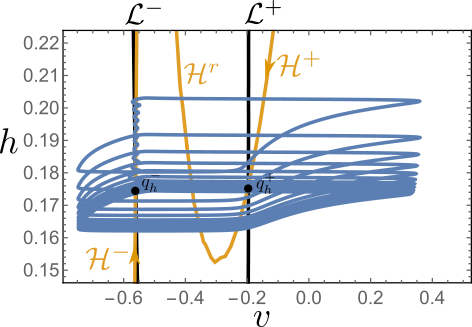}
				\caption{$ {I} = 26.03452346 $}
			\end{subfigure}
			\\
			\centering
			\begin{subfigure}[b]{0.45\textwidth}
				\centering
				\includegraphics[scale = 0.45]{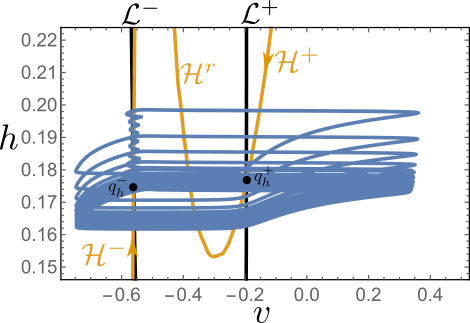}
				\caption{$ {I} = 26.1209956 $}
			\end{subfigure}
			~
			\begin{subfigure}[b]{0.45\textwidth}
				\centering
				\includegraphics[scale = 0.45]{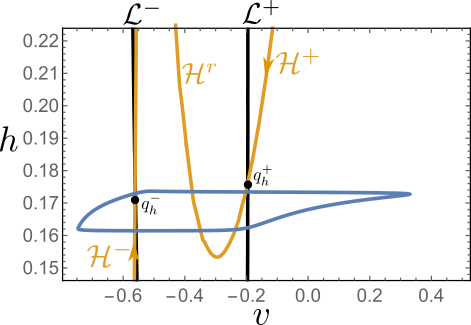}
				\caption{$ {I} = 26.2 $}
			\end{subfigure}
			\caption{Singular geometry of Equation~\eqref{em1hinter} which underlies the qualitative properties of the time series, and the associated MMO trajectories, illustrated  in \figref{HStimeseries}: as $\bar{I}$ is increased, the orbitally connected folded singularities $q_h^\mp$ become orbitally remote. Therefore, Equation~\eqref{hhmod-fast} with $ \varepsilon =0.0083 $, $\tau_h=40$,  and varying values of $ {I} =k_vg_{Na}\bar{I}$, transitions from mixed-mode dynamics with double epochs {of slow dynamics} for $ \bar{I}<\bar{I}_h^a $ (panels (a) and (b)) to single epochs thereof for $ \bar{I}>\bar{I}_h^a $ (panels (d) and (e)), via {``exotic"} MMO trajectories when $ \bar{I}$ is close to $\bar{I}_h^a $  (panel (c)). Finally, for $ \bar{I}>\bar{I}_h^r $, \eqref{hhmod-fast} features relaxation oscillation (panel (f)).}
			\figlab{HSphase}
		\end{figure}
		
		\subsection{{{``Transitional''} MMOs}}
		In \cite{kaklamanos2022bifurcations}, it is clearly stated that MMO trajectories of {{the three-dimensional reduction} in \eqref{3redfast} with $\varepsilon, \delta_h>0$ small } are not, strictly speaking, perturbations of individual singular cycles; therefore, other (and more complicated) phenomena might occur {in terms of} the values of $ \varepsilon $ and $ \delta_h $ in {\eqref{hhmod-fast} and \eqref{3redfast}}. In particular, during the transition from oscillation with double epochs of {slow dynamics} to single-epoch MMOs, i.e., for $ \bar{I} $-values close to the value $\bar I_h^a$ for which the folded singularities $ q^\mp_h $ are orbitally aligned, Equation~\eqref{3redfast} exhibits {``exotic"} mixed-mode dynamics upon variation of $ \bar{I} $ in which SAO epochs ``above" and ``below" are separated by LAOs, as illustrated in \figref{HStimeseries}(c). {We refer to those MMO trajectories as ``transitional".} 
		
		{Correspondingly, in \figref{HSphase}(c), a trajectory that leaves the vicinity of $\mc{L}^+$, jumping without having interacted with $\mc{H}^+$ on $\mc{S}^{a^+}$, experiences delayed loss of stability along $\mc{H}^-$ on $\mc{S}^{a^-}$ before connecting to a segment which is attracted to $\mc{H}^+$ on $\mc{S}^{a^+}$. After jumping from the vicinity of $\mc{L}^+$, the trajectory does not interact with $\mc{H}^-$ on $\mc{S}^{a^-}$; a jump occurs at $\mc{L}^-$, and the cycle starts anew. As the exit points of trajectories from vicinities of $\mc{H}^{\mp}$ vary with the corresponding entry points, cf. Appendix~\ref{partpert}, trajectories may or may not reach $\mc{H}^{\mp}$ after jumping from $\mc{L}^{\pm}$ when the distance between the folded singularities $ q_h^\mp$ in the $h$-direction is comparable to the displacement in the $h$-direction near $\mc{S}^{a^\mp}$ in \eqref{3redslow}. We quantify that displacement  more precisely in the following subsection, referring the reader to \cite{kaklamanos2022bifurcations} for details.}

		We remark that the above is not surprising, as it is explicitly stated in \cite{kaklamanos2022bifurcations} that for a fixed ``remote'' geometry in {the three-dimensional reduction} in \eqref{3redfast}, one can find $\varepsilon$ and $\delta_h$ sufficiently small such that MMOs with single epochs of {slow dynamics} occur. In our example here, $\varepsilon$ and $\delta_h$ are fixed and $\bar{I}$ is varied, which means that a geometric configuration arises which allows for the phenomenon described above to occur, i.e., when $\varepsilon$ and $\delta_h$ are not sufficiently small in relation to these values of $\bar{I}$, in spite of the folded singularities $q^\mp_h$ being remote. 
		
		By numerical sweeping, we find that for $\gamma=0.083$, $\varepsilon = 0.1$, and $\delta_h = 0.025$, with the values of the remaining parameters as in \eqref{parbars}, Equation~\eqref{hhmod-fast} features such MMOs with double epochs {of slow dynamics} interspersed with LAOs until approximately $\bar{I} \simeq 26/(k_v g_{Na})$, which is close to $\bar{I}_h^a$, as expected.

\subsection{{Single epochs {of slow dynamics} and relaxation}}
		
		
		{Next, we discuss MMOs with single epochs of {slow dynamics} and relaxation oscillation in the case where the folded singularities $q_{h}^\mp$ of {{the three-dimensional reduction} in \eqref{3redfast}} are remote.}
		
		{We first describe how single-epoch MMO trajectories can be constructed for $\bar{I}>\bar{I}_h^a$}. Consider the limit of $ \varepsilon=0 $ in \eqref{em1red}, with $ \delta_h>0 $ small. Eliminating time, away from $ \mc{H} $ we can then approximate the flow on $ \mc{S}^{a^\mp} $ by
		\begin{align*}
		\frac{\tn{d}h}{\tn{d}v} = -\delta_h\frac{\partial_v  \left[V(v,m_\infty(v), h,n) \right]\lvert_{n = \nu(v,h)} H(v,h)}{\partial_n  \left[V(v,m_\infty(v), h,n) \right]\lvert_{n = \nu(v,h)} N(v,\nu(v,h))},
		\end{align*}
		where we have omitted $ \mc{O}\lp \delta_h\rp $-terms in the denominator; recall that {$\partial_n \left[V(v,m_\infty(v),h,n)\right]<0$, by \eqref{parders}}. 
		
		We now denote by $P(\mc{L}^\mp)\in \mc{S}^{a^\pm}$ the projection of $\mc{L}^\mp$ onto $\mc{S}^{a^\pm}$ along the fast fibres of \eqref{3redlay}, and we define the points
		\begin{align*}
		(v_{\tn{min}}, h_{q_h}^-, n_{\tn{min}}) := P(\mc{L}^+) \cap \lb h = h_{q_h}^-\rb\quad\text{and}\quad 
		(v_{\tn{max}}, h_{q_h}^-, n_{\tn{max}}) := P(\mc{L}^-) \cap \lb h = h_{q_h}^-\rb,
		\end{align*}
		as illustrated in \figref{umax}. A trajectory that leaves the vicinity of $ \mc{L}^- $ at a point with $ h $-coordinate close to $ h_{q_h}^- $ returns to a point with $h$-coordinate $ \hat{h} $ in that vicinity after one large excursion, where
		\begin{align}
		\begin{aligned}
		\hat{h}(h,\bar{I}) = h +\delta_h \lp \mc{G}^-_h(h, \bar{I})+\mc{G}^+_h(h, \bar{I})\rp+\mc{O}(\delta_h^2),
		\end{aligned}\eqlab{hret}
		\end{align}
		and where we denote
		\begin{align}
	\begin{aligned}
	\mc{G}^-_h(h, \bar{I}) &:= -\int_{v_{\tn{min}}}^{v_{q^-_h} } \frac{\partial_v  \left[V(v,m_\infty(v), h,n) \right]\lvert_{n = \nu(v,h)} H(v,h)}{\partial_n  \left[V(v,m_\infty(v), h,n) \right]\lvert_{n = \nu(v,h)} N(v,\nu(v,h))}\tn{d}v\quad\text{and} \\
	\mc{G}^+_h(h, \bar{I}) &:= -\int_{v_{\tn{max}}}^{v_{q^+_h} } \frac{\partial_v  \left[V(v,m_\infty(v), h,n) \right]\lvert_{n = \nu(v,h)} H(v,h)}{\partial_n  \left[V(v,m_\infty(v), h,n) \right]\lvert_{n = \nu(v,h)} N(v,\nu(v,h))}\tn{d}v.
	\end{aligned}
	\eqlab{GRints}
	\end{align}

		\begin{figure}[ht!]
			\centering
			\begin{subfigure}[b]{0.45\textwidth}
				\centering
				\includegraphics[scale = 0.37]{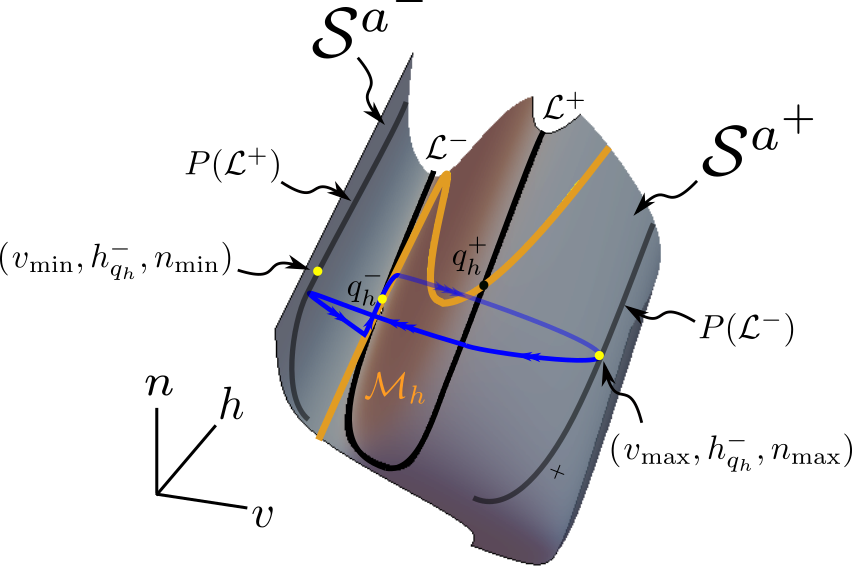}
				\caption{}
			\end{subfigure}
			\qquad
			\begin{subfigure}[b]{0.45\textwidth}
				\centering
				\includegraphics[scale = 0.37]{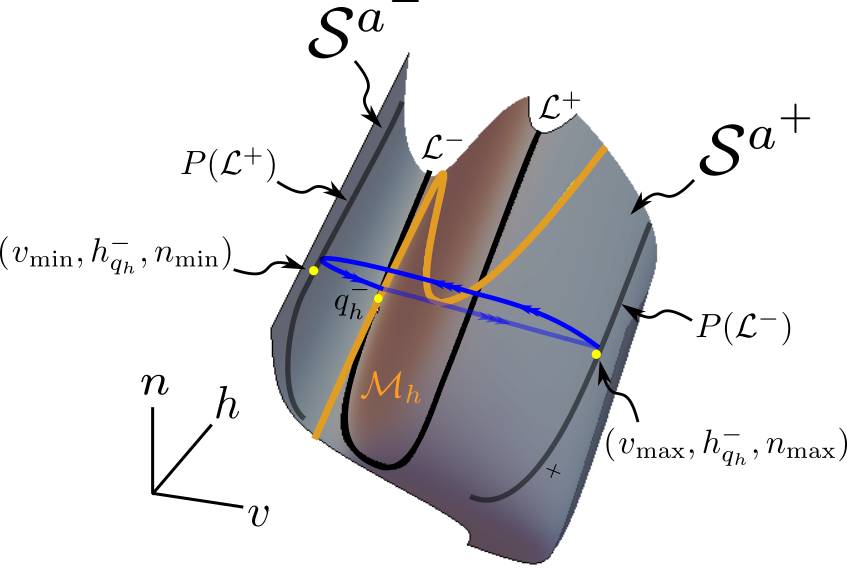}
				\caption{}
			\end{subfigure}
			\caption{{Schematic illustration of the points $(v_{\tn{min}}, h_{q_h}^-, n_{\tn{min}})$ and $(v_{\tn{max}}, h_{q_h}^-, n_{\tn{max}})$ on ${\mc{M}_2}$ in the two scenarios where (a) an MMO trajectory or (b) relaxation oscillation emerges.}}
			\figlab{umax}
		\end{figure}
		
		{\begin{remark}
				Referring back to \remref{nernst}, it is now apparent that, in order to apply GSPT, one should consider a sufficiently large compact subset of $\mc{M}_2$ which contains $P(\mc{L^\mp})$.
		\end{remark}}
		Equation~\eqref{hret} encodes the displacement in the $h$-direction after a large excursion, i.e., after the trajectory has been attracted to $\mathcal{S}^{a^+}$ and then again repelled away therefrom in the vicinity of $\mathcal{L}^+$ before its return to $\mathcal{S}^{a^-}$ near $\mathcal{L}^-$; see \cite{krupa2008mixed,kaklamanos2022bifurcations} for details. We define the displacement function as
		\begin{align*}
		\Delta_h(h,\bar{I}) = \mc{G}^-_h(h, \bar{I})+\mc{G}^+_h(h, \bar{I}).
		\end{align*}
		{While the integrals in \eqref{GRints} cannot be solved analytically, calculations using computer algebra software show that}
		\begin{align}
		\begin{aligned}
		\Delta_h({h_{q_h}^-},\bar{I})
		\begin{cases}
		<0 &\quad\tn{if }\bar{I}\in (\bar{I}_h^a,\bar{I}_h^r), \\
		=0 &\quad\tn{if }\bar{I}=\bar{I}_h^r, \\
		>0 &\quad\tn{if }\bar{I}>\bar{I}_h^r,
		\end{cases}
		\end{aligned}
		\eqlab{hbalance}
		\end{align}
		where {we numerically obtain}
		\begin{align*}
		\bar{I}_h^r \simeq \frac{29.705}{{k_v g_{Na}}}.
		\end{align*} 
		We then have the following result for the ``full", four-dimensional system in \eqref{hhmod-fast}.
		
		{\begin{proposition}
				Fix $\bar{I}\in(\bar{I}_h^a,\bar{I}_h^r)$ and let $\delta_n = 1$. Then, there exist $\gamma_0, \varepsilon_0$, and $\delta_{0}$ positive and sufficiently small such that the HH equations in \eqref{hhmod-fast} feature MMOs with single epochs of {slow dynamics} ``below'' for all $(\gamma, \varepsilon, \delta_h)\in(0,\gamma_0)\times(0, \varepsilon_0)\times(0, \delta_{0})$.
				\proplab{hssingle}
		\end{proposition}}
		
		\begin{proof}
			By \thmref{3reduction}, there exists $\gamma_0$ sufficiently small such that, for $\gamma\in(0,\gamma_0)$, the dynamics of \eqref{hhmod-fast} near its attracting slow manifold ${\mathcal{M}_{1\gamma}}$ is captured by the dynamics of \eqref{3redfast}. Moreover, if $\bar{I}\in (\bar{I}_h^a,\bar{I}_h^r)$ -- i.e., if {$\Delta_h({h_q^-},\bar{I})<0$}, by \eqref{hbalance} -- then  trajectories of \eqref{3redfast} with $\varepsilon=0$ and $\delta_h>0$ sufficiently small that leave $\mc{L}^-$ at $q_h^-$ will return at a point on $\mc{S}^{a^-}$ from which they will be attracted close to $\mc{H}^-$, i.e., in the funnel of $q_h^-$. As is well known \cite{wechselberger2005existence, brons2006mixed,krupa2008mixed,letson2017analysis,kaklamanos2022bifurcations}, in that case SAOs occur in the vicinity of $q_h^-$ in \eqref{3redfast} for $\varepsilon, \delta_h>0$ small. The above implies that, for fixed $\bar{I}\in (\bar{I}_h^a,\bar{I}_h^r)$, there exist $\varepsilon, \delta_h>0$ sufficiently small such that Equation~\eqref{3redfast} features MMOs with single epochs of {slow dynamics} \cite[Theorem 3]{kaklamanos2022bifurcations}. 
		\end{proof}
		
		\propref{hssingle} states that for every fixed $\bar{I}\in(\bar{I}_h^a,\bar{I}_h^r)$, there exist $\gamma, \varepsilon, \delta_h>0$ sufficiently small such that {the three-dimensional reduction} in \eqref{3redfast} -- and, by extension, Equation~\eqref{hhmod-fast} -- exhibits MMOs with single epochs of {slow dynamics}. Conversely, for fixed $ \gamma,\varepsilon,\delta_h>0$ small, the three-dimensional system in \eqref{3redfast} exhibits MMOs with single epochs of {slow dynamics} for $\bar{I}\in (\bar{I}_h^a+\mc{O}(\varepsilon,\delta_h),\bar{I}_h^r +\mc{O}(\delta_h, \delta_h\varepsilon\ln\varepsilon)) $; see \cite{szmolyan2004relaxation,kaklamanos2022bifurcations} for details. It then follows that the full four-dimensional system in \eqref{hhmod-fast} exhibits single-epoch MMOs for $\bar{I}\in (\bar{I}_h^a+\mc{O}(\gamma,\varepsilon,\delta_h), \bar{I}_h^r +\mc{O}(\gamma,\delta_h, \delta_h\varepsilon\ln\varepsilon)) $. By numerical sweeping, we find that for $\gamma=0.083$, $\varepsilon = 0.1$, and $\delta_h = 0.025$, with the remaining parameters defined as in \eqref{parbars}, Equation~\eqref{hhmod-fast} features MMOs with single epochs {of slow dynamics} until approximately $\bar{I} \simeq 26.127/(k_v g_{Na})$, which is consistent with the above error estimates.
		
		{We now proceed by showing the transition from mixed-mode dynamics to relaxation oscillation. To establish existence of the latter, we need to show that the assumptions of \cite[Theorem 4]{szmolyan2004relaxation} are satisfied; in the context of the HH equations, these assumptions are summarised in \cite[Section 3.1]{rubin2007giant} for the two-timescale case -- we reiterate that the three-timescale case which we study here is not considered therein. For a prototypical three-timescale system, the corresponding assumptions are outlined in \cite[Theorem 2]{kaklamanos2022bifurcations}.} {To that end, we rewrite \eqref{hret} as
		\begin{align*}
		\hat h(h,\bar{I})-h =  \delta_ h \Psi(h,\bar{I},\delta_h),
		\end{align*}
		where we define 
		\begin{align*}
		\Psi(h,\bar{I},\delta_h) = \Delta_h(h,\bar{I})+\mc{O}(\delta_h);
		\end{align*}
		recall \eqref{hret}. In particular, we then have that
		\begin{align}
		\Psi(h_{q_h}^-,\bar{I},0) = 0 = \Delta(h_{q_h}^-,\bar{I}), \eqlab{Del0}
		\end{align} 
		 by \eqref{hbalance}. Moreover, we use computer algebra software to evaluate the partial derivatives of $\Psi({h},\bar{I}, \delta_h)$ with respect to $h$ and $\bar{I}$: numerically, we calculate that
		\begin{align}
		\frac{\partial \Psi}{\partial h}(h_{q_h}^-,\bar{I}_h^r,0)\simeq -15.0377< 0 \eqlab{Delh}
		\end{align}
		and
		\begin{align}
		\frac{\partial \Psi}{\partial \bar{I}}(h_{q_h}^-,\bar{I}_h^r,0)\simeq -278.4470<0. \eqlab{DelI}
		\end{align}}
		
		{\begin{proposition}
				There exists $\bar{I}_h^*>\bar{I}_h^r$ such that, for $\bar{I}\in({I}_h^r, \bar{I}_h^*) $ and $\delta_n = 1$, there exist $\gamma_0, \varepsilon_0$, and $\delta_{0}$ positive and sufficiently small such that the HH equations in \eqref{hhmod-fast} feature relaxation oscillation for all $(\gamma, \varepsilon, \delta_h)\in(0,\gamma_0)\times(0, \varepsilon_0)\times(0, \delta_{0})$.
				\proplab{hsrelax}
		\end{proposition}}
		
		\begin{proof}
			
			By \thmref{3reduction}, there exists $\gamma_0>0$ sufficiently small such that, for $\gamma\in(0,\gamma_0)$, the dynamics of \eqref{hhmod-fast} near its attracting slow manifold ${\mathcal{M}_{1\gamma}}$ is captured by the dynamics of \eqref{3redfast}. The existence of stable relaxation oscillation in \eqref{3redfast} is established by showing that a set of assumptions introduced in \cite{szmolyan2004relaxation, rubin2007giant} is satisfied, similarly to the proof of \cite[Theorem 2]{kaklamanos2022bifurcations}.
			
			{Specifically, the first set of assumptions states that ${\mc{M}_2}$ is $S$-shaped with two fold curves $\mc{L}^\mp$, and that the reduced flow of \eqref{em1hinter} on ${\mc{M}_2}$ is transverse to both $\mc{L}^\mp$ and $P\lp \mc{L}^\pm\rp$. These assumptions are satisfied for Equation~\eqref{3redfast} away from $q_h^\mp$ and from the tangential connection of $\mc{L}^\mp$, as was also shown in \cite{rubin2007giant}.}
			
			{The other set of assumptions states that, in the singular limit of $\varepsilon=0$ with $\delta_h>0$ sufficiently small, there exists a hyperbolic and stable singular cycle. To that end, one shows that in this limit there exists a one-dimensional map  
				\begin{align*}
				\pi^-:\ \mc{L}^-\lvert_{h\in(h_{q_h}^-, h_{q_h}^+)}\to\mc{L}^-\lvert_{h\in(h_{q_h}^-, h_{q_h}^+)}, \   h\mapsto \hat{h} ,
				\end{align*}  which has a hyperbolic and stable fixed point. Here, the map $ \pi^- $ is given by \eqref{hret}}. 
			From the implicit function theorem, it follows that there exist $\bar{I}^*_h>\bar{I}_h^r$ and a function $ \tilde{h}(\bar{I},\delta_h)$, with $\tilde{h}(\bar{I}_h^r,0) = h_{q_h}^-$, such that
			$\Psi(\tilde{h}(\bar{I},\delta_h),\bar{I},\delta_h) =0$, i.e., the map $\pi^-$ has a fixed point for $\delta_h>0$ sufficiently small and $\bar{I}\in(\bar{I}_h^r,\bar{I}_h^*)$. The $h$-coordinate of this fixed point on $\mc{L}^-$ is given by $h=\tilde{h}(\bar{I},\delta_h)$;
			from \eqref{Delh} and \eqref{DelI}, it follows that, for fixed $\delta_h>0$ sufficiently small,  $\tilde{h}(\bar{I}, \delta_h)$ is a decreasing function of $\bar{I}$. {Specifically, we find that $\frac{\partial \tilde{h}}{\partial\bar{I}}(\bar{I}_h^r, \delta_h)\simeq -\frac{278.447}{15.0377}$, by \eqref{Delh} and \eqref{DelI}. On the other hand, we calculate that $\frac{\partial h_{q_h}^-}{\partial\bar{I}}(v_{q_h}^-, \bar{I}_h^r )=-279.416<\frac{\partial \tilde{h}}{\partial\bar{I}}(\bar{I}_h^r, \delta_h)$ for the $h$-coordinate $h_{q_h}^-=\eta(v_{q_h}^-, n_\infty(v_{q_h}^-))$ of the folded singularity $q_h^-$ in a first approximation. It follows that $\tilde{h}(\bar{I},\delta_h)>h_{q_h}^-$ for $\bar{I}\in(\bar{I}_h^r,\bar{I}_h^*)$, i.e., that the trajectory of \eqref{3redfast} corresponding to the fixed point of $\pi^-$ reaches the vicinity of $\mc{L}^-$ outside the funnel of $h_{q_h}^-$; see again \figref{umax} for an illustration of the underlying geometry.}
		
			Finally, from \eqref{hret}, we have that
			\begin{align*}
			\frac{\partial \hat{h}}{\partial h} (\tilde{h}(\bar{I},\delta_h),\bar{I}) = 1+\delta_h \frac{\partial \Psi}{\partial h}(\tilde{h}(\bar{I},\delta_h),\bar{I},\delta_h) \in (0,1)
			\end{align*}
			for all $\delta_h>0$ sufficiently small; therefore, the fixed point is attracting. It follows that, for $\varepsilon=0$ and $\delta_h>0$ sufficiently small, \eqref{3redfast} admits a stable singular cycle $ \Gamma_h $, while for $\varepsilon,\delta_h>0$ sufficiently small, the system in \eqref{3redfast} exhibits stable relaxation oscillation \cite{szmolyan2004relaxation, rubin2007giant, kaklamanos2022bifurcations}.
		\end{proof}
		
		\propref{hsrelax} states that there exists $\bar{I}^*_h>\bar{I}_h^r$ such that, if $\bar{I}\in (\bar{I}_h^r,\bar{I}^*_h)$, then there exist $\gamma, \varepsilon, \delta_h>0$ sufficiently small such that {the three-dimensional reduction} in \eqref{3redfast} -- and, by extension, Equation~\eqref{hhmod-fast} -- exhibits MMOs with single epochs of {slow dynamics}. Conversely, for fixed $ \gamma,\varepsilon,\delta_h>0$ small, the three-dimensional system in \eqref{3redfast} admits relaxation oscillation for $\bar{I}\in (\bar{I}_h^r+\mc{O}(\varepsilon,\delta_h),\bar{I}^*_h+\mc{O}(\varepsilon,\delta_h)) $; 
		see \cite{krupa2008mixed,kaklamanos2022bifurcations} for details. It then follows that the full four-dimensional system in \eqref{hhmod-fast} exhibits relaxation oscillation for  $\bar{I}\in (\bar{I}_h^r+\mc{O}(\gamma, \varepsilon,\delta_h),\bar{I}^*_h+\mc{O}(\gamma, \varepsilon,\delta_h)) $.

		
		Numerically, we observe that the transition from relaxation oscillation to steady state, i.e., the cessation of oscillatory dynamics when $\bar{I}$ is close to $\bar{I}_{h}^+$ -- recall Equation~\eqref{Imp} and the discussion that follows it -- happens via stable small-amplitude oscillations that are a product of the supercritical Hopf bifurcation near $\bar{I}_{h}^+$. We also remark that, in \cite{doi2001complex}, chaotic MMOs were documented during the transition from MMOs with single SAO epochs to relaxation; cf.~\figref{HStimeseries}(e). The study of the underlying mechanisms in the three-timescale regime considered here is included in plans for future work. 
		\begin{remark}
		    We note that the numerical values in \eqref{Delh} and \eqref{DelI} may seem large. However, we emphasise that the applied current $\bar{I}$ is itself a very small quantity due to the scaling in \eqref{scaling} \cite{rubin2007giant} which has been widely used in the study of the HH model, Equation~\eqref{original}. Therefore, the resulting variation in $\bar{I}$ is also small, in relation, for instance, to \eqref{DelI}. An investigation of potential alternative scalings of $h$ and $\bar{I}$ which might yield more moderate values for the quantities corresponding to \eqref{Delh} and \eqref{DelI} is left for future work. 
		    \remlab{scaling}
		\end{remark}
		
\section{The $n$-slow regime}
		\seclab{nslow}
		
		In this section, we consider the regime where the variable $n$ is the slowest variable in \eqref{hhmod-fast}, which is
		realised for $\delta_n>0$ sufficiently small and $\delta_h=1$, that is, when $ \tau_n $ is large and $ \tau_h =\mc{O}(1)$ in the original HH model, Equation~\eqref{original}. In that regime, {the three-dimensional reduction} in \eqref{3redfast} on $ {\mc{M}_2} $ is a slow-fast system in the standard form of GSPT. As the analysis of this regime is fundamentally similar to that of the $h$-slow regime presented in \secref{hslow}, our exposition here is deliberately more condensed.}
	
	\begin{figure}[h!]
		\centering
		\begin{subfigure}[b]{0.45\textwidth}
			\centering
			\includegraphics[scale = 0.5]{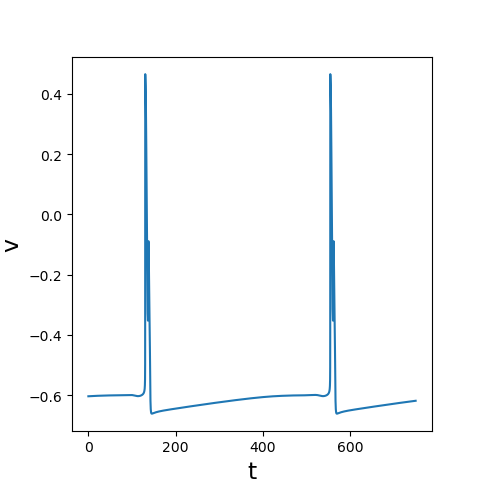}
			\caption{$ {I} = 9 $}
		\end{subfigure}
		~
		\centering
		\begin{subfigure}[b]{0.45\textwidth}
			\centering
			\includegraphics[scale = 0.5]{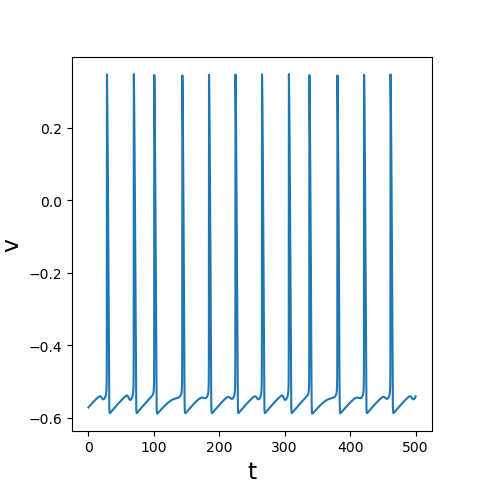}
			\caption{$ {I} = 64.5 $}
		\end{subfigure}
		\\
		\begin{subfigure}[b]{0.45\textwidth}
			\centering
			\includegraphics[scale = 0.5]{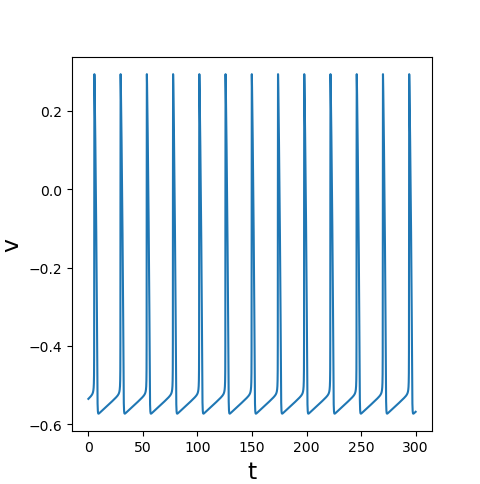}
			\caption{$ {I} = 90 $}
		\end{subfigure}
		~
		\begin{subfigure}[b]{0.45\textwidth}
			\centering
			\includegraphics[scale = 0.5]{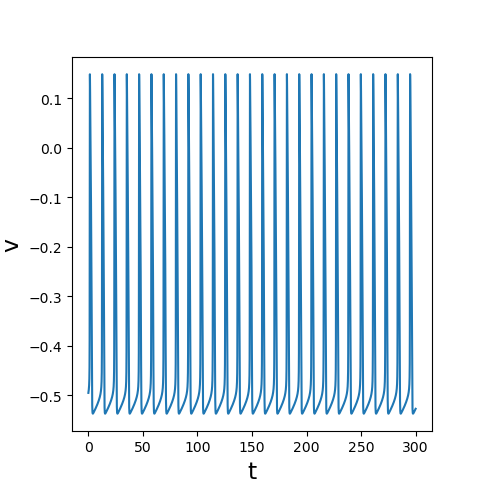}
			\caption{$ {I} = 150 $}
		\end{subfigure}
		\caption{Time series of the $ v $-variable in the dimensionless HH model, Equation~\eqref{hh}, for $ \varepsilon =0.0083 $, $\tau_n=100$, and varying values of $ {I} =k_vg_{Na}\bar{I}$: {panel (a) illustrates an MMO trajectory with double epochs of {slow dynamics}; panel (b) shows an MMO trajectory with a single epoch of {slow dynamics}. Finally, panels (c) and (d) display relaxation oscillation of different amplitudes.}}
		\figlab{NStimeseries}
	\end{figure}
	
	As was done in \secref{hslow}, we will again classify the mixed-mode dynamics of Equation~\eqref{hhmod-fast} with $ \gamma,\varepsilon$, and $\delta_n>0 $ small {in terms of} $\bar I$, by applying results of \cite{kaklamanos2022bifurcations} to {the three-dimensional reduction} in \eqref{3redfast}. We will show that, in the parameter regime defined in \eqref{parbars}, there exist values $ 0<\bar{I}_{n}^-<\bar{I}_n^a<\bar{I}_{n}^+ $ of $\bar I$ that distinguish between the various types of oscillatory dynamics in \eqref{hhmod-fast} for $ \gamma$, $\varepsilon$, and $\delta_n$ positive and sufficiently small and $ \delta_h = \mc{O}(1) $. The resulting classification is illustrated in \figref{ibifn}, with the corresponding mixed-mode oscillatory dynamics as shown in \figref{NStimeseries}. {The geometric configuration in phase space for $n$ slow is analogous to the one illustrated in \figref{fig11} for the $h$-slow case.}
	
	\begin{figure}[ht!]
		\begin{subfigure}[b]{1\textwidth}
			\centering
			\includegraphics[scale = 0.15]{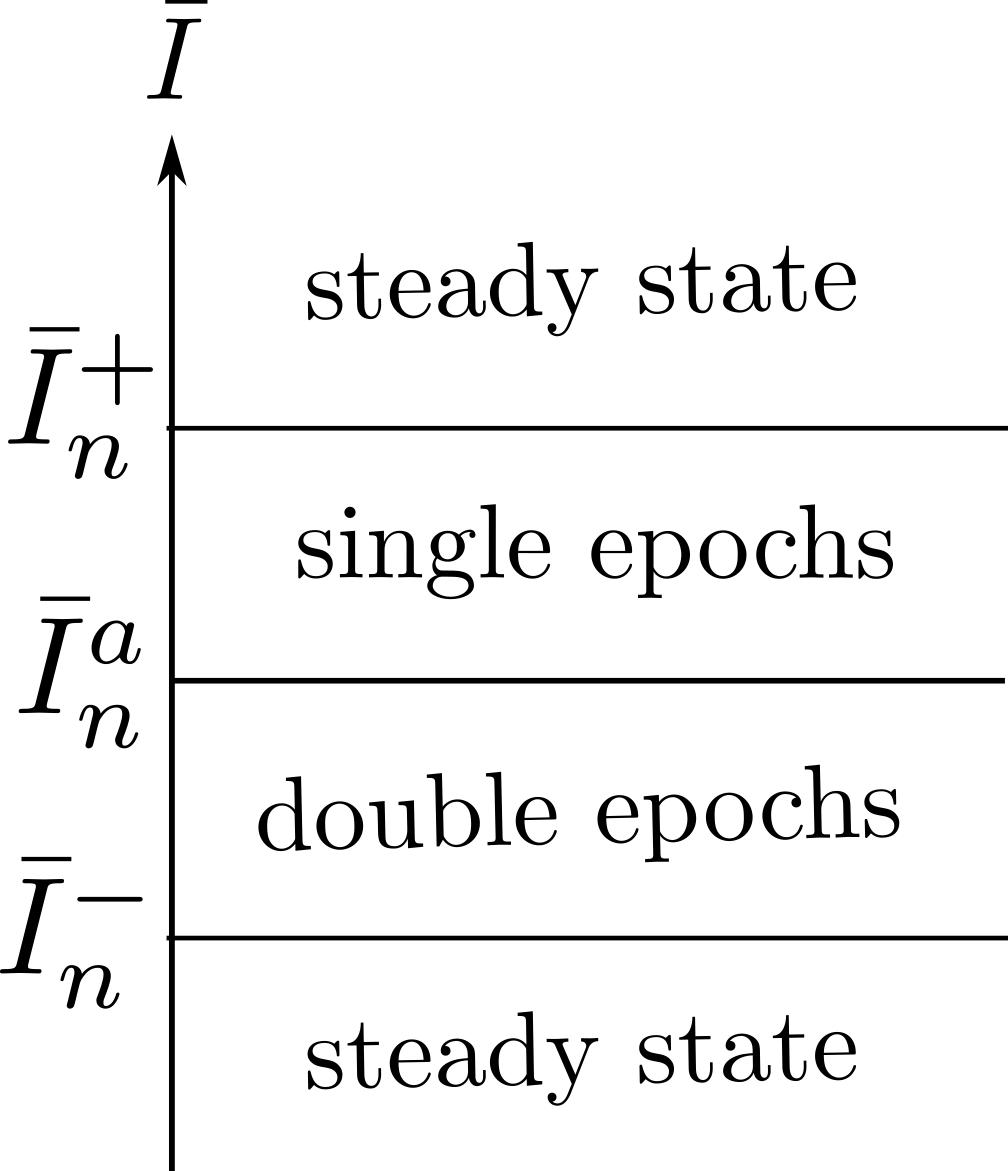}
		\end{subfigure}
		\caption{Bifurcations of oscillatory dynamics {in terms of} the parameter $ {\bar I} $ in the $ n $-slow regime, for $\varepsilon$ and $\delta_n$ positive and sufficiently small.}
		\figlab{ibifn}
	\end{figure}
	
	The values $ \bar{I}_n^-$ and $ \bar{I}_n^+$ again distinguish between oscillatory dynamics and steady-state behaviour. The value $ \bar I_n^a $ separates MMO trajectories with double epochs of {slow dynamics} from those with single epochs under the perturbed flow of {the three-dimensional reduction} in \eqref{3redfast}, with $\varepsilon$ and $\delta_n$ sufficiently small; see panel (a), as well as panels (d) and (e), of \figref{NStimeseries}, respectively.

	\subsection{Singular geometry}
	
	In the singular limit of $ \varepsilon =0 $, with $ {\delta_n}>0 $ sufficiently small and $ \delta_h=1 $, the {reduced flow of \eqref{3redslow}} on $ \mc{S}^a $ in Equation~\eqref{em1redn} reads
		\begin{align}
		\begin{aligned}
		\begin{split}
		{v}' &= \partial_h  \left[V(v,m_\infty(v), h,n) \right]\lvert_{h = \eta(v,n)} H(v,\eta(v,n))+\mc{O}\lp \delta_n\rp, 
		\end{split}  \\
		{n}' &=-\delta_n\partial_v  \left[V(v,m_\infty(v), h,n) \right]\lvert_{h = \eta(v,n)} N(v,n) 
		\end{aligned}
		\eqlab{em1ninter}
		\end{align}
	in the intermediate formulation, whereas on the 
	slow timescale {$ s_n = \delta_n t$}, we can write
	
		\begin{align}
		\begin{aligned}
		\begin{split}
		\delta_n\dot{v} &= \partial_h  \left[V(v,m_\infty(v), h,n) \right]\lvert_{h = \eta(v,n)} H(v,\eta(v,n))+\mc{O}\lp \delta_n\rp, 
		\end{split}  \\
		\dot{n} &=-\partial_v  \left[V(v,m_\infty(v), h,n) \right]\lvert_{h = \eta(v,n)} N(v,n); 
		\end{aligned}
		\eqlab{em1nslow}
		\end{align}
	hence, we again obtain a slow-fast system in the standard form of GSPT \cite{fenichel1979geometric}. 
	
	Setting $ \delta_n=0$ in Equation~\eqref{em1ninter} gives the one-dimensional layer problem
		\begin{align}
		\begin{aligned}
		{v}' &= \partial_h  \left[V(v,m_\infty(v), h,n) \right]\lvert_{h = \eta(v,n)} H(v,\eta(v,n)),  \\
		{n}' &=0. 
		\eqlab{em1nlay}
		\end{aligned}
		\end{align}
	Solutions of \eqref{em1nlay} yield \textit{intermediate fibres} with $ n $ constant. {Equilibria of \eqref{em1nlay}, which are given by
		\begin{align*}
		H(v,\eta(v,n))=0,
		\end{align*}
		define the {critical} manifold $ \mathcal{M}_{n} $, which is the subset of ${\mc{M}_2}$ that can be expressed as
		\begin{align}
		\begin{aligned}
		h_\infty(v)-\eta(v,n)=0;
		\end{aligned} \eqlab{algconMn}
		\end{align}
		recall \eqref{H4} and \eqref{Fn}. We reiterate that $\partial_h  \left[V(v,m_\infty(v), h,n) \right]\neq 0$ by \eqref{parders}, and we remark that the algebraic constraint in \eqref{algconMn} is equivalent to
		\begin{align}
		V(v,m_\infty(v),h_\infty(v),n) =0, \eqlab{acNV}
		\end{align}
		where $ V(v,m, h,n) $ is as defined in \eqref{V4}.}
	
	The manifold $ \mathcal{M}_n $  is normally hyperbolic on the set $ \mathcal{N}$ where 
	\begin{align*}
	\partial_v  \left[V(v,m_\infty(v), h_\infty(v),n) \right] \neq0,
	\end{align*}
	losing normal hyperbolicity on the set $ \mc{F}_{\mc{M}_n} =\mc{M}_n\backslash\mc{N} $ where 
	\begin{align*}
	\partial_v  \left[V(v,m_\infty(v), h_\infty(v),n) \right] =0.
	\end{align*}
	Computations show that
	\begin{align*}
	\mc{F}_{\mc{M}_n} = \begin{cases}
	\lb p_n^-, p_n^+\rb & \quad\tn{for } \bar{I} < \bar{I}_n^p, \\
	\emptyset & \quad\tn{for }\bar{I} > \bar{I}_n^p,
	\end{cases}
	\end{align*}
	where 
	\begin{align} 
	\bar{I}_n^p \simeq \frac{81.7}{k_v g_{Na}}; \eqlab{Ipn}
	\end{align}
	cf. \figref{ns-sing} and \figref{ns280}. When $\mathcal{M}_n$ admits two fold points $p_n^\mp$, these separate the normally hyperbolic portion $ \mc{N} $ of $\mc{M}_n$ into two outer  branches $ \mc{N}^{\mp}$, whereon $\partial_v  \left[V(v,m_\infty(v), h_\infty(v),n) \right] <0$, and a middle branch $ \mc{N}^r$, whereon $ \partial_v  \left[V(v,m_\infty(v), h_\infty(v),n) \right] >0$; see {\figref{ns-sing}}. 
	
	{Computations show that, in the parameter regime given by \eqref{parbars}, the points $ p_n^\mp $ lie on $ \mathcal{S}^r $, when they exist. In terms of the desingularised flow in \eqref{em1ninter}, the subsets $\mc{N}^\mp\cap\mc{S}^{a^\mp}$ are characterised as stable, while the branches $\mc{N}^\mp\cap\mc{S}^{r}$ are characterised as unstable; see \cite{kaklamanos2022bifurcations} for details.}
	
	\begin{figure}[ht!]
		\centering
		\begin{subfigure}[b]{0.45\textwidth}
			\centering
			\includegraphics[scale = 0.35]{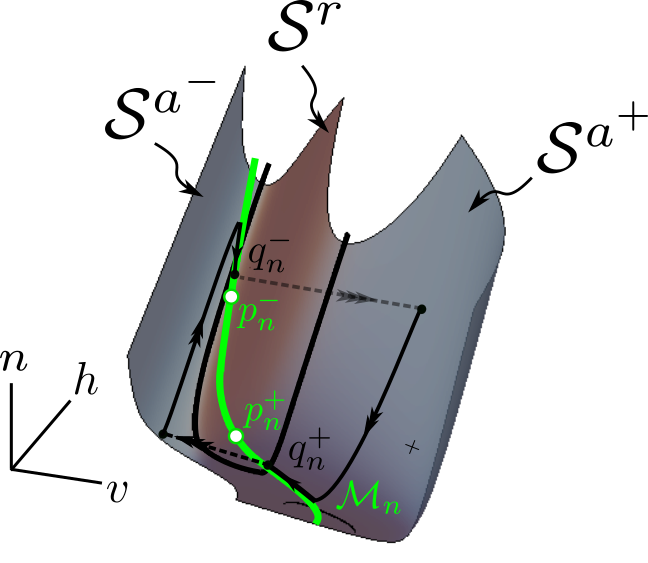}
		\caption{}
		\end{subfigure}
		~
		\begin{subfigure}[b]{0.45\textwidth}
			\centering
			\includegraphics[scale = 0.15]{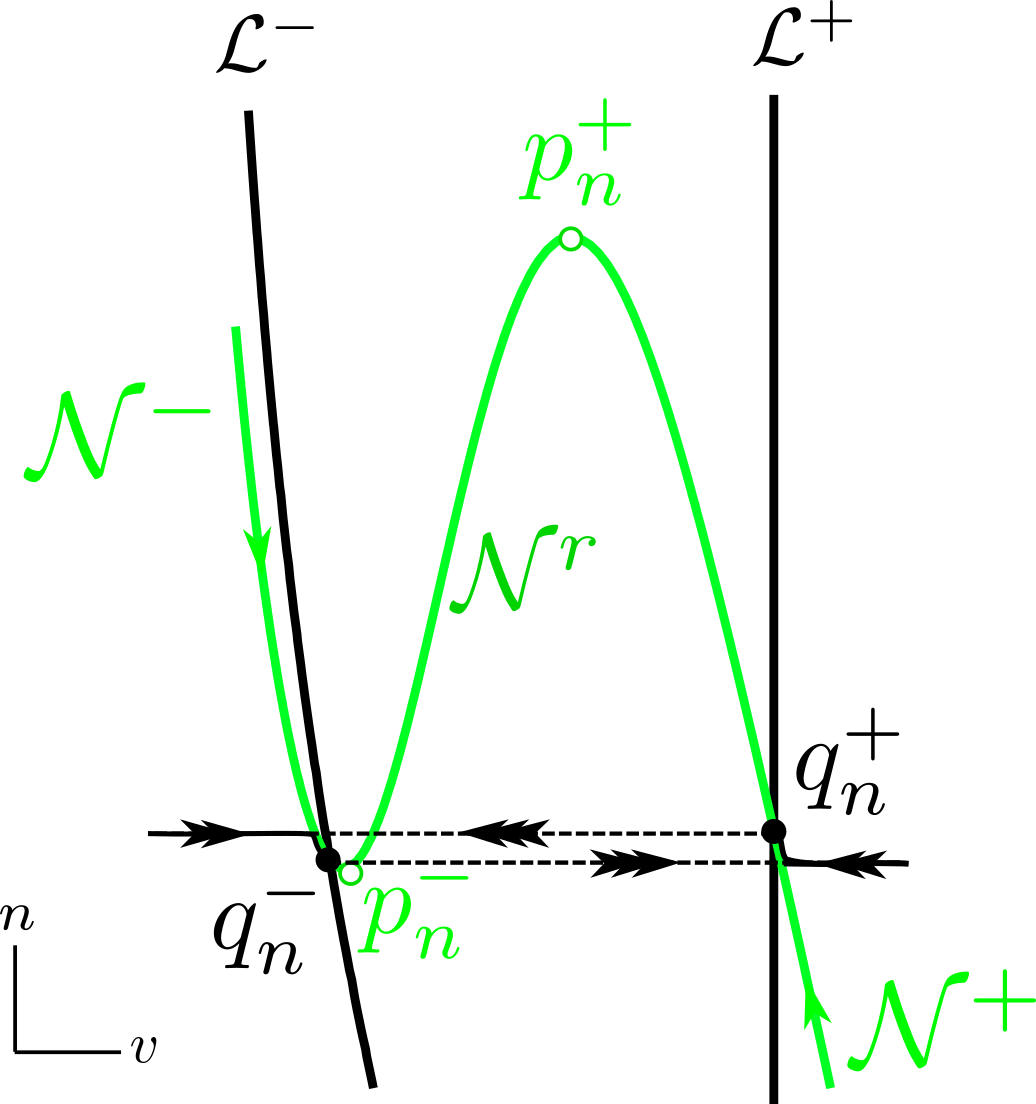}
			\caption{}
		\end{subfigure}
		\caption{Singular geometry of {the three-dimensional reduction} in \eqref{3redfast} for $\bar{I} = 7/(k_vg_{Na})$. (a) In the double singular limit of $ \varepsilon=0=\delta_n $, the manifold $ \mathcal{N} $ consists of two outer branches $ \mathcal{N}^{\mp} $ which are attracting in $\mc{S}^{a^\mp}$, and separated by a repelling branch $ \mathcal{N}^r $. The folded singularities $ p_n^\mp $ of \eqref{3redfast} lie on $ \mc{S}^r $; the flow on {$ \mathcal{N}^{\mp} \cap \mc{S}^a$} is directed towards $ q^\mp_n $. (b) Singular cycles are obtained by concatenating fast, intermediate, and slow segments of \eqref{3redlay}, \eqref{em1nlay}, and \eqref{redMn}, respectively.}
		\figlab{ns-sing}
	\end{figure}
	
	In the double singular limit of $ \varepsilon=0=\delta_n $, the \textit{folded singularities} $ Q_n=\lb q_n^-,q_n^+\rb $ \cite{szmolyan2001canards}  of \eqref{3red} are given by 
	\begin{align*}
	q^\mp_n = \mathcal{M}_n\cap\mathcal{L}^\mp;
	\end{align*} 
	see {again \figref{ns-sing}}. In the following, we will write
	\begin{align*}
	q^\mp_n = \lp v_{q_n^\mp}, h_{q_n^\mp}, n_{q_n^\mp}\rp;
	\end{align*}
	{the above coordinates are obtained by solving \eqref{em1}, \eqref{em1V}, and \eqref{acNV}.}
	As in \secref{hslow}, the geometries of $ \mathcal{M}_n $, $ \mathcal{F}_{\mathcal{M}_n} $, and $ Q_n $ depend on the rescaled applied current $\bar I$.
	
	Next, setting $ \delta_n=0$ in the slow formulation, Equation~\eqref{em1nslow}, we obtain the one-dimensional flow on $ \mathcal{M}_{n} $:
	
		\begin{align}
		\begin{aligned}
		0 &=  \partial_h  \left[V(v,m_\infty(v), h,n) \right]\lvert_{h = \eta(v,n)} H(v,\eta(v,n)),  \\
		\dot{n} &=-\partial_v  \left[V(v,m_\infty(v), h,n) \right]\lvert_{h = \eta(v,n)} N(v,n). 
		\end{aligned}
		\eqlab{em1nred}
		\end{align}
	
	Differentiating the algebraic constraint in \eqref{acNV} implicitly, via
	\begin{align*}
	-\partial_v  \left[V(v,m_\infty(v), h_\infty(v),n \right]\dot{v} &= \partial_n   \left[V(v,m_\infty(v), h_\infty(v),n \right]\dot{n},
	\end{align*}
	and using {the $h$-equation of \eqref{em1nred}} and \eqref{Fn}, we find the following expression for the {one-dimensional} flow on $ \mc{N}^{\mp} $:
	\begin{align}
	\dot{v} &=\frac{\partial_v  \left[V(v,m_\infty(v), h,n) \right]\partial_n  \left[V(v,m_\infty(v), h_\infty(v),n) \right]}{\partial_v  \left[V(v,m_\infty(v), h_\infty(v),n \right]} \bigg\lvert_{\{h = h_\infty(v),n=\nu(v,h{_\infty(v)})\}}N(v,\nu(v,h{_\infty(v)})). \eqlab{redMn}
	\end{align}

	Using computer algebra software, we can confirm that, {when a true equilibrium exists in $\mc{S}^r$}, the {one-dimensional}  flow on $ \mc{N}^{\mp} $ is directed towards $ q_n^\mp $, respectively; cf. \figref{ns-sing} and \figref{ns280} {below}. The resulting singular geometry is summarised in \figref{ns-sing}.  {Combining orbit segments from the layer flow of  Equation~\eqref{3redlay}, the intermediate fibres of \eqref{em1nlay} on $ {\mc{M}_2} $, and the {one-dimensional} dynamics of \eqref{redMn} on $ \mc{M}_n $, one can again construct {singular cycles} as in \secref{hslow}.}
	
	We emphasise that, according to the above, Equation~\eqref{3red}, with $ n $ the slow variable, satisfies assumptions that are analogous to P1 through P3 in \secref{hslow}.
	
	\begin{remark}
		We remark that considering the reduced flow {of \eqref{3redslow}} on $\mc{S}^a$ in \eqref{em1red}, instead of \eqref{em1redn} with $\delta_n>0$ small, gives a slow-fast system in the non-standard form of GSPT, {in analogy to the $h$-slow regime described in \remref{hs-ns}}.
		
	\end{remark}
	
	\subsection{Perturbed dynamics and MMOs}

	By standard GSPT \cite{fenichel1979geometric} and according to \cite{kaklamanos2022bifurcations}, there exist invariant manifolds $ \mathcal{S}^{a}_{\varepsilon\delta_n} $ and $\mathcal{S}^{r}_{\varepsilon\delta_n}$ that are diffeomorphic to their unperturbed, normally hyperbolic counterparts $ \mathcal{S}^{a}$ and $\mathcal{S}^{r}$, respectively, and that lie $ \mathcal{O}\lp\varepsilon,\delta_n\rp $-close to them in the Hausdorff distance, for $ \varepsilon$ and $\delta_n$ positive and sufficiently small. The perturbed manifolds $ \mathcal{S}^{a}_{\varepsilon\delta_n} $ and $\mathcal{S}^{r}_{\varepsilon\delta_n}$ are locally invariant under the flow of {the three-dimensional reduction} in \eqref{3redfast}. Moreover, for $\varepsilon $ and $\delta_n$ sufficiently small, there exist invariant manifolds $ \mathcal{N}^{\mp}_{\varepsilon\delta_n} $ and $\mathcal{N}^{r}_{\varepsilon\delta_n}$ that are diffeomorphic, and $ \mathcal{O}\lp\delta_n\rp $-close in the Hausdorff distance, to their unperturbed, normally hyperbolic counterparts $ \mathcal{N}^{\mp}$ and $\mathcal{N}^{r}$, respectively. The perturbed branches $ \mathcal{N}^{\mp}_{\varepsilon\delta_n} $ and $\mathcal{N}^{r}_{\varepsilon\delta_n} $ are locally invariant under the flow of \eqref{3redfast}.  
	
	{In the $n$-slow case,} orbits for the perturbed flow of {the three-dimensional reduction} in \eqref{3redfast} with $ \varepsilon $ and $\delta_n$ positive and sufficiently small can be constructed {in analogy to the $h$-slow case,} by combining perturbations of fast, intermediate, and slow segments of singular trajectories, as described above. In the following, we demonstrate how the various $ \bar{I} $-values that distinguish between qualitatively different firing patterns in Equation~\eqref{3redfast} with $ \varepsilon,\delta_n >0$ small, as shown in \figref{ibifn}, are obtained.

	\subsection{Onset and cessation of oscillatory dynamics}
	
	In the singular limit of $\varepsilon=0=\delta_n$ in {the three-dimensional reduction} \eqref{3redfast}, the {one-dimensional} flow on $\mc{N}$, as defined in Equation~\eqref{redMn}, has a stable equilibrium point on $\mc{S}^{r}$ for $\bar{I}\in (\bar{I}_n^-, \bar{I}_n^+)$, where
	we numerically obtain
	\begin{align} 
	\bar{I}_n^- \simeq \frac{4.8}{k_v g_{Na}}\quad\text{and}\quad 
	\bar{I}_n^+ \simeq \frac{280}{k_v g_{Na}}; \eqlab{Impn}
	\end{align} 
	cf.~\figref{ns280}. Moreover, we again find $ \bar{I}_n^-<\bar{I}_n^p<\bar{I}_n^+$; recall \eqref{Ipn}. The {one-dimensional} flow in \eqref{redMn} has a stable equilibrium point on $\mc{N}^{-}$ for $\bar{I}<\bar{I}_n^-$, and on the unique normally hyperbolic branch $\mc{N}$ for $\bar{I}> \bar{I}_n^+$, see again \figref{ns280}; correspondingly, that equilibrium crosses $ \mc{L}^- $ for $ \bar{I} = \bar{I}_n^- $ and $ \bar{I} = \bar{I}_n^+ $. {(While it may not be obvious from \figref{ns280} whether $\mc{N}$ and $\mc{F}_{\mc{M}_2}$ intersect, such intersections are of no interest to us once the stable equilibrium enters $\mc{S}^a$, as trajectories will ultimately converge to that equilibrium.)}
	
	\begin{figure}[ht!]
		\centering
		\begin{subfigure}[b]{0.3\textwidth}
			\centering
			\includegraphics[scale = 0.22]{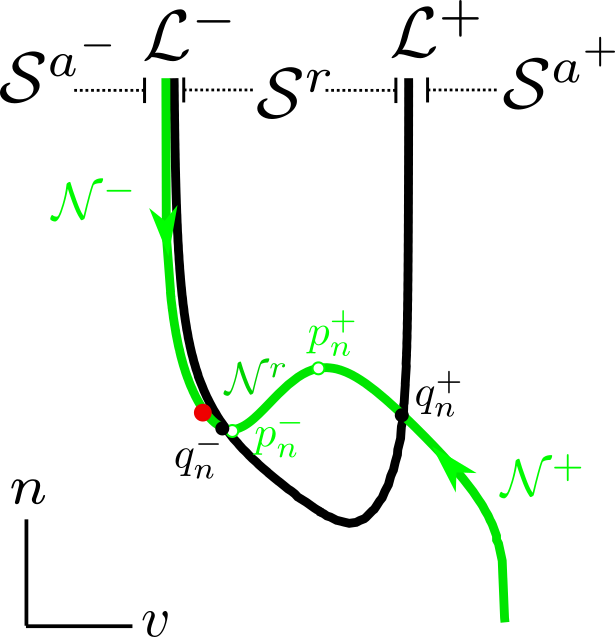}
			\caption{$\bar{I}<\bar{I}_n^-$}
		\end{subfigure}
		~
		\begin{subfigure}[b]{0.3\textwidth}
			\centering
			\includegraphics[scale = 0.22]{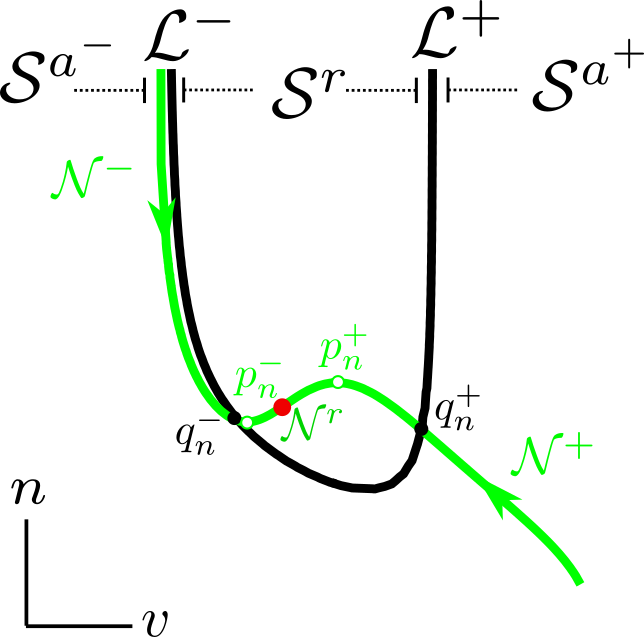}
			\caption{$\bar{I}\in (\bar{I}_n^-,\bar{I}_n^+)$}
		\end{subfigure}
		~
		\begin{subfigure}[b]{0.3\textwidth}
			\centering
			\includegraphics[scale = 0.22]{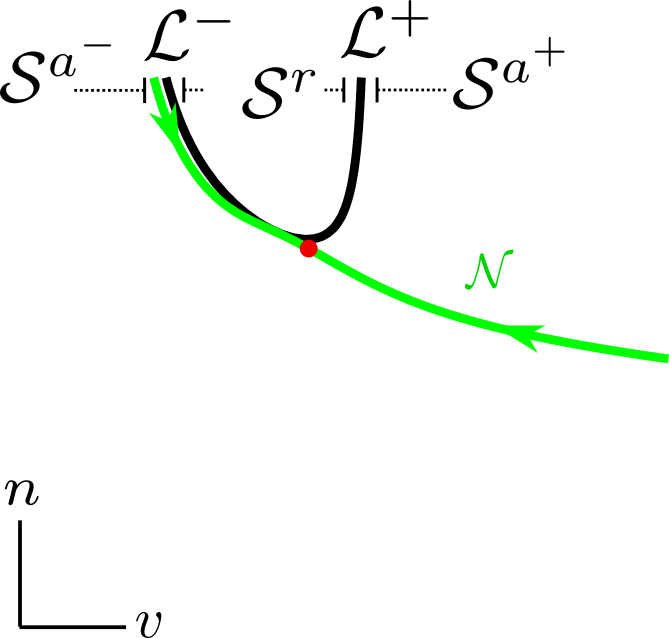}
			\caption{$\bar{I}>\bar{I}_n^+$}
		\end{subfigure}
		\caption{A stable equilibrium of the {one-dimensional} flow in \eqref{redMn} (red dot) lies in the attracting portion $\mc{S}^a$ of ${\mc{M}_2}$ for $\bar{I}<\bar{I}_n^-$ and $\bar{I}>\bar{I}_n^+$, as seen in panels (a) and (c), respectively, and in the repelling portion $\mc{S}^r$ for $\bar{I}\in (\bar{I}_n^-,\bar{I}_n^+)$; cf.~panel (b). As $ \bar{I} $ increases, the distance between $ \mc{L}^- $ and $ \mc{L}^+ $ decreases, with the leftmost branch of $ \mc{N} $ ``moving downward''. For $\bar{I} >\bar{I}_n^p$, the manifold $\mc{M}_n$ has no fold points; it then consists of a unique branch $\mc{N}$ which is attracting on $\mc{S}^a$ and repelling on $\mc{S}^r$.}
		\figlab{ns280}
	\end{figure}
	
	For $\bar{I}\in (\bar{I}_{n}^-, \bar{I}_{n}^+)$ and away from the tangential connection of $\mc{L}^\mp$, the geometry of {the three-dimensional reduction} in \eqref{3redfast} is again similar to that of the prototypical system introduced in \cite{kaklamanos2022bifurcations} {in the $n$-slow case}. {For $\varepsilon,\delta_n>0$ sufficiently small, \eqref{3redfast} undergoes Hopf bifurcations for $\bar{I}$-values that are $\mc{O}(\varepsilon,\delta_n)$-close to $\bar{I}_n^\mp$, which implies that the ``full" four-dimensional system, Equation~\eqref{hhmod-fast}, undergoes Hopf bifurcations for $\bar{I}$-values  that are $\mc{O}(\gamma,\varepsilon,\delta_n)$-close to $\bar{I}_n^\mp$. The Hopf bifurcation of Equation~\eqref{hhmod-fast} near $\bar{I}_n^-$ is subcritical, while the one close to $\bar{I}_n^+$ is supercritical \cite{doi2001complex,rubin2007giant}. It then follows that, for fixed $\bar{I}\in (\bar{I}_n^-, \bar{I}_n^+)$, there exist $\gamma,\varepsilon,\delta_n>0$ sufficiently small such that Equation~\eqref{3redfast} -- and, by extension, Equation~\eqref{hhmod-fast} -- features global mixed-mode dynamics.} By numerical sweeping, we find that for $\gamma=0.083$, $\varepsilon = 0.1$, $\delta_n = 0.01$, and $\delta_h=1$, with the values of the remaining parameters as in \eqref{parbars}, the onset of oscillatory dynamics in \eqref{hhmod-fast} occurs at approximately $\bar{I} \simeq 6.6/(k_v g_{Na})$, while the cessation of oscillatory dynamics is observed at approximately $\bar{I} \simeq 268/(k_v g_{Na})$; these values are close to $\bar{I}_{n}^-$ and $\bar{I}_{n}^+$, respectively, as expected. 

 \begin{remark}
  {We note that $\mathcal N^-$ appears almost aligned with $\mathcal L^{-}$ in \figref{ns280}. While we assume in the present work that the angle between the two manifolds is $\mathcal O(1)$, it may be of interest to investigate this apparent alignment further.}
  \end{remark}
	\subsection{Double and single epochs {of slow dynamics}}
	
	We characterise the pair of folded singularities $ q_n^\mp $ based on whether they are ``connected'' by singular cycles or not, in accordance with \defnref{orbconHH}, which is equivalent to the following classification in terms of the $n$-coordinates $n_{q_n}^\mp$ of the folded singularities $q_n^{\mp}$.
	\begin{definition}\textcolor{white}{white text}
		\begin{enumerate}
			\item The folded singularities $ q_n^\mp $ are \textit{orbitally connected} if $ n_{q_n}^->n_{q_n}^+ $; cf.~\figref{ns-rel}(a).
			\item The folded singularities $ q_n^\mp $ are \textit{orbitally aligned} if $ n_{q_n}^-=n_{q_n}^+ $; cf.~\figref{ns-rel}(b).
			\item The folded singularities $ q_n^\mp $ are \textit{orbitally remote} if $ n_{q_n}^-<n_{q_n}^+ $; cf.~\figref{ns-rel}(c). 
		\end{enumerate}
		\defnlab{orbconNN}
	\end{definition}
	\begin{figure}[ht!]
		\centering
		\begin{subfigure}[b]{0.3\textwidth}
			\centering
			\includegraphics[scale = 0.15]{pics/hh/nslow7.png}
			\caption{ $ {I} < \bar{I}_n^a $}
		\end{subfigure}
		\begin{subfigure}[b]{0.3\textwidth}
			\centering
			\includegraphics[scale = 0.15]{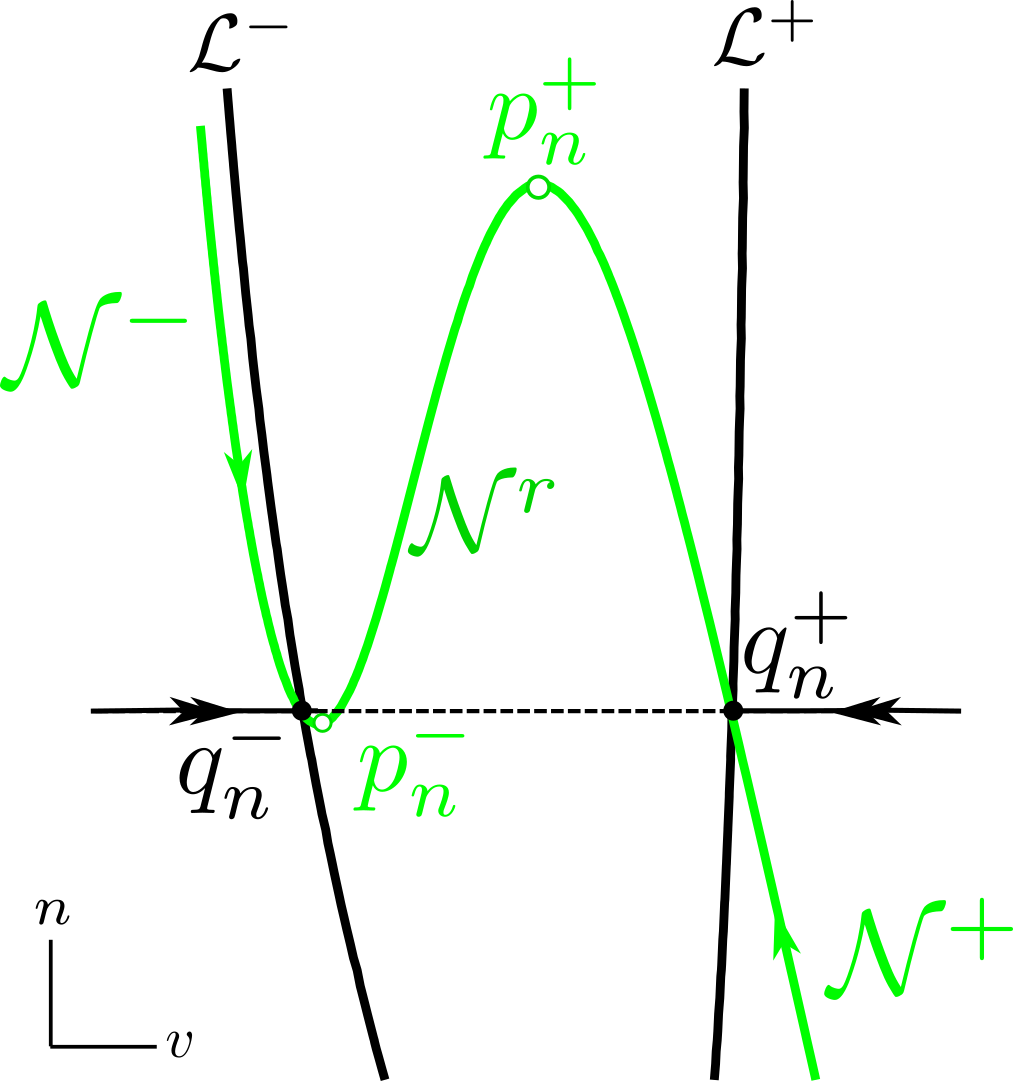}
			\caption{ $ \bar{I} = \bar{I}_n^a  $}
		\end{subfigure}
		\begin{subfigure}[b]{0.3\textwidth}
			\centering
			\includegraphics[scale = 0.15]{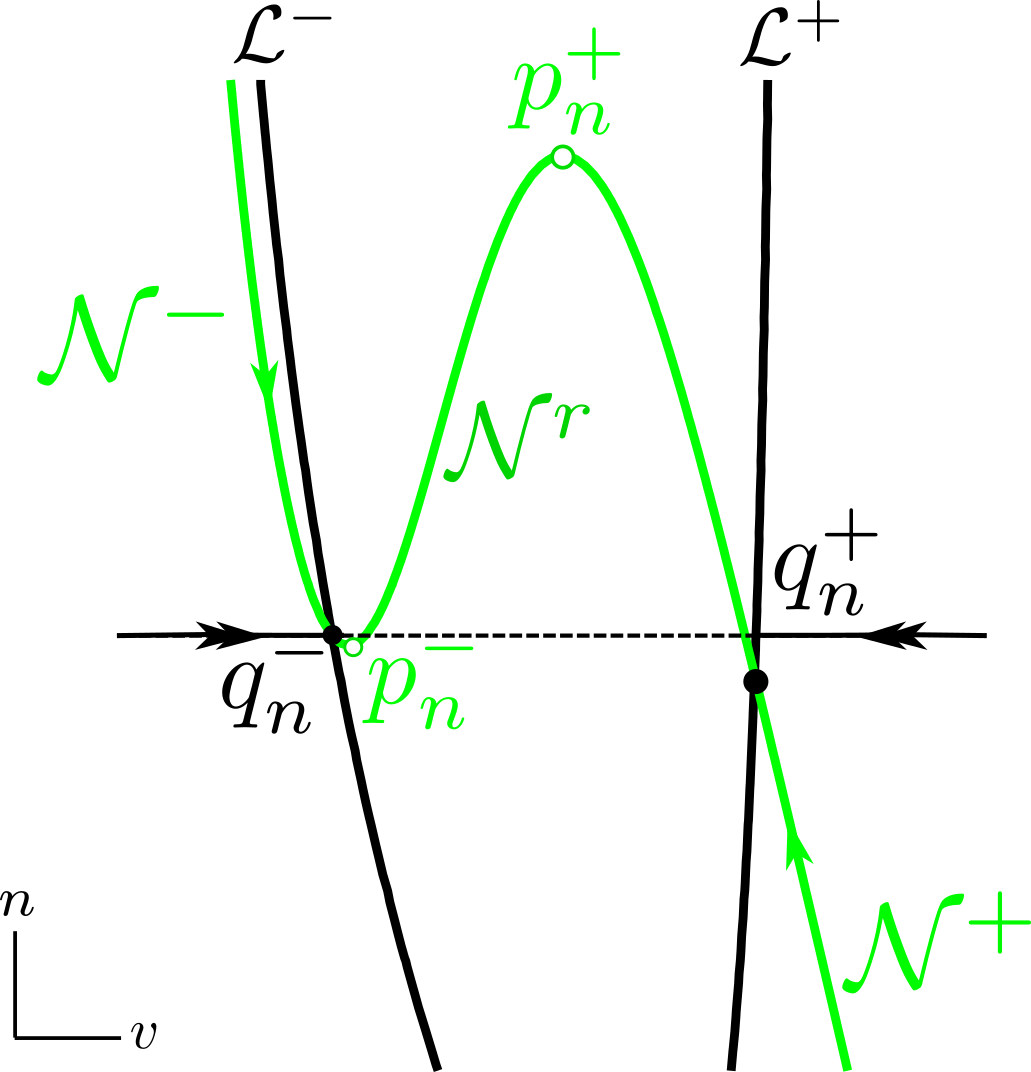}
			\caption{$ \bar{I} > \bar{I}_n^a $}
		\end{subfigure}
		\caption{Orbital connection, or lack thereof, in the limit of $ \varepsilon=0=\delta_n $, in accordance with \defnref{orbconNN}. (a) For $ \bar{I} < \bar{I}_n^a  $, the folded singularities $ q^\mp_n  $ are orbitally connected, in the sense that there exists a singular cycle that passes through both of them and that also contains slow segments on both $ \mc{N}^{\mp} $. (b) For $ \bar{I} = \bar{I}_n^a  $, $ q^\mp_n $ are orbitally aligned, in that there exists a singular cycle that passes through both of them but that contains no slow segments on either $ \mc{N}^{\mp}$ for $ \varepsilon=0=\delta_n $. (c) For $ \bar{I} > \bar{I}_n^a  $, $ q^\mp_n $ are orbitally remote, in that a singular cycle that passes through one of the folded singularities does not pass through the other.} 
		\figlab{ns-rel}
	\end{figure}
	
	Due to the $\bar I$-dependence of the location of $q_n^\mp$, the classification in \defnref{orbconNN} hence encodes the position of the folded singularities $q_n^\mp$ relative to one another {in terms of} the rescaled applied current $ \bar{I} $.  It follows that, in the parameter regime given by \eqref{parbars}, there exists a unique value $\bar{I}_n^a$  such that the singularities $q_n^-$ and $q_n^+$ are orbitally aligned. Numerically, that value is found to equal 
	\begin{align*}
	\bar{I}_n^a \simeq \frac{10.1}{{k_v g_{Na}}}; 
	\end{align*} 
	recall \eqref{parbars}. Correspondingly, $ q^\mp_n $ are orbitally connected for $\bar I<\bar{I}_n^a$ and orbitally remote for $\bar I>\bar{I}_n^a$.  {We have the following result:}
	{\begin{proposition}
			Fix $\bar{I}\in(\bar{I}_n^-,\bar{I}_n^a)$ and let $\delta_h = 1$. Then, there exist $\gamma_0, \varepsilon_0$, and $\delta_{0}$ positive and sufficiently small such that the HH equations in \eqref{hhmod-fast} feature MMOs with double epochs of {slow dynamics} for all $(\gamma, \varepsilon, \delta_n)\in(0,\gamma_0)\times(0, \varepsilon_0)\times(0, \delta_{0})$.
			\proplab{nsdouble}
	\end{proposition}}
	
	\begin{proof}
		The proof is analogous to the proof of \propref{hsdouble} in \secref{hslow}.
	\end{proof}
	
	\propref{nsdouble} states that for every fixed $\bar{I}\in(\bar{I}_n^-,\bar{I}_n^a)$, there exist $\gamma, \varepsilon, \delta_n>0$ sufficiently small such that {the three-dimensional reduction} in \eqref{3redfast} -- and, by extension, Equation~\eqref{hhmod-fast} -- exhibits MMOs with double epochs of {slow dynamics}. Conversely, for fixed $ \varepsilon,\delta_n>0 $ sufficiently small, \eqref{3redfast} admits MMOs with double epochs of {slow dynamics} for $\bar{I}\in (\bar{I}_n^-+\mc{O}(\varepsilon, \delta_n), \bar{I}_n^a+\mc{O}(\varepsilon, \delta_m))$; by extension, for $ \gamma,\varepsilon $, and $ \delta_{n}>0 $ sufficiently small, the mixed-mode dynamics of the full four-dimensional system in \eqref{hhmod-fast} is characterised by double epochs of {slow dynamics} for $\bar{I}\in (\bar{I}_n^-+\mc{O}(\gamma,\varepsilon, \delta_n), \bar{I}_n^a+\mc{O}(\gamma,\varepsilon, \delta_n))$; cf. panels (a) and (b) of \figref{NStimeseries} and \figref{NSphase}. By numerical sweeping, we find that for $\gamma=0.083$, $\varepsilon = 0.1$, and $\delta_n = 0.01$, with the values of the remaining parameters as in \eqref{parbars}, \eqref{hhmod-fast} features double-epoch MMOs until approximately $\bar{I} \simeq 9.4/(k_v g_{Na})$, which is consistent with the above error estimates.
	

	
	\begin{figure}[ht]
		\centering
		\begin{subfigure}[b]{0.45\textwidth}
			\centering
			\includegraphics[scale = 0.5]{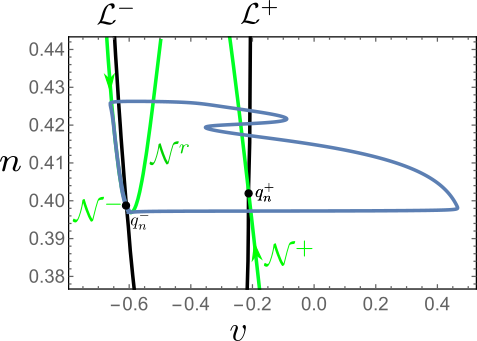}
			\caption{$ {I} = 9 $}
		\end{subfigure}
		~
		\begin{subfigure}[b]{0.45\textwidth}
			\centering
			\includegraphics[scale = 0.5]{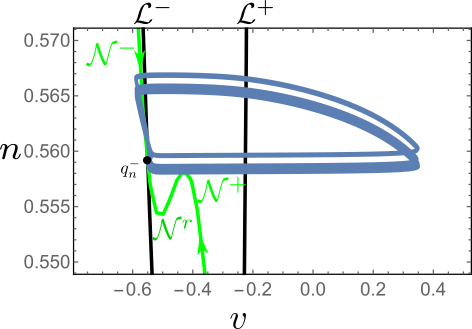}
			\caption{$ {I} = 64.5 $}
		\end{subfigure}
		\\
		\centering
		\begin{subfigure}[b]{0.45\textwidth}
			\centering
			\includegraphics[scale = 0.5]{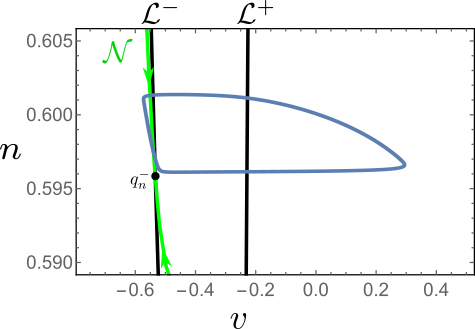}
			\caption{$ {I} = 90 $}
		\end{subfigure}
		~
		\begin{subfigure}[b]{0.45\textwidth}
			\centering
			\includegraphics[scale = 0.5]{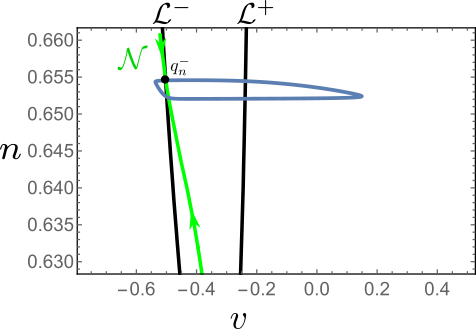}
			\caption{$ {I} = 150 $}
		\end{subfigure}
		\caption{Singular geometry of Equation~\eqref{em1ninter} which underlies the qualitative properties of the time series, and the associated MMO trajectories, illustrated  in \figref{NStimeseries}: as $\bar{I}$ is varied, the orbitally connected folded singularities $q_n^\mp$ become orbitally remote. Therefore, Equation~\eqref{hh1f} with $ \varepsilon =0.0083 $, $\tau_n=100$, and varying values of $ {I} =k_vg_{Na}\bar{I}$, transitions from {oscillatory dynamics with double epochs of {slow dynamics}} for $ \bar{I}<\bar{I}_n^a $ (panel (a)) to single epochs thereof for $ \bar{I}>\bar{I}_n^a $  (panels (b), (c), and (d)).}
		\figlab{NSphase}
	\end{figure}
	
	In contrast to the $ h $-slow regime, we do not find an $ \bar{I} $-value for which MMOs with single epochs {of slow dynamics} turn into two-timescale relaxation oscillation. To that end, we consider Equation~\eqref{em1ninter} with $ \delta_n>0 $ small. Eliminating time, away from $ \mc{N} $ we can then approximate the flow on $ \mc{S}^{a^\mp} $ by
	\begin{align}
	\frac{\tn{d}n}{\tn{d}v} = -\delta_n\frac{\partial_v  \left[V(v,m_\infty(v), h,n) \right]\lvert_{h = \eta(v,n)} N(v,n)}{\partial_h  \left[V(v,m_\infty(v), h,n) \right]\lvert_{h = \eta(v,n)} H(v,\eta(v,n))},
	\end{align}
	where we have omitted $ \mc{O}\lp \delta_n\rp $-terms in the denominator; recall that $\partial_h \left[V(v,m_\infty(v),h,n)\right]<0$, by \eqref{parders}. 
	
	We recall that we denote by $P(\mc{L}^\mp)\in \mc{S}^{a^\pm}$ the projection of $\mc{L}^-$ onto $\mc{S}^{a^\pm}$ along the fast fibres of \eqref{3redlay}, and we now define the points
	\begin{align*}
	(v_{\tn{min}}, h_{\tn{min}},  n_{q_n}^- ) := P(\mc{L}^+) \cap \lb h = n_{q_n}^-\rb\quad\text{and}\quad 
	(v_{\tn{max}}, h_{\tn{max}},  n_{q_n}^- ) := P(\mc{L}^-) \cap \lb n = n_{q_h}^-\rb,
	\end{align*}
	as illustrated in \figref{umax} {for the $h$-slow regime}. A trajectory that leaves the vicinity of $ \mc{L}^- $ at a point with $ n $-coordinate close to $ n_{q_n}^- $ returns to a point with $n$-coordinate $ \hat{n} $ in that vicinity after one large excursion, where
	\begin{align}
	\begin{aligned}
	\hat{n} = n +\delta_n \lp \mc{G}^-_n(n_{q^-_n}, \bar{I})+\mc{G}^+_n(n_{q^-_n}, \bar{I})\rp+\mc{O}(\delta_n^2),
	\end{aligned}\eqlab{nret}
	\end{align}
	and where we now denote
	\begin{align*}
	\mc{G}^-_n(n, \bar{I}) &:= -\int_{{v}_{\tn{min}}}^{{v}_{q^-_n} } \frac{\partial_v  \left[V\lp v,m_\infty(v), h,n\rp \right]\lvert_{h = \eta\lp v,n\rp} N\lp v,n\rp }{\partial_h  \left[V\lp v,m_\infty(v), h,n\rp \right]\lvert_{h = \eta\lp v,n\rp} H(v,\eta\lp v,n\rp)}\tn{d}v\quad\text{and} \\
	\mc{G}^+_n(n, \bar{I}) &:= -\int_{{v}_{\tn{max}}}^{{v}_{q^+_n} } \frac{\partial_v  \left[V\lp v,m_\infty(v), h,n\rp \right]\lvert_{h = \eta\lp v,n\rp} N\lp v,n\rp }{\partial_h  \left[V\lp v,m_\infty(v), h,n\rp \right]\lvert_{h = \eta\lp v,n\rp} H(v,\eta\lp v,n\rp)}\tn{d}v;
	\end{align*}
	see \cite{krupa2008mixed,kaklamanos2022bifurcations} for details. 
	
	Equation~\eqref{nret} encodes the displacement in the $n$-direction after a large excursion, i.e., after the trajectory has been attracted to $\mathcal{S}^{a^+}$ and then again repelled away therefrom in the vicinity of $\mathcal{L}^+$ before its return to $\mathcal{S}^{a^-}$ near $\mathcal{L}^-$; see \cite{krupa2008mixed,kaklamanos2022bifurcations} for details. We define the displacement function as
	\begin{align*}
	\Delta_n(n,\bar{I}) = \mc{G}^-_n(n, \bar{I})+\mc{G}^+_n(n, \bar{I});
	\end{align*}
	then, direct calculation shows that 
	\begin{align*}
	\Delta_n({n_q^-},\bar{I})
	>0 \quad\tn{for }\bar{I}\in (\bar{I}_n^a, \bar{I}_n^+),
	\end{align*}
	which implies that trajectories are attracted to $\mc{N}$ for $\bar{I}\in(\bar{I}_n^a, \bar{I}_n^+)$. We conclude that no two-timescale relaxation oscillation occurs in \eqref{3redfast} for $n$ slow; rather, for $\gamma, \varepsilon, \delta_n>0$ sufficiently small, \eqref{hhmod-fast} features trajectories with single epochs of {slow dynamics} when $\bar{I}\in (\bar{I}_n^a+\mc{O}(\gamma, \varepsilon, \delta_n), \bar{I}_n^++\mc{O}(\gamma, \varepsilon, \delta_n))$, as is also evidenced numerically in panels (c) and (d) of \figref{NSphase}, {where trajectories seem to be attracted to $\mc{N}^-$ -- or to the unique branch $\mc{N}$ if $\bar{I}>\bar{I}_n^p$}.

	We remark that, for the $ I $-values considered in \cite{doi2001complex}, the equilibrium of \eqref{em1ninter} as depicted in {\figref{ns280}} lies on $ \mc{N}^r$, which is why in \cite{doi2001complex} the $ n $-slow regime is found to be {fundamentally} different to the $ h $-slow one; moreover, it is postulated there that chaotic dynamics can be realised only in the $h$-slow regime. Here, we show that, for sufficiently large $ \bar{I} $-values, that equilibrium point lies on the {unique branch $ \mc{N} $ in $\mc{S}^r$, recall \eqref{Ipn} and \eqref{Impn}}. Motivated also by the time series in \figref{NStimeseries}(b), we remark that it would be worth investigating whether chaotic dynamics can be realised in the corresponding $ \bar{I} $-regime. The search for such dynamics and the potential study of the corresponding generating mechanisms is included in plans for future work. 
	
	\section{Conclusions}
	\seclab{HHconclusion}
	
	In this work, we have proposed a novel and global three-dimensional reduction of the four-timescale Hodgkin-Huxley (HH) equations in  \eqref{original} \cite{doi2001complex,HHmain} that is based on a scaled system, Equation~\eqref{hh}, first proposed in  \cite{rubin2007giant}; cf. \thmref{3reduction}. Specifically, we have reduced \eqref{hh} to a globally normally hyperbolic -- and, in fact, normally attracting -- slow manifold, instead of to a local centre manifold, as was done in \cite{rubin2007giant, rubin2008selection}, {and we have shown that the  two reductions feature the same critical manifold ${\mc{M}_2}$.} We have then concluded that, depending on the values of the parameters in the system relative to each other, our reduction is itself a three-timescale singular perturbation problem. We emphasise that this three-timescale structure is also apparent in the {three-dimensional reduction of} Equation~\eqref{RW} obtained in \cite{rubin2007giant, rubin2008selection} if either $ \delta_h>0 $ is assumed to be sufficiently small and $ \delta_n = \mc{O}(1) $, or if $ \delta_n>0 $ is sufficiently small and $ \delta_h = \mc{O}(1) $; however, those regimes were not considered in \cite{rubin2007giant,rubin2008selection}.
	
	The timescale separation achieved in our {three-dimensional reduction}, Equation~\eqref{3red}, allows for a further iterative dimension reduction to a hierarchy of invariant manifolds via the systematic study of the associated layer problems and reduced flows. By decomposing the global dynamics of \eqref{3red} into segments that evolve on different timescales, we are able to classify the resulting mixed-mode dynamics -- as well as the dynamics of the ``full" (transformed) HH model, Equation~\eqref{hhmod-fast}. In particular, we can explain transitions between oscillatory trajectories with different qualitative properties, in accordance with the geometric mechanisms proposed in \cite{kaklamanos2022bifurcations}.
	
	Specifically, for the regimes where either $ h $ or $ n $ is taken to be the slowest variable in \eqref{hhmod-fast}, we have classified the different firing patterns which are correspondingly illustrated in \figref{HStimeseries} and \figref{NStimeseries}, and we have described the transitions between them, in the framework of multiple-timescale GSPT; cf. \figref{HSphase} and \figref{NSphase}, respectively. (We reiterate that these patterns had previously been documented in \cite{doi2001complex}.) {In both regimes, we have characterised mixed-mode dynamics with either double or single epochs of {slow dynamics} upon variation of the rescaled applied current $\bar{I}$. For the parameter values considered here, we have observed (two-timescale) relaxation oscillation and ``exotic" MMO trajectories for $h$ slow only.} {We reiterate that the methodology presented here extends to a wide variety of multiple-scale models from mathematical physiology that can be expressed in an HH-type formalism \cite{cloete2020dual,diekman2021circadian,forger2017biological, iosub2015calcium, kugler2018early, pages2019model, takano2021highly}. {In such models, the physiological properties of the particular system under consideration determine whether small parameters -- and, hence, a separation of timescales -- are present in the resulting, singularly perturbed system of differential equations.}}
	
	The local mechanisms which generate small-amplitude (``sub-threshold") oscillations in \eqref{3red} can also be analysed via the approach that is documented in \cite{kaklamanos2022bifurcations} in our three-timescale context; cf. Appendix~\ref{partpert} for details. We remark that it would seem natural to attempt a normal form transformation in order to make the singular geometry of \eqref{3redfast} independent of the external applied current $\bar{I}$. However, we have found that such a transformation would encode transitions between MMO trajectories with different qualitative properties in the reduced flow on $\mc{S}$, resulting in a slow-fast system in the non-standard form of GSPT \cite{wechselberger2020geometric}. Moreover, we have adopted the widely accepted scaling in \eqref{scaling} \cite{rubin2007giant} throughout; an investigation of alternative scalings is left for future work; recall, in particular, the proof of \propref{hsrelax} and \remref{scaling}.

	Finally, in \cite{doi2001complex} it was demonstrated numerically that Equation~\eqref{original} features chaotic dynamics in the $ h$-slow regime; however, the underlying chaos-generating mechanisms were not explained. {Moreover, no chaotic trajectories were found in the $ n $-slow regime, and it was postulated that this difference is due to the fact that the unique equilibrium point of the system lies on $ \mc{H}^{-}\cap\mc{S}^r $ (in our notation) for $h$ slow, while in the $ n $-slow regime it is located on $ \mc{N}^r $ for the values of $I$ considered by them; cf. \figref{hs280} and \figref{ns280}, respectively. Here, we have shown that by considering larger values of $ I $, that equilibrium can equally be made to lie on {the unique branch $ \mc{N} $ in $\mc{S}^r$ }} in the $ n $-slow regime; cf.~\figref{ns280}(c) {and recall \eqref{Ipn} and \eqref{Impn}}. Hence, motivated by the time series in 	panel (c) of \figref{NStimeseries}, where slow segments with and without SAOs alternate ``below'', we suggest that it would be worth investigating whether chaotic dynamics is equally possible in that regime, since our analysis shows that, ultimately, the two regimes are not fundamentally different, as claimed in \cite{doi2001complex} {-- an observation which has also been made by Rubin and Wechselberger \cite{rubin2007giant,rubin2008selection}}. A more systematic analysis of the mechanisms that are responsible for chaotic behaviour in \eqref{original}, and the potential relation thereof to Shilnikov-type homoclinic phenomena or to period-doubling bifurcations of small-amplitude periodic orbits within the framework of three-timescale GSPT, for both $ h $ slow and $ n $ slow, is left for future work. 
	
	\section*{Acknowledgements} {The content of this work was part of the first author's PhD thesis, completed between 2018 and 2021
at the University of Edinburgh \cite{kaklamanos2022mixed}.} The authors would like to thank Martin Krupa and Martin Wechselberger for insightful discussions and constructive feedback. {The authors are also grateful to two anonymous reviewers for comments and suggestions which greatly improved the original manuscript.}
	
	\appendix
	
	\section{Local dynamics and SAO-generating mechanisms}
	\label{partpert}
	The local mechanisms underlying the {slow} dynamics of {the three-dimensional reduction} of the Hodgkin-Huxley (HH) model in \eqref{3redfast} are similar to the ones described in \cite{kaklamanos2022bifurcations} for the extended prototypical model introduced there, and are only discussed in brief here.
	
	We reiterate the fast formulation in \eqref{3redfast} for convenience,
	\begin{align}
	\begin{aligned}
	{v}' &= \frac{m_\infty(v)-\mu(v,h,n)}{t_m(v)\partial_v\mu(v,h,n)} -\varepsilon\delta_hH(v,h)\frac{\partial_h\mu(v,h,n)}{\partial_v\mu(v,h,n)}-\varepsilon\delta_nN(v,n)\frac{\partial_n\mu(v,h,n)}{\partial_v\mu(v,h,n)}, \\
	&=: U(v,h,n;\varepsilon, \delta_h, \delta_n) \\ 
	{h}' &= \varepsilon \delta_hH(v,h), \\
	{n}' &= \varepsilon \delta_nN(v,n), 
	\end{aligned}
	\eqlab{3app}
	\end{align}
	and we consider the regime where $h$ is slow, with $ \varepsilon,\delta_h>0 $ sufficiently small and $ \delta_n =\mc{O}(1)$. In the singular limit of $ \delta_h=0 $ with $ \varepsilon>0 $ and $ \delta_n=1 $, we obtain the system
		\begin{align}
		\begin{aligned}
		{v}' &=\frac{m_\infty(v)-\mu(v,h,n)}{t_m(v)\partial_v\mu(v,h,n)} -\varepsilon N(v,n)\frac{\partial_n\mu(v,h,n)}{\partial_v\mu(v,h,n)}= U(v,h,n;\varepsilon, 0, 1), \\ 
		{h}' &= 0, \\
		{n}' &= \varepsilon N(v,n), 
		\end{aligned}
		\eqlab{3appsub}
		\end{align}
	the equilibria of which are located on $ \mc{M}_h $. Linearisation of \eqref{3appsub} about $ \mc{M}_h $ gives the Jacobian matrix
	\begin{align}
	\textbf{A}_h = \begin{pmatrix}
	\partial_vU(v,h,n;\varepsilon, 0, 1) & \partial_nU(v,h,n;\varepsilon, 0, 1) \\
	\varepsilon\partial_vN(v,n) & \varepsilon\partial_nN(v,n)
	\end{pmatrix}
	\eqlab{jacAh}\bigg\lvert_{\mc{M}_h};
	\end{align}
	the eigenvalues of $ \textbf{A}_h $ determine the regimes where $ \mc{M}_h $ is either focally or nodally attracting, respectively repelling, under the two-dimensional flow of \eqref{3appsub}. In particular, there exists a point on $ \mc{M}_h $ at which \eqref{3appsub} undergoes a Hopf bifurcation; whether this Hopf bifurcation is subcritical or supercritical depends on the value of $ \bar I $. 
	
	In the fully perturbed system, Equation~\eqref{3app}, with $ \varepsilon,\delta_h>0 $ sufficiently small, trajectories entering the focally attracting regime on $ \mc{M}_h $ give rise to bifurcation delay-type SAOs, while attraction to the nodally attracting regime implies the absence of SAOs -- or the presence of only a few thereof that occur close to the jump point; cf. panels (a) and (b) of \figref{HStimeseries}, respectively. In both cases, trajectories jump away from $ \mc{M}_h $ not immediately after entering the repelling region thereon but, rather, after {an $\mathcal O(1)$ delay}. 
 The corresponding jump-off point can {-- under the assumption of analyticity --} be calculated via a way-in/way-out function, as is discussed in more detail in \cite{hayes2016a,kaklamanos2022bifurcations,krupa2010local,letson2017analysis}; {cf. \figref{HStimeseries} and \figref{HSphase}}. We note that this delay-type mechanism is a residual of the 2-fast/1-slow nature of Equation~\eqref{3app}, as well as that a degenerate node where the stability of $ \mc{M}_h $ turns from focal to nodal attraction lies $ \mc{O}(\sqrt{\varepsilon}) $-close to the Hopf bifurcation point of \eqref{3appsub} on $ \mc{M}_h $. 
	
	Finally, the small periodic orbits that emerge at the Hopf bifurcation in \eqref{3appsub} cease to exist at a transverse intersection between $ \mc{S}_{\varepsilon\delta_h}^a $ and $ \mc{S}_{\varepsilon\delta_h}^r $ that occurs $ \mc{O}(\varepsilon) $-close to the Hopf bifurcation point. These periodic trajectories of the partially perturbed system in \eqref{3appsub} with $ \varepsilon,\delta_h>0 $ sufficiently small are associated with \textit{secondary canards} \cite{wechselberger2005existence} and \textit{sector of rotations} \cite{krupa2008mixed}; therefore, if trajectories of the fully perturbed system are attracted to the corresponding region on $ \mc{M}_h $, \textit{canard-induced SAOs} arise; see \cite{krupa2008mixed,de2016sector,desroches2012mixed} for details. These are a residual of the 1-fast/2-slow nature of Equation~\eqref{3app}. 
	
	A detailed quantitative analysis of the above classification is beyond the scope of this work. However, we emphasise that, although the reduction in \eqref{RW} proposed by Rubin and Wechselberger \cite{rubin2007giant} and \eqref{3app} admit the same critical manifolds $ {\mc{M}_2} $ and $ \mc{M}_h $, there are quantitative differences in terms of the regions where the stability of $ \mc{M}_h $ changes; these are reflected in the $\varepsilon$-dependence of the first row in $ {\textbf{A}}_h $ in \eqref{jacAh}, whereas that row would be independent of $ \varepsilon $ if one considered \eqref{RW} instead. 
	
	Finally, we remark that the analysis of the $ n $-slow regime is analogous, with a corresponding Jacobian matrix
	\begin{align*}
	\textbf{A}_n = \begin{pmatrix}
	\partial_vU(v,h,n;\varepsilon, 1,0) & \partial_hU(v,h,n;\varepsilon, 1, 0) \\
	\varepsilon\partial_vH(v,h) & \varepsilon\partial_hH(v,h)
	\end{pmatrix}
	\bigg\lvert_{\mc{M}_n}
	\end{align*}
	about $ \mc{M} _n$, and where the function $ U(v,h,n;\varepsilon, 1,0) $ contains $ \varepsilon H(v,h) $-dependent terms instead of $ \varepsilon N(v,n) $-dependent ones.

	\bibliographystyle{siam}
	\bibliography{refs}

\begin{thebibliography}{10}

\bibitem{brons2006mixed}
{\sc M.~Br{\o}ns, M.~Krupa, and M.~Wechselberger}, {\em Mixed mode oscillations
  due to the generalized canard phenomenon}, Fields Institute Communications,
  49 (2006), pp.~39--63.

\bibitem{cloete2020dual}
{\sc I.~Cloete, P.~J. Bartlett, V.~Kirk, A.~P. Thomas, and J.~Sneyd}, {\em Dual
  mechanisms of ca2+ oscillations in hepatocytes}, Journal of Theoretical
  Biology, 503 (2020), p.~110390.

\bibitem{de2016sector}
{\sc P.~De~Maesschalck, E.~Kutafina, and N.~Popovi{\'c}}, {\em
  Sector-delayed-{H}opf-type mixed-mode oscillations in a prototypical
  three-time-scale model}, Applied Mathematics and Computation, 273 (2016),
  pp.~337--352.

\bibitem{desroches2012mixed}
{\sc M.~Desroches, J.~Guckenheimer, B.~Krauskopf, C.~Kuehn, H.~M. Osinga, and
  M.~Wechselberger}, {\em Mixed-mode oscillations with multiple time scales},
  SIAM Review, 54 (2012), pp.~211--288.

\bibitem{diekman2021circadian}
{\sc C.~O. Diekman and N.~Wei}, {\em Circadian rhythms of early
  afterdepolarizations and ventricular arrhythmias in a cardiomyocyte model},
  Biophysical Journal, 120 (2021), pp.~319--333.

\bibitem{doi2001complex}
{\sc S.~Doi, S.~Nabetani, and S.~Kumagai}, {\em Complex nonlinear dynamics of
  the hodgkin--huxley equations induced by time scale changes}, Biological
  Cybernetics, 85 (2001), pp.~51--64.

\bibitem{fenichel1979geometric}
{\sc N.~Fenichel}, {\em Geometric singular perturbation theory for ordinary
  differential equations}, Journal of Differential Equations, 31 (1979),
  pp.~53--98.

\bibitem{forger2017biological}
{\sc D.~B. Forger}, {\em Biological clocks, rhythms, and oscillations: the
  theory of biological timekeeping}, MIT Press, 2017.

\bibitem{hayes2016a}
{\sc M.~G. Hayes, T.~J. Kaper, P.~Szmolyan, and M.~Wechselberger}, {\em
  Geometric desingularization of degenerate singularities in the presence of
  fast rotation: A new proof of known results for slow passage through hopf
  bifurcations}, Indagationes Mathematicae, 27 (2016), pp.~1184--1203.

\bibitem{hek2010geometric}
{\sc G.~Hek}, {\em Geometric singular perturbation theory in biological
  practice}, Journal of Mathematical Biology, 60 (2010), pp.~347--386.

\bibitem{HHmain}
{\sc A.~L. Hodgkin and A.~F. Huxley}, {\em A quantitative description of
  membrane current and its application to conduction and excitation in nerve},
  The Journal of Physiology, 117 (1952), pp.~500--544.

\bibitem{iosub2015calcium}
{\sc R.~Iosub, D.~Avitabile, L.~Grant, K.~Tsaneva-Atanasova, and H.~J.
  Kennedy}, {\em Calcium-induced calcium release during action potential firing
  in developing inner hair cells}, Biophysical journal, 108 (2015),
  pp.~1003--1012.

\bibitem{izhikevich2007dynamical}
{\sc E.~M. Izhikevich}, {\em Dynamical systems in neuroscience}, MIT press,
  2007.

\bibitem{jelbart2020two}
{\sc S.~Jelbart and M.~Wechselberger}, {\em Two-stroke relaxation oscillators},
  Nonlinearity, 33 (2020), p.~2364.

\bibitem{kaklamanos2022mixed}
{\sc P.~Kaklamanos}, {\em Mixed-mode oscillations in singularly perturbed
  three-timescale systems}, PhD thesis, University of Edinburgh, 2022.

\bibitem{kaklamanos2022bifurcations}
{\sc P.~Kaklamanos, N.~Popovi{\'c}, and K.~U. Kristiansen}, {\em Bifurcations
  of mixed-mode oscillations in three-timescale systems: {A}n extended
  prototypical example}, Chaos: An Interdisciplinary Journal of Nonlinear
  Science, 32 (2022), p.~013108.

\bibitem{krupa2008mixed}
{\sc M.~Krupa, N.~Popovi{\'c}, and N.~Kopell}, {\em Mixed-mode oscillations in
  three time-scale systems: a prototypical example}, SIAM Journal on Applied
  Dynamical Systems, 7 (2008), pp.~361--420.

\bibitem{krupa2001extending}
{\sc M.~Krupa and P.~Szmolyan}, {\em Extending geometric singular perturbation
  theory to nonhyperbolic points---fold and canard points in two dimensions},
  SIAM Journal on Mathematical Analysis, 33 (2001), pp.~286--314.

\bibitem{krupa2010local}
{\sc M.~Krupa and M.~Wechselberger}, {\em Local analysis near a folded
  saddle-node singularity}, Journal of Differential Equations, 248 (2010),
  pp.~2841--2888.

\bibitem{kugler2018early}
{\sc P.~K{\"u}gler, A.~H. Erhardt, and M.~Bulelzai}, {\em Early
  afterdepolarizations in cardiac action potentials as mixed mode oscillations
  due to a folded node singularity}, PLoS One, 13 (2018), p.~e0209498.

\bibitem{letson2017analysis}
{\sc B.~Letson, J.~E. Rubin, and T.~Vo}, {\em Analysis of interacting local
  oscillation mechanisms in three-timescale systems}, SIAM Journal on Applied
  Mathematics, 77 (2017), pp.~1020--1046.

\bibitem{pages2019model}
{\sc N.~Pages, E.~Vera-Sig{\"u}enza, J.~Rugis, V.~Kirk, D.~I. Yule, and
  J.~Sneyd}, {\em A model of ca2+ dynamics in an accurate reconstruction of
  parotid acinar cells}, Bulletin of Mathematical Biology, 81 (2019), p.~1394.

\bibitem{rubin2007giant}
{\sc J.~Rubin and M.~Wechselberger}, {\em Giant squid-hidden canard: the 3{D}
  geometry of the {H}odgkin-{H}uxley model}, Biological Cybernetics, 97 (2007),
  pp.~5--32.

\bibitem{rubin2008selection}
\leavevmode\vrule height 2pt depth -1.6pt width 23pt, {\em The selection of
  mixed-mode oscillations in a {H}odgkin-{H}uxley model with multiple
  timescales}, Chaos: An Interdisciplinary Journal of Nonlinear Science, 18
  (2008), p.~015105.

\bibitem{szmolyan2001canards}
{\sc P.~Szmolyan and M.~Wechselberger}, {\em Canards in $\mathbb{R}^3$},
  Journal of Differential Equations, 177 (2001), pp.~419--453.

\bibitem{szmolyan2004relaxation}
\leavevmode\vrule height 2pt depth -1.6pt width 23pt, {\em Relaxation
  oscillations in $\mathbb{R}^3$}, Journal of Differential Equations, 200
  (2004), pp.~69--104.

\bibitem{takano2021highly}
{\sc T.~Takano, A.~M. Wahl, K.-T. Huang, T.~Narita, J.~Rugis, J.~Sneyd, and
  D.~I. Yule}, {\em Highly localized intracellular ca2+ signals promote optimal
  salivary gland fluid secretion}, Elife, 10 (2021), p.~e66170.

\bibitem{wechselberger2005existence}
{\sc M.~Wechselberger}, {\em Existence and bifurcation of canards in
  $\mathbb{R}^{3}$ in the case of a folded node}, SIAM Journal on Applied
  Dynamical Systems, 4 (2005), pp.~101--139.

\bibitem{wechselberger2020geometric}
\leavevmode\vrule height 2pt depth -1.6pt width 23pt, {\em Geometric singular
  perturbation theory beyond the standard form}, Springer Nature, 2020.

\end{thebibliography}
\end{document}